\documentclass[openany]{book} %%[12pt,a4paper,openany,titlepage]
\usepackage[total={15cm,23cm},top=2cm,left=3cm]{geometry}
\usepackage[latin1,utf8]{inputenc}
\usepackage[spanish]{babel}
\usepackage{amsmath}
\usepackage{amsfonts}
\usepackage{amssymb}
\usepackage{amsthm}
\usepackage{graphicx}
\usepackage{epstopdf}
\usepackage{float}
\usepackage{subfig}
\usepackage{cancel}
\theoremstyle{plain}
\newtheorem{defi}{Definición}[section]
\newtheorem{theorem}{Teorema}[section]
\newtheorem{obs}{Observación}
\newtheorem{lema}{Lema} [chapter]
\theoremstyle{definition}
\newtheorem{ejem}{Ejemplo}[chapter]
\numberwithin{equation}{chapter}
\spanishdecimal{.}

\begin{document}

\begin{titlepage}

\begin{center}
\begin{figure}[htb]
\begin{center}
\includegraphics[width=3.5cm]{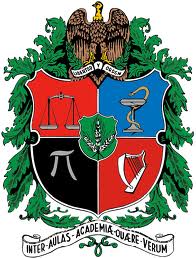}
\end{center}
\end{figure}
\rule{100mm}{0.1mm}\\
\vspace*{0.1in}

\begin{Large}
Universidad Nacional de Colombia \\ 
Facultad de Ciencias \\ 
Departamento de Matemáticas \\ 
\vspace*{0.6in}
\end{Large}

\begin{LARGE}

\textbf{ Solución numérica de ecuaciones diferenciales parciales parabólicas} \\
\textit{\Large Método de elementos finitos aplicado 
a las ecuaciones de Stokes y de Advección-Difusión}

\end{LARGE}

\vspace*{0.7in}

\begin{large}

Trabajo de grado presentado por: \\ 
Jonathan David Galeano Vargas \\
\vspace*{0.7in}
Como requisito para el grado de matemático en la Universidad Nacional de Colombia\\

\vspace*{0.6in}
\rule{80mm}{0.1mm}\\
\vspace*{0.1in}

Revisado por: \\
Juan Galvis \\
\end{large}
\vspace*{0.4in}
Bogotá, Colombia \\
2013-II
\end{center}

\end{titlepage}

\tableofcontents

\part{Introducción al Método de Elementos Finitos }

\chapter{Introducción}

Al estudiar ecuaciones diferenciales parciales, existen distintos métodos para resolverlas numéricamente, pero en este trabajo nos enfocaremos en el ``Método de los elementos finitos (MEF)", el cual, en general, requiere de tres pasos:

\begin{enumerate}

\item Debemos construir una \textit{``formulación débil (o variacional)"} de la ecuación diferencial en estudio, esto implica que, la ecuación debe ser planteada en un espacio de funciones adecuado, que frecuentemente es un espacio de Hilbert, y luego verificar la existencia de la solución de la formulación débil para este problema.

\item Luego de escribir el problema en un espacio de funciones adecuado, debemos introducir espacios de elementos finitos, es decir, debemos aproximar el problema obtenido en el primer paso a un espacio de dimension finita, de aquí, surge un sistema linear, ya que introducimos unas bases para hallar una aproximación de elementos finitos.

\item Por ultimo, debemos resolver el sistema lineal obtenido en el paso 2, el cual se soluciona por un método directo o iterativo.

\end{enumerate}

A continuación presentaremos algunos preliminares, usados mas adelante, para aclarar algunos resultados del trabajo.
 
\section{Preliminares}

En esta sección introduciremos nociones necesarias para el adecuado desarrollo del trabajo, que son indispensables para el estudio del \textit{Método de los Elementos Finitos}. Mas detalles sobre las temáticas de esta sección se pueden consultar en \cite{Adams, Evans}.

\subsection{Espacios de funciones}
\label{EP}

Para obtener una forma precisa de una formulación débil de alguna ecuación diferencial parcial (EDP), o para estudiar el comportamiento del MEF, es necesario introducir algunos espacios de funciones.

Para esto, comenzamos recordando la noción de norma.

\begin{defi}

Dado un espacio vectorial V, una función $\parallel \cdot \parallel : V \rightarrow \mathbb{R}^+ \cup \{0\} $ define una norma, si satisface las siguientes propiedades:

\begin{enumerate}
\item Dado $v \in V$, $\parallel v \parallel = 0$ si y solamente si $v=0.$
\item Para todo $\alpha \in \mathbb{R}$ y $v\in V$ se tiene que,
$$\parallel\alpha v\parallel=|\alpha| \parallel v \parallel. $$
\item Dados cualesquiera $u,v \in V$ $$||u+v|| \leq ||u|| + ||v||.$$
\end{enumerate}
\end{defi}

De aquí, decimos que el par $(V,||\cdot||)$ forma un espacio vectorial normado.

\subsubsection{Espacios por compleción}

En esta sección veremos como definir espacios de funciones por compleción, esto sera necesario para poder abordar el problema débil de una ecuación diferencial parcial.

Presentamos la definición de un espacio vectorial normado y completo.

\begin{defi}

Decimos que un espacio vectorial normado es completo si para cualquier sucesión de Cauchy $\{v_n\}_{n=1}^\infty \subset V$, existe $u \in V$ tal que el $\displaystyle \lim_{n \to \infty} v_n = v.$

\end{defi}

Ahora, el siguiente teorema nos muestra que se puede completar de manera única un espacio vectorial normado.

\begin{theorem}
Sea $(V,||\cdot||)$ un espacio vectorial normado, entonces existe un único espacio vectorial completo $(H,||\cdot||)$ tal que

\begin{itemize}
\item $V\subset H.$
\item Dado cualquier elemento $v \in H$, existe una sucesión $\{v_n\}_{n=1}^\infty \subset V$, tal que 
$$\displaystyle \lim_{n \to \infty} v_n = v.$$
\end{itemize}

En este caso decimos que $V$ es denso en H; o dicho de mejor manera que $H$ es el cierre (o compleción) de $V.$ 
\end{theorem}

Para una demostración de este teorema véase \cite{Lima}, para verla en detalle. Usando esta herramienta, podemos construir los siguientes espacios, denominados \textit{Espacios de funciones tipo Sobolev}. Véase \cite{Evans} para mas detalles. 

Pero antes definiremos los espacios $C(\Omega)$, $C^\infty(\Omega)$ y $C_0^\infty(\Omega)$, que van a ser los espacios que se van a completar por medio de esta herramienta.

\begin{obs}
Definimos $C(\Omega)$, $C^\infty(\Omega)$ y $C_0^\infty(\Omega)$ como,
\begin{itemize}
\item El Espacio $C(\Omega)$ se define como el conjunto de las funciones reales continuas con dominio $\Omega$.
\item El espacio $C^\infty(\Omega)$ es el conjunto de las funciones continuamente diferenciables, es decir, si existen todas sus derivadas en $\Omega$.
\item El espacio $C_0^\infty(\Omega)$ se define como,

$$C_0^\infty(\Omega)=\{v \in C^\infty(\Omega) : v=0 \text{ en } \Gamma \subset \partial \Omega \}.$$
\end{itemize}
\end{obs}

Ahora presentaremos los espacios de funciones que obtenemos por medio de esta herramienta.

\subsubsection*{Espacio $L^2$ }

Sea la norma $L^2$ de una función $v: \Omega \rightarrow \mathbb{R} $ dada por

$$ ||v||_0 = \left(\int_\Omega v(x)^2 dx \right) ^{1/2}. $$

El espacio $L^2(\Omega)$ esta dado por el compleción del espacio $C(\Omega)$ con relación a la norma $L^2$, el espacio $L^2(\Omega)$ es el conjunto de las funciones cuadrado integrable, es decir si $v \in C^\infty (\Omega) $, entonces

\begin{center}
$v \in L^2(\Omega)$ \textit{si y solo si} $\displaystyle \int_\Omega v^2 < \infty. $
\end{center}

Es decir, $v$ es cuadrado integrable si la integral de la función al cuadrado es finita.

\subsubsection*{Espacios $H^1$, $H_0^1$}

Definimos la norma $H^1$ como,

$$ ||v||_1 = \left(\int_\Omega v(x)^2 + |\nabla v(x)|^2 dx \right) ^{1/2}. $$

Los espacios $H^1(\Omega)$ y $H_0^1(\Omega)$, son la compleción de los espacios $C^\infty(\Omega)$ y $C_0^\infty(\Omega)$ con relación a la norma $||\cdot||_1$.
El espacio $H^1(\Omega)$ es un espacio de \textit{Hilbert}, apropiado para la mayoría de problemas considerados en este trabajo. Note que,

\begin{center}
$H^1(\Omega) := \{v \in L^2(\Omega) \mid v' \in L^2(\Omega) \}.$
\end{center}

También usaremos el espacio de funciones $H_{0}^{1}(\Omega)$, definido por,

\begin{center}
$H_{0}^{1}(\Omega) := \{v \in H^1(\Omega) \mid  v=0  \mbox{ en } \Gamma \subset \partial \Omega \},$
\end{center}
tal que, $H_{0}^{1}(\Omega) \subset H^1(\Omega).$

\subsubsection*{Espacio $H^2$}

Finalmente definimos la norma $H^2(\Omega)$ de la siguiente manera,

$$ ||v||_2 = \left(\int_\Omega v(x)^2 + |\nabla v(x)|^2 + |\Delta v(x)|^2 dx \right) ^{1/2}. $$

Donde el espacio $H^2(\Omega)$ es definido como la compleción del espacio $C^\infty(\Omega)$ con respecto a la norma $||\cdot||_2$. Decimos que $H_2$ es definido como,
\begin{center}
$H^2(\Omega) := \{v \in L^2(\Omega) \mid v',v'' \in L^2(\Omega) \}.$
\end{center}
Podemos ver que, $$H^2(\Omega) \subset H^1(\Omega) \subset L^2(\Omega).$$

\begin{ejem}
Para la ecuación de Laplace, el espacio $H^1(\Omega)$ es usado para construir su formulación débil. 
\end{ejem}

\subsubsection{El espacio de funciones lineales por partes}

En esta sección introduciremos $V \subset \mathbb{R}^d$ $d=1,2$, el espacio de funciones lineales por partes, el cual es uno de los ejemplos mas simples de espacios de elementos finitos. Primero, debemos conocer la noción de \textit{Espacio de funciones de dimensión finita}.

\begin{defi}[Espacio de dimensión finita]

Un espacio de funciones $V$ tiene dimensión finita $n$, si existe un conjunto de funciones $\phi_i:\Omega \rightarrow \mathbb{R}$, $i \in \{1,\ldots,n\},$ tal que cualquier función $\varphi \in V$ puede ser escrita como una única combinación lineal de las funciones $\phi_i$. Esto es, existen n constantes $\alpha_i$, tal que, $\varphi = \sum_{i=1}^n \alpha_i\phi_i. $
\end{defi}

Primero consideraremos el caso unidimensional $(d=1)$, es decir cuando $\Omega \subset \mathbb{R}.$

\subsubsection*{Caso unidimensional}

Como $\Omega \subset \mathbb{R}$, podemos suponer, sin perdida de generalidad que $\Omega = (a,b)$. Entenderemos como una partición de $\Omega$, una división de este dominio en subintervalos de la siguiente forma $[a,b]=\cup_{i=1}^N[z_{i-1},z_i], \text{ donde } z_0=a, z_N=b \text{ y } z_i<z_{i+1}$. Los valores $z_i$ son llamados vértices de la partición. 

Vamos a definir el espacio de las funciones lineales por partes sobre la partición introducida como,
$$V=\{f \in C(\Omega); f|_{(z_i,z_{i+1})} \text{ es lineal}\},$$ 
lo que quiere decir que, $f|_{(z_i,z_{i+1})} = a_ix + b_i $, siendo $a_i$, y $b_i$ constantes apropiadas.

Podemos ver que el conjunto de funciones $\phi_i \in V$, $i \in \{0,1,\ldots,N\}$ donde, 
$$\phi_i(z_j)=\left\{ \begin{matrix}
1 & \text{Si } i=j \\
0 & \text{Si } i\neq j
\end{matrix} \right. $$
forman una base para el espacio $V$. En otras palabras, para cualquier función $v \in V$ puede ser escrita como,

$$v(x)=c_0\phi_0(x) + c_1\phi_1(x)+\ldots+c_N\phi_N(x)=\displaystyle\sum_{i=0}^N c_i\phi_i(x),$$
donde $c_i=v(x_i)$ con $x_i=z_i$ vértices de la partición.

\subsubsection*{Caso bidimensional }

Para el caso en dos dimensiones, es decir, $\Omega \subset \mathbb{R}^2$, supongamos que $\Omega$  es un dominio poligonal, y así simplificar la introducción del espacio de funciones lineales por partes en dos dimensiones. Comenzaremos introduciendo la noción de triangulación.

\begin{defi}
Una triangulación de un dominio poligonal $\Omega \subset \mathbb{R}^2$ es una subdivisión de $\Omega$ en un conjunto finito de triángulos que satisfacen la siguiente propiedad: Cualquier vértice de la partición no pertenece al interior de cualquier arista de la partición.
\label{triang}
\end{defi}

Sea $\mathcal{T}$ una triangulación de $\Omega$. Definimos el espacio de funciones lineales por partes sobre $\mathcal{T}$,  como el espacio formado por las funciones $v$ tales que $v$ restringida a cada triángulo de $\mathcal{T}$ es lineal en las variables $x$ y  $y$. Esto es, para todo triangulo $T \in \mathcal{T}$, tenemos que
$$v|_T(x,y)=ax+by+c,$$ donde $a,b,c$ son constantes apropiadas, y dependen del triángulo $T$. Estas constantes corresponden a los tres grados de libertad que tenemos para cada función lineal por partes restringida al triángulo $T$. Vemos también que la base para $\mathbb{P}^1(\mathcal{T}^h)$ es $\{1,x,y\}$ y como tiene tres grados de libertad dim $\mathbb{P}^1(\mathcal{T}^h)=3$.

El espacio de funciones lineales por partes sobre $\mathcal{T}$ lo notaremos como $V$, y podemos ver que, una base para el espacio $V$, es dada por el conjunto de funciones $\phi_i \in V, \, i \in \{0,1,\ldots,N\},$ donde $N$ representa el numero de vértices de $\mathcal{T}$ y cada función base es de la forma,
$$\phi_i(z_j)= \begin{cases}
1 & \text{Si } i=j \\
0 & \text{Si } i \neq j \\
\text{Funcion lineal} & \text{En otro caso.}
\end{cases}$$ 

Vamos a definir una característica que tienen las triangulaciones de estos espacios de funciones lineales por partes.

\begin{defi}
Sea $\mathcal{T}^h$ una familia de triangulaciones de $\Omega \subset \mathbb{R}^2$, $0<h \leq 1.$
\begin{itemize}
\item Decimos que una familia de triangulaciones $\mathcal{T}^h$ es de \textit{aspecto regular} si existe una constante $C>0$ independiente de $h$ tal que $\rho(K)\leq C$ para todo elemento
$K \in \mathcal{T}^h$, donde $\rho(K)=diam(K)/r_k$, donde $r_k$ es el radio del mayor circulo contenido en $K$ y $diam(K)$ es el radio del menor circulo que contiene a $K.$
\item Decimos que una la familia de triangulaciones $\mathcal{T}^h$ es \textit{cuasi uniforme} si existe una constante $C$ independiente de $h$ tal que $diam(K) \geq Ch$ para todo elemento $K \in \mathcal{T}^h$. O visto de otra manera, si existen constantes $c_1$ y $c_2$ tales que,
$$\max_{T \in \mathcal{T}^h}\{\text{diametro de } B_T \leq c_1h \quad \text{y} \quad \min_{T \in \mathcal{T}^h}\{\text{diametro de } b_T \geq c_2h,$$
donde $B_T$ es el menor circulo que contiene a $T$, y $b_T$ es el mayor circulo contenido en $T$.
\end{itemize}
\end{defi} 

Dada una familia de triangulaciones $\mathcal{T}^h$ de $\Omega$ vamos a denotar $\mathbb{P}^1(\mathcal{T}^h)$ el espacio de funciones lineales por partes, definido en esta sección.

\begin{equation*}
\mathbb{P}^1(\mathcal{T}^h)=\left\{
\begin{aligned}
 v \in C(\Omega) \quad : & \ v|_K \text{ es un polinomio en dos variables de grado total 1  } \\
&\text{ para todo elemento de la triangulación } \mathcal{T}^h. 
\end{aligned}
\right\}
\label{P1}
\end{equation*}

\subsubsection{El espacio de funciones cuadráticas por partes}
\label{EP2} 
Dada una familia de triangulaciones $\mathcal{T}^h$ de $\Omega$ denotaremos $\mathbb{P}^2_0(\mathcal{T}^h)$ el espacio de funciones cuadráticas por partes, el cual definiremos a continuación. Para mas detalles, véase \cite{Cla}.

Para  $\Omega \subset \mathbb{R}^2$, supongamos que $\Omega$  es un dominio poligonal, para así simplificar la introducción del espacio de funciones cuadráticas por partes usando la noción de triangulación definida en \eqref{triang}.

Sea $\mathcal{T}$ una triangulación de $\Omega$. Definimos el espacio de funciones cuadráticas por partes sobre $\mathcal{T}$,  como el espacio formado por las funciones $v$ tales que $v$ restringida a cada triángulo de $\mathcal{T}$ es de grado total$\leq 2$ en las variables $x$ y  $y$, es decir que la suma de los grados de cada variable no sea mayor que 2. Esto es, para todo triangulo $T \in \mathcal{T}$, tenemos que
$$v|_T(x,y)=a+bx+cy+dx^2+exy+fy^2,$$ donde $a,b,c,d,e,f$ son constantes apropiadas, y dependen del triángulo $T$. Estas constantes corresponden a los seis grados de libertad que tenemos para cada función cuadrática por partes restringida al triángulo $T$. Vemos también que la base para $\mathbb{P}^2(\mathcal{T}^h)$ es $\{1,x,y,x^2,xy,y^2\}$ y como tiene seis grados de libertad, entonces dim $\mathbb{P}^2(\mathcal{T}^h)=6$.

\begin{figure}[H]
  \centering
    \includegraphics[width=0.61\textwidth]{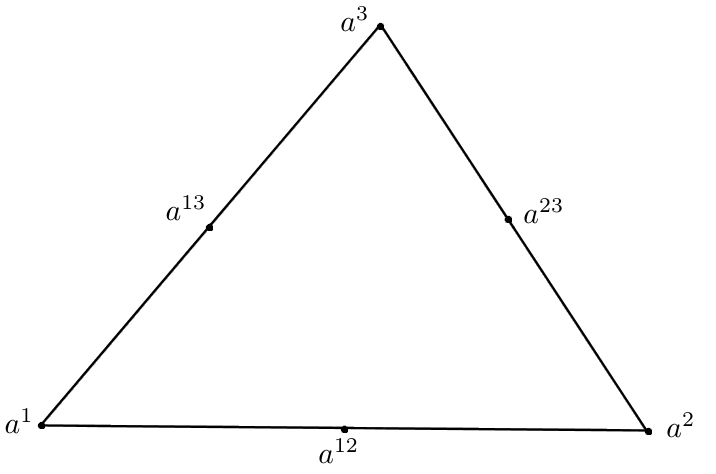}
  \caption{Distribución de los nodos para el espacio $\mathbb{P}^2$}
  \label{TP21}
\end{figure}

Como tiene 6 grados de libertad, vamos a determinar los nodos en cada elemento $K$ (triángulo) correspondientes a los seis grados de libertad, que serán los tres vértices $a^i$, $i=1,2,3.$ y los tres puntos medios de los lados de $K$ $a^{ij}$, $i<j$, $i,j=1,2,3.$ véase la figura \ref{TP21}.

Ahora, tenemos el siguiente teorema, demostrado en \cite{Cla}.
\begin{theorem}
Una función $\varphi \in \mathbb{P}^2(\mathcal{T}^h$ es únicamente determinada por los siguientes grados de libertad:
\begin{equation}
\begin{aligned}
\varphi(a^i), \quad \quad & i=1,2,3, \\
\varphi(a^{ij}), \quad \quad & i<j, \ i,j=1,2,3.
\end{aligned}
\end{equation}
\end{theorem}

Podemos ver que, una función base para el espacio de funciones cuadráticas por partes, tiene la representación
 
$$\varphi = \sum_{i=1}^3\varphi(a^i)\lambda_i(2\lambda_i-1) + \sum_{\substack{i,j=1 \\ i<j}}^3\varphi(a^{ij})4\lambda_i\lambda_j.$$ 
donde $\lambda_i$ es función base de $\mathbb{P}^1(\mathcal{T}^h)$. Por lo que las funciones base para $\mathbb{P}^2(\mathcal{T}^h)$ corresponde al particular grado de libertad, por ejemplo, el valor en el vértice $a^i$, es una función $\Psi \in \mathbb{P}^2(\mathcal{T}^h)$ tal que $\Psi(a^i)=1$ y $\Psi$ se anula en el resto de nodos $a^i,a^{ij}.$

Entonces $\mathbb{P}^2(\mathcal{T}^h)$ se define como,
\begin{equation*}
\mathbb{P}^2(\mathcal{T}^h)=\left\{
\begin{aligned}
 v \in C(\Omega) \quad : & \ v|_K \text{ es un polinomio en dos variables de grado total 2  } \\
&\text{ para todo elemento de la triangulación } \mathcal{T}^h.
\end{aligned}
\right\}
\label{P2}
\end{equation*}

\section{Condiciones de Frontera}

Existen distintas condiciones de frontera, que se usan como su nombre lo indica, para restringir el valor de una EDP en su frontera, de acuerdo a la situación del problema. Estas condiciones se pueden introducir independientes una de la otra, o se pueden combinar, pero en este trabajo solo trabajaremos con la condición de Dirichlet, y la condición de Neumann, las cuales presentaremos a continuación.

\subsection{Condición de Dirichlet}

Esta condición indica que la solución, restringida a la frontera, es requerida con un valor determinado, la enunciaremos para los casos en una y dos dimensiones.

\subsubsection*{Una dimensión}

Para una dimension, sabemos que $\Omega = (a,b) \subset \mathbb{R}$ y que su frontera son los puntos $a$ y $b$.  

\begin{equation*}
\begin{aligned}[b]
&\textit{Encontrar } u : (a,b) \to \mathbb{R} \textit{ tal que:} \\
&\left\{
\begin{aligned}
&\text{EDP en estudio}  \\
&u(x) = g(x),
\end{aligned}
\quad
\begin{aligned}
&\textnormal{para } x\in (a,b) \\
&\textnormal{para } x=a, x=b.
\end{aligned} \right.
\end{aligned}
\end{equation*} 

\subsubsection*{Dos dimensiones}

Para dos dimensiones tenemos que $x=(x_1,x_2)$, $\Omega \subset \mathbb{R}^2$, y frontera $\partial \Omega$ entonces,

\begin{equation*}
\begin{aligned}[b]
&\textit{Encontrar } u : \Omega \to \mathbb{R} \textit{ tal que:} \\
&\left\{
\begin{aligned}
&\text{EDP en estudio}  \\
&u(x) = g(x),
\end{aligned}
\quad
\begin{aligned}
&\textnormal{para } x\in \Omega \\
&\textnormal{para } x \in \Gamma \subset \partial \Omega.
\end{aligned} \right.
\end{aligned}
\end{equation*}

\subsection{Condición de Neumann}

Esta condición indica que la solución, restringida a la frontera, es requerida con un valor determinado, la enunciaremos para los casos en una y dos dimensiones.

\subsubsection*{Una dimension}

Para una dimension, sabemos que $\Omega = (a,b) \subset \mathbb{R}$ y que su frontera son los puntos $a$ y $b$. 

\begin{equation*}
\begin{aligned} [t]
&\textit{Encontrar } u : (a,b) \to \mathbb{R} \textit{ tal que:} \\
&\left\{
\begin{aligned} 
&\text{EDP en estudio}  \\
&\kappa(x) u'(x)= g(x), \\
&g(a)-g(b)={\textstyle\int_a^b}f \\
&{\textstyle\int_{a}^{b}} u=0 
\end{aligned}
\quad
\begin{aligned}
&\textnormal{para } x\in (a,b) \\
&\textnormal{para } x=a, x=b \\ 
&\textnormal{(Condición de Compatibilidad) } \\ 
&\textnormal{(Condición de solubilidad) } \\ 
\end{aligned}\right.
\end{aligned} 
\end{equation*} 

\subsubsection*{Dos dimensiones}
Para dos dimensiones tenemos que $x=(x_1,x_2)$, $\Omega \subset \mathbb{R}^2$, y frontera $\partial \Omega$ entonces,
\begin{equation*}
\begin{aligned} [t]
&\textit{Encontrar } u : \Omega \to \mathbb{R} \textit{ tal que:} \\
&\left\{
\begin{aligned} 
&\text{EDP en estudio}  \\
&-\nabla u(x)\cdot \vec{\eta}(x)= h(x), \\
&{\textstyle\int_\Omega}u(x) = 0 \\
&{\textstyle\int_{\partial\Omega}}h(x)dS={\textstyle\int_{\Omega}}f(x)dx
\end{aligned}
\quad
\begin{aligned}
&\textnormal{para } x\in \Omega \\
&\textnormal{para } x \in \partial \Omega \\ 
&\textnormal{(Condición de solubilidad)} \\ 
&\textnormal{(Condición de Compatibilidad) } \\ 
\end{aligned}\right.
\end{aligned} 
\end{equation*} 

\section{Algunas EDP's}

Algunas EDP's que se usaran para introducir el MEF serán enunciadas a continuación.

\subsection{EDP básica}

Consideremos la siguiente EDP, en esta sección esta enunciada con la condición de Dirichlet, pero esta condición cambia de acuerdo al problema. La enunciaremos en los casos de una y dos dimensiones.

\subsubsection*{Una dimensión}

Para una dimensión, tomamos $\Omega = (a,b) \subset \mathbb{R}$, y consideramos la EDP básica como,

\begin{equation}
(\text{Formulación Fuerte})
\begin{aligned}[b]
&\textit{Encontrar } u : (a,b) \to \mathbb{R} \textit{ tal que:} \\
&\left\{
\begin{aligned}
-(\kappa(x)u'(x))' &= f(x),  \\
u(x) &= g(x),
\end{aligned}
\quad
\begin{aligned}
&\textnormal{para } x\in (a,b), \\
&\textnormal{para } x=a, x=b,
\end{aligned} \right.
\label{EE1D}
\end{aligned}
\end{equation}

donde $0<\kappa_{min}\leq \kappa(x) \leq \kappa_{max}$ para todo $x \in (a,b)$.

\subsubsection*{Dos dimensiones}
Para dos dimensiones tenemos que $x=(x_1,x_2)$, $\Omega \subset \mathbb{R}^2$, y frontera $\partial \Omega$ entonces,
\begin{equation}
(\text{Formulación Fuerte})
\begin{aligned}[b]
&\textit{Encontrar } u : \Omega \to \mathbb{R} \textit{ tal que:} \\
&\left\{
\begin{aligned}
-\text{div }(\kappa(x)\nabla u(x)) &= f(x),  \\
u(x) &= g(x),
\end{aligned}
\quad
\begin{aligned}
&\textnormal{para } x\in \Omega, \\
&\textnormal{para } x\in \partial\Omega,
\end{aligned} \right.
\label{EE2D}
\end{aligned}
\end{equation}

donde $0<\kappa_{min}\leq \kappa(x) \leq \kappa_{max}$ para todo $x \in \Omega$.
\subsection{EDP de Laplace}

Esta EDP la usaremos, aludiendo a su fácil estudio, y comodidad para trabajar, acá definiremos su denominada \textit{Formulación Fuerte}, para el caso de una y dos dimensiones, con condición de Dirichlet, la cual puede variar de acuerdo al problema.

\subsubsection*{Una dimensión}

Para una dimensión, tomamos $\Omega = (a,b) \subset \mathbb{R}$, y consideramos la EDP de Laplace como,

\begin{equation}
(\text{Formulación Fuerte})
\begin{aligned}[b]
&\textit{Encontrar } u : (a,b) \to \mathbb{R} \textit{ tal que:} \\
&\left\{
\begin{aligned}
-u''(x) &= f(x),  \\
u(x) &= g(x),
\end{aligned}
\quad
\begin{aligned}
&\textnormal{para } x\in (a,b), \\
&\textnormal{para } x=a, x=b.
\end{aligned} \right.
\end{aligned}
\label{EL1D}
\end{equation}

\subsubsection*{Dos dimensiones}
Para dos dimensiones tenemos que $x=(x_1,x_2)$, $\Omega \subset \mathbb{R}^2$, y frontera $\partial \Omega$ entonces,

\begin{equation}
(\text{Formulación Fuerte})
\begin{aligned}[b]
&\textit{Encontrar } u : \Omega \to \mathbb{R} \textit{ tal que:} \\
&\left\{
\begin{aligned}
-\Delta u(x) &= f(x),  \\
u(x) &= g(x),
\end{aligned}
\quad
\begin{aligned}
&\textnormal{para } x\in \Omega, \\
&\textnormal{para } x \in \partial \Omega.
\end{aligned} \right.
\end{aligned}
\label{EL2D}
\end{equation}

\section{Sobre la existencia de soluciones débiles}
\subsection{Teorema de Lax-Milgram}
\label{LM}

El teorema de Lax-Milgram nos garantiza la existencia de soluciones débiles para las EDP's elípticas. Pero para esto debemos conocer la noción de forma bilineal y funcional lineal.

\begin{defi}
Sea V un espacio de funciones dotado de una norma $||\cdot||.$ Decimos que una aplicación $f: V \to \mathbb{R}$ es un funcional lineal limitado si (la aplicación del funcional $f$ a una función $v \in V$ se denota como $\langle f,v \rangle$): 
\begin{itemize}
\item Para toda $v_1,v_2 \in V$ y $c_1,c_2 \in \mathbb{R}, \, \langle f,c_1v_1 + c_2v_2 \rangle = c_1\langle f,v_1 \rangle + c_2 \langle f,c_2v_2 \rangle.$
\item Existe una constante $c \in \mathbb{R}$ tal que $|\langle f,v \rangle| \leq c||v||.$
\end{itemize}
\noindent También que una aplicación $\mathcal{A}: V \times V \to \mathbb{R}$ es una forma bilineal si,
\begin{equation*}
\begin{aligned}
\mathcal{A}(c_1v_1 + c_2v_2, d_1w_1 + d_2w_2) = & c_1d_1\mathcal{A}(v_1,w_1) + c_1d_2\mathcal{A}
(v_1,w_2)c_2d_1\mathcal{A}(v_2,w_1)  \\ & + c_2d_2\mathcal{A}(v_2,w_2)
\end{aligned}
\end{equation*}
para todo $c_i,d_i \in \mathbb{R}$ y  $w_i,v_i \in V, \, i \in \{1,2\}.$
\end{defi}

Y ahora introducimos el teorema de Lax-Milgram, usado para garantizar la existencia de una solución de una formulación débil de una EDP elíptica

\begin{theorem}[Teorema de Lax-Milgram]

Sea $V$ un espacio de funciones completo dotado de una norma $||\cdot||$, y $B:V\times V \rightarrow \mathbb{R}$
es una forma bilineal, que satisface
\begin{equation}
|B(u,v)|\leq \alpha ||u||||v||, \quad \text{para toda } u,v \in V 
\end{equation}
y es elíptica, esto es 
\begin{equation}
\beta||u||^2\leq B(u,u) \quad \text{para toda } u \in V
\end{equation}
para constantes $\alpha,\beta>0$ apropiadas. Entonces para cada funcional lineal limitado $f:V\rightarrow\mathbb{R}$ existe una única función $u\in V$, tal que
$$B(u,v)=\langle f,v \rangle \quad \text{para toda } v\in V.$$ 
\end{theorem}

Los detalles de esta demostración están en \cite{Evans}, al aplicar este teorema, podemos garantizar la existencia de la solución de una EDP elíptica, pero antes debemos elegir el espacio de funciones $V$ que va a ser utilizado, el cual va a depender de la condición de contorno impuesta al problema. 

\subsection{Condición de Babuska-Brezzi}
\label{Brezzi}
Este teorema es mas general, que el teorema de Lax-Milgram, y se usa en el caso en que alguna de las formas bilineales no sea simétrica, o cuando es una \textit{formulación punto de silla} como es el caso de la ecuación de Stokes, y la ecuación de advección-difusión.

La condición dice lo siguiente,

\begin{theorem}[Condición de Babuska-Brezzi]

Sean H y Q espacios de Hilbert, $\mathcal{A} : H \times  H \to \mathbb{R}$, $ \mathcal{B} : H \times Q \to R$ formas bilineales acotadas y $V := \{v \in H : \mathcal{B}(v,\mu) = 0 \text{ para todo } \mu \in Q\}$.
Suponga que:
\begin{itemize}
\item $\mathcal{A}$ es V-elíptica: $\mathcal{A}(v,v)\geq \alpha ||v||_H^2.$
\item $\mathcal{B}$ satisface la condición de Babuska-Brezzi:
$$\sup_{v\in H, v\neq 0} \dfrac{\mathcal{B}(v,\mu)}{||v||_H} \geq ||\mu||_Q \quad \quad \text{para todo } \mu \in Q.$$
\end{itemize} 
Entonces, para todo $(F,G) \in H' \times Q'$ existe un único $(u;\lambda) \in H \times Q$ tal que,
$$\mathcal{A}(u,v) + \mathcal{B}(v,\mu) = F(v) \quad \quad \text{para todo } \mu \in H.$$
$$\mathcal{B}(u,\mu) = G(\mu) \quad \quad \text{para todo } \mu \in Q.$$
Además, existe $ C > 0$ independiente de $(u,\lambda)$ tal que,
$$||(u,\lambda)|| \leq  C \{ ||F||_{H'} + ||G||_{Q'}\}.$$
\end{theorem}
Este teorema permite encontrar la solución de formulaciones de punto de silla, se encuentra demostrado en \cite{Babuska,Brezzi}.

\chapter{Método de los elementos finitos en una dimensión}

Veamos el MEF en una ecuación diferencial en una dimensión (tomaremos, sin perdida de generalidad, el intervalo (a,b)), encontraremos los espacios adecuados para resolver el problema, sus formulaciones, y estudiaremos algunos ejemplos.

\section{Formulación Débil}

Para deducir la formulación débil de una ecuación diferencial en el intervalo (a,b) tenemos que,

\begin{enumerate}
\item Suponer que existe $u$ solución de nuestra ecuación diferencial, multiplicar los dos lados de la ecuación por una \textit{función de prueba} $v \in C_0^\infty(a,b)$ e integrar a ambos lados de la igualdad.

\item Usar la fórmula de integración por partes y las condiciones iniciales (de contorno) para obtener expresiones que necesiten solamente de derivadas del menor orden posible.

\item Cambiar $v \in C_0^\infty(a,b)$ por $v \in V$ donde el espacio de funciones de prueba $V$ sea el mas grande posible, también escoger $U$ (espacio de funciones de forma), tal que la solución $u$ este naturalmente en $U$.   
\end{enumerate}

Los espacios de funciones $U$ y $V$ mas adecuados para estas condiciones son los \textit{Espacios de funciones tipo Sobolev (Espacios de Hilbert)} definidos en $(a,b)$. Véase \cite{Adams,Evans}. Así pues, los espacios $U$ y $V$ adecuados deben cumplir:

\begin{itemize}

\item La elección de estos espacios $U$ y $V$ es independiente una de otra, aunque por comodidad, se puede tomar $U=V$, y su ventaja es que se obtienen sistemas lineales simétricos, y de aquí podemos llegar, a calcular mas fácil su solución. 

\item Las integrales que se obtuvieron deben estar definidas en estos espacios.

\item Para toda $u \in U$ se tiene que poder imponer la condición de contorno en estudio.

\end{itemize}

Para la elección de los espacios de funciones, podemos hacer $\Omega = (a,b)\subset \mathbb{R}$ y tomamos los espacios de la Sección \ref{EP}. 

Por ejemplo, los espacios de Hilbert, $H^1(a,b)$, $H_0^1(a,b)$, $H^2(a,b).$

Al obtener dicha \textit{formulación débil}, podemos garantizar la existencia de la solución de esta formulación, que se conoce como \textit{solución débil} de la ecuación original.
La solución del problema original se denomina \textit{solución fuerte}, y en este caso se cumplen las igualdades en el problema para todo $x \in (a,b)$, y por lo tanto se pueden calcular sus derivadas.

\begin{obs}[Véase \cite{Cla}]
Toda solución fuerte de una ecuación es también solución débil. Si los coeficientes de la ecuación son regulares y una solución débil es regular, entonces esta es solución fuerte. Regular, quiere decir que la solución tiene tantas derivadas continuas como el orden de la ecuación diferencial.
\end{obs}

Para probar la existencia de \textit{soluciones débiles} se usan resultados como el Teorema de Lax-Milgram, enunciado en la sección \ref{LM} y demostrado en \cite{Evans}.

\begin{ejem}

Vamos a construir la formulación débil, para la EDP definida en \eqref{EE1D}, tal como,

\begin{equation}
(\text{Formulación Fuerte})
\begin{aligned}[b]
&\textit{Encontrar } u : (a,b) \to \mathbb{R} \textit{ tal que:} \\
&\left\{
\begin{aligned}
-(\kappa(x)u'(x)) &= f(x),  \\
u(x) &= g(x),
\end{aligned}
\quad
\begin{aligned}
&\textnormal{para } x\in (a,b), \\
&\textnormal{para } x=a, x=b,
\end{aligned} \right.
\end{aligned}
\label{EEFF}
\end{equation}

donde $0<\kappa_{min}\leq \kappa(x) \leq \kappa_{max}$ para todo $x \in (a,b).$

Para construir la \textit{formulación débil} de \eqref{EEFF}, primero multiplicamos los dos lados de la igualdad en \eqref{EEFF} por la función de prueba $v \in C_0^\infty(a,b)$ fija pero arbitraria, luego integramos a ambos lados y obtenemos

\begin{equation}
\begin{aligned}[b]
&\textit{Encontrar } u : (a,b) \to \mathbb{R} \textit{ tal que:} \\
&\left\{
\begin{aligned}
-\int_a^b (\kappa(x)u'(x))'v(x)dx = \int_a^b f(x)v(x) dx, \quad
&\textnormal{para toda } v \in C_0^\infty(a,b) \\
u(x) = g(x), \quad
&\textnormal{para } x=a, x=b.
\end{aligned} \right.
\end{aligned}
\label{EEFF1}
\end{equation}

Luego, al usar la fórmula de integración por partes, y como $v \in C_0^\infty(a,b)$, entonces $v(a)=v(b)=0$, por lo que

\begin{align*}
-\int_a^b (\kappa(x)u'(x))v(x)dx &= \int_a^b \kappa(x)u'(x)v'(x)dx -[\kappa(b)u'(b)v(b)- \kappa(a)u'(a)v(a)] \\
&= \int_a^b \kappa(x)u'(x)v'(x)dx -[u'(b)0- u'(a)0] \\
&= \int_a^b \kappa(x)u'(x)v'(x)dx.
\end{align*}

Entonces, \eqref{EEFF1} se escribe como

\begin{equation}
\begin{aligned}[b]
&\textit{Encontrar } u : (a,b) \to \mathbb{R} \textit{ tal que:} \\
&\left\{
\begin{aligned}
\int_a^b u'(x)v'(x)dx = \int_a^b f(x)v(x) dx, \quad
&\textnormal{para toda } v \in C_0^\infty(a,b)  \\
u(x) = g(x), \quad
&\textnormal{para } x=a, x=b. 
\end{aligned} \right.
\end{aligned} 
\label{EEFD1}
\end{equation}
Para que la ecuación \eqref{EEFD1}, sea la formulación débil de \eqref{EEFF}, falta escoger los espacios de funciones $U=H^1(a,b)$ y $V=H_0^1(a,b)$ adecuados para la solución, donde $u \in U$, $v \in V$ y $C_0^\infty(a,b) \subset V$. Introducimos la notación,

\begin{equation}
\mathcal{A}(u,v)= \int_a^b \kappa(x)u'(x)v'(x)dx \quad \text{para toda } v\in V \text{ y } u\in U,
\label{fb}
\end{equation}
donde $\mathcal{A}:U\times V \to \mathbb{R}$ es una forma bilineal, y

\begin{equation}
\mathcal{F}(v)= \int_a^b f(x)v(x)dx \quad \text{para toda } v\in V,
\end{equation}
donde $\mathcal{F}: V \to \mathbb{R}$ es un funcional lineal. Entonces decimos que la formulación débil de \eqref{EEFF} es la ecuación \eqref{EEFD1}, pero puesta en los espacios U y V adecuados. La formulación débil puede ser escrita como,

\begin{equation}
(\text{Formulación Débil})
\begin{aligned}[b]
&\textit{Encontrar } u \in U \textit{ tal que:} \\
&\left\{
\begin{aligned}
\mathcal{A}(u,v) &= \mathcal{F}(v),  \\
u(x) &= g(x),
\end{aligned}
\quad
\begin{aligned}
&\textnormal{para toda } v \in V \\
&\textnormal{para } x=a, x=b.
\end{aligned} \right.
\end{aligned}
\label{EEFD}
\end{equation}
\end{ejem}

\begin{obs}
La formulación débil \eqref{EEFD} se puede reducir al caso $g(x)=0$ en $x=a$ y $x=b$. Es decir si $u_g$ es una función continuamente diferenciable tal que $u_g(a)=g(a)$ y $u_g(b)=g(b)$, entonces $u = u^* + u_g$ donde $u^* \in V $ y se resuelve el problema

\begin{equation}
\begin{aligned}
&\textit{Encontrar } u^* \in V \textit{ tal que:} \\
&\left\{ \mathcal{A}(u^*,v) = \mathcal{F}(v)- \mathcal{A}(u_g,v)\quad \textnormal{para toda } v \in V.
\right.
\end{aligned}
\label{EEFD'}
\end{equation}
\end{obs}

\section{Formulación de Galerkin}

Asumiendo que las formulaciones débiles que trabajaremos son de la forma 

\begin{equation}
\begin{aligned}
&\textit{Encontrar } u^* \in V \textit{ tal que:} \\
&\left\{ \mathcal{A}(u^*,v) = \mathcal{F}(v) \quad \textnormal{para toda } v \in V,
\right.
\end{aligned}
\label{FG'}
\end{equation}
donde $V$ es el espacio de dimensión infinita $V=H_0^1(a,b)$. Ahora, debemos cambiar el espacio de dimension infinita de \eqref{FG'} por espacios de dimensión finita (espacios de elementos finitos).

\bigskip

Considere $V^h$, $h>0$, subespacios de dimension finita tal que $V^h \subset V$. El parámetro $h$ indica el tamaño del espacio, a menor valor $h$, mayor es la dimensión de $V^h$. Lo que se necesita es que $V^h \xrightarrow [h \to 0]\, V$ en algún sentido. Con este $V_h \subset V$ de dimensión finita, la formulación de Galerkin es

\begin{equation}
\begin{aligned}
&\textit{Encontrar } u^*_h \in V^h \textit{ tal que:} \\
&\left\{ \mathcal{A}(u^*_h,v_h) = \mathcal{F}(v_h) \quad \textnormal{para toda } v_h \in V_h.
\right.
\end{aligned}
\label{FG}
\end{equation}
Como el espacio $V_h$ es de dimensión finita, podemos denotar

$$\{\phi_1,\phi_2, \ldots , \phi_{N_h} \} \text{ donde } N_h \text{ es la dimensión de } V_h,$$
como una base de $V_h$, y de aquí, escribir $u_h^*$ como combinación lineal de esta base,

$$u_h^*=x_1\phi_1 + x_2\phi_2 + \cdots + x_{N_h}=\displaystyle \sum_{i=1}^{N_h}x_i\phi_i.$$
Y para toda $v_h \in V^h$ se tiene que 

$$\mathcal{A}(u_h^*,v_h)= \mathcal{A}\left(\displaystyle \sum_{i=1}^{N_h}x_i\phi_i , v_h\right) = \sum_{i=1}^{N_h}x_i\mathcal{A}(\phi_i ,v_h),$$
ya que $\mathcal{A}$ es una forma bilineal. Ahora para verificar \eqref{FG}, basta con hacerlo para las \textit{funciones teste} $v_h=\phi_i,i=1, \ldots, N_h$, funciones de la base de $V_h$, entonces

\begin{equation}
\begin{aligned}
&\textit{Encontrar } u^*_h = \displaystyle \sum_{i=1}^{N_h}x_i\phi_i  \textit{ tal que:} \\
&\left\{ \sum_{i=1}^{N_h}x_i\mathcal{A}(\phi_i ,\phi_j) = \mathcal{F}(\phi_j), \quad j=1,\ldots, N_h.
\right.
\end{aligned}
\label{FG1}
\end{equation}

La forma bilineal $\mathcal{A}$ tiene una representación matricial, tal que, la matriz $A=[a_{ij}]$ de dimensión $N_h \times N_h$ es de la forma

\begin{equation}
a_{ij}=\mathcal{A}(\phi_i,\phi_j), \quad i,j=1,\ldots,N_h.
\label{A}
\end{equation} 
También definiremos $b=b_j$ como el vector

\begin{equation}
b=\begin{bmatrix}
b_1 \\
b_2 \\
\vdots \\
b_{N_h}
\end{bmatrix} \in \mathbb{R}^{N_h}, 
\quad b_j=\mathcal{F}(\phi_j), \quad j=1,\ldots,n,
\label{b}
\end{equation} 
y el vector

$$\boldsymbol{u_h} = \begin{bmatrix}
x_1 \\
x_2 \\
\vdots \\
x_{N_h}
\end{bmatrix} \in \mathbb{R}^{N_h}. $$
Con esto, podemos escribir una formulación matricial equivalente a \eqref{FG} y \eqref{FG1} que es

\begin{equation}
\begin{aligned}
&\textit{Encontrar } u^*_h = \displaystyle \sum_{i=1}^{N_h}x_i\phi_i  \textit{ tal que:} \\
&\left\{ \sum_{i=1}^{N_h}a_{ij}x_i = b_j, \quad j=1,\ldots, N_h,
\right.
\end{aligned}
\label{FM1}
\end{equation}
que en notación matricial, se reduce al sistema lineal

\begin{equation}
\textit{Encontrar } \boldsymbol{u_h} \in \mathbb{R}^{N_h} \textit{tal que: } A\boldsymbol{u_h} = b
\label{FM}
\end{equation}
donde la matriz $A$ y el vector $b$ son los definidos arriba por la forma bilineal $\mathcal{A}$ y el funcional lineal $\mathcal{F}$ respectivamente, entonces la solución de \eqref{FG} es dada por la combinación lineal de las funciones base $\phi_i$ con el vector $\boldsymbol{u_h}$ solución de \eqref{FM}, de aquí

$$ u_h^*=\sum_{i=1}^{N_h}x_i\phi_i.$$
Ahora necesitamos un espacio $V_h$ adecuado, $V_h$ puede depender o no de una partición (triangulación), y puede que dicha partición dependa o no del dominio en el que se encuentra la EDP en estudio, su elección y la de su base es importante para el método, ya que de esto depende que la matriz resultante sea cómoda para trabajar, así que las funciones base deben tener soporte pequeño. Por lo que usaremos el espacio de funciones lineales por partes basada en una triangulación del dominio de la EDP. 

\subsection{El espacio de elementos finitos de funciones lineales por partes}

Sea $(a,b)$ el dominio en donde esta definida una EDP. Llamaremos $\mathcal{T}^h$ a la partición (triangulación) de el intervalo $(a,b)$ , la cual es formada por pequeños intervalos, llamados \textit{elementos}, de la forma $K_i=(x_i,x_{i+1}),$ $i= 1,\ldots, M-1;$ donde los vértices $\{x_i\}_{i=1}^M $ satisfacen
$$ a=x_1<x_2<\ldots < x_{M-1}<x_M=b$$
y los diámetros de los elementos son de tamaño proporcional a $h>0$ (máximo diámetro de los elementos), donde $\{x_1,x_M\}$ son denominados vértices frontera y $\{ x_i \}_{i=2}^{M-1}$ son denominados vértices interiores.

\begin{obs}
Notemos que los elementos $K_i$ cumplen,
\begin{equation*}
\begin{aligned}
K_i \cap K_j =& \emptyset \quad \text{Si } i\neq j \\
\bar{K_i} \cap \bar{K}_j =& \begin{cases}
K_i & \text{Si } i=j \\
\emptyset & \\
\text{Vértice común de } K_i \text{ y } K_j
\end{cases} \\
\bigcup_{i=0,\ldots,N} \bar{K_i} =& [a,b]
\end{aligned}
\end{equation*}
\end{obs}
\noindent El espacio de funciones lineales por partes asociado a la triangulación $\mathcal{T}^h$, notado por  $\mathbb{P}^1(\mathcal{T}^h)$ es,
$$\mathbb{P}^1(\mathcal{T}^h)=\left\{ v \in C(a,b) \quad : \quad v|_K \text{ es un polinomio de grado 1 para todo } K \in \mathcal{T}^h \right\}$$
y también definimos 

$$\mathbb{P}^1_0(\mathcal{T}^h)=\left\{ v \in \mathbb{P}^1(\mathcal{T}^h)\quad : \quad v(x_1)=v(x_M)=0 \right\}.$$

Por lo que usaremos $V = \mathbb{P}^1(\mathcal{T}^h)$ en \eqref{FG} y esta solución sera la aproximación de elementos finitos de \eqref{FG'} en una partición de $(a,b)$ y se denota por $u^h$. Ahora, si $v^h \in \mathbb{P}^1(\mathcal{T}^h)$, entonces

$$ v^h(x) = \displaystyle \sum_{i=1}^{M} v^h(x_i) \phi_i (x) $$
donde $\phi_i, \, i=1,\ldots, M$, son las funciones base y están definidas por,

$$\phi_i(x)= \left\{ \begin{array}{llc}

1, & \text{Si } x=x_i, \\
0, & \text{Si } x=x_j, j \neq i \\
\text{extension lineal}, & \text{Si } x \text{ no es un vértice.}

\end{array} \right. $$

\begin{obs}
El espacio $\mathbb{P}^1(\mathcal{T}^h) $ es generado por las funciones $\{\phi_i\}_{i=1}^M$ y $\mathbb{P}^1_0(\mathcal{T}^h)$ por $\{\phi_i\}_{i=2}^{M-1}$.
\end{obs}

\begin{lema}
Sea $K_i = (x_i,x_{i+1})$ un elemento de la triangulación $\mathcal{T}^h$ de $(a,b)$, entonces

$$\begin{array}{lll}

\phi_i (x)=\dfrac{x_{i+i}-x}{x_{i+1}-x_i}, & \phi_{i+1} (x)=\dfrac{x-x_i}{x_{i+1}-x_i}, & \text{ para todo } x \in K_i.

\end{array} $$

\end{lema}

\begin{obs}
De acuerdo a este lema, decimos que el método de los elementos finitos lineales por partes en una dimensión, tiene dos grados de libertad por cada elemento.

\end{obs}

De lo anterior podemos ver que el vector $\boldsymbol{u}^h \in \mathbb{R}^M$ es dado por los valores de la aproximación $u^h$ en los vértices de la triangulación $\mathcal{T}^h$, es decir,

\begin{equation}
\boldsymbol{u}^h = [u^h(x_1),u^h(x_M), \ldots, u^h(x_M)]^T \in \mathbb{R}^M.
\label{u}
\end{equation}

Y obtenemos un sistema lineal con funciones del espacio de elementos finitos de funciones lineales por partes, dado por $\boldsymbol{u}^h$ como en \eqref{u}, $\boldsymbol{b}$ de \eqref{b}, y $A$ la matriz definida en \eqref{A}.

\begin{equation}
A\boldsymbol{u}^h=\boldsymbol{b}.
\end{equation}

\begin{obs}

Dependiendo de las condiciones de contorno se elige la matriz a trabajar, puede ser, escogiendo todos los vértices, se obtiene una matriz de $M \times M$ (Matriz de Neumann), o si excluimos los vértices de frontera, obtenemos una matriz de dimensión $(M-2) \times(M-2)$ (Matriz de Dirichlet), y de aquí facilitar la solución de dicho sistema.

\end{obs}

Ahora, enunciaremos un lema, que es muy útil para formar la matriz de elementos finitos, pero para eso usaremos la forma bilineal de \eqref{fb}, y tomamos funciones de elementos finitos por partes.

\begin{lema}
\label{MM}

Sea $\mathcal{A}$ la forma bilineal definida en \eqref{fb} y $\mathcal{T}^h$ una triangulación de $(a,b)$. Sean $K_i = (x_i,x_{i+1}), \, i=1,\ldots,M-1$ los elementos de la triangulación. Defina $R_i$ como la matriz $2\times M$ de restricción al elemento $K_i$, es decir, la matriz $R_i$ tiene todas sus entradas nulas, excepto las posiciones $(1,i)$ y $(2,i+1)$ donde su valor es 1. Sea $A_{K_i}$ la matriz local definida por 

\begin{equation}
A_{K_i}=\begin{bmatrix}
\mathcal{A}_{K_i}(\phi_i,\phi_i) & \mathcal{A}_{K_i}(\phi_i,\phi_{i+1}) \\
\mathcal{A}_{K_i}(\phi_{i+1},\phi_i) & \mathcal{A}_{K_i}(\phi_{i+1},\phi_{i+1})
\end{bmatrix}.
\label{AKi}
\end{equation}
donde $\mathcal{A}_{K_i}$, la restricción de la forma bilineal $\mathcal{A}$ al elemento $K_i$, es dada por,

\begin{equation}
{A}_{K_i}(v,w)=int_{K_i}\kappa(x)v'(x)w'(x)dx.
\end{equation}
Y sea $\boldsymbol{b}_{K_i}$, el lado directo local, tal que

\begin{equation}
\boldsymbol{b}_{K_i}=\begin{bmatrix}
\mathcal{F}_{K_i}(\phi_i) \\
\mathcal{F}_{K_i}(\phi_{i+1})
\label{bi}
\end{bmatrix},
\end{equation}
donde $\mathcal{F}_{K_i}$, es la restricción de $\mathcal{F}$ al elemento $K_i$, es dada por,

$$\mathcal{F}_{K_i}(v)=\int_{K_i} f(x)v(x)dx.$$
Tenemos entonces que
\begin{equation}
A=\sum_{i=1}^{M-1}R_i^TA_{K_i}R_i
\label{MA}
\end{equation}
y

\begin{equation}
\boldsymbol{b}=\sum_{i=1}^{M-1}R_i^T\boldsymbol{b}_{K_i}.
\label{vb}
\end{equation}

\end{lema}

Este lema permite montar la matriz $A$ de \eqref{MA}, usando los aportes de cada elemento, y acomodándolos con las matrices $R_i$.

\section{Ejemplos}

Presentaremos algunos ejemplos de como solucionar una EDP por medio del método de elementos finitos, se mostrara el procedimiento, pero algunos cálculos fueron hechos en un programa de MatLab.

\subsection*{Ecuación de Laplace}

Queremos aproximar la solución de la ecuación de Laplace \eqref{EL1D} usando elementos finitos, con los siguientes parámetros, y condiciones de contorno:

\begin{equation}
(\text{Formulación Fuerte})
\begin{aligned}[b]
&\textit{Encontrar } u : [0,1] \to \mathbb{R} \textit{ tal que:} \\
&\left\{
\begin{aligned}
&-u''(x) = -1, \qquad 0<x<1 \\ 
& u(0)=u(1)=1.
\end{aligned} \right.
\end{aligned}
\label{FF}
\end{equation}

\subsubsection*{Formulación Débil}

Para construir la \textit{formulación débil} de \eqref{FF}, primero multiplicamos los dos lados de la igualdad en \eqref{FF} por la función de prueba $v \in C_0^\infty(a,b)$ fija pero arbitraria, luego integramos a ambos lados y obtenemos

\begin{equation}
\begin{aligned}[b]
&\textit{Encontrar } u : [0,1] \to \mathbb{R} \textit{ tal que:} \\
&\left\{
\begin{aligned}
&-\int_0^1 u''(x)v(x)dx &= \int_0^1 f(x)v(x) dx,  \\
& u(0)=u(1)=1,
\end{aligned} \right.
\end{aligned}
\label{FF1}
\end{equation}
luego al usar la fórmula de integración por partes, y como $v \in C_0^\infty(a,b)$, entonces $v(0)=v(1)=0$, por lo que

\begin{align*}
-\int_0^1 u''(x)v(x)dx &= \int_a^b u'(x)v'(x)dx -[u'(1)v(1)- u'(0)v(0)] \\
&= \int_0^1 u'(x)v'(x)dx -[u'(b)0- u'(a)0] \\
&= \int_0^1 u'(x)v'(x)dx.
\end{align*}
Falta cambiar los espacios de funciones que requiere la ecuación de Laplace, el cual es un espacio de \textit{Hilbert}, y $C_0^\infty(a,b)$ no es un espacio de Hilbert. Para esto, escogemos $V=H_0^1(a,b)$ y $U=H^1(a,b)$.

Ahora notamos,

\begin{equation}
\mathcal{A}(u,v)= \int_a^b u'(x)v'(x)dx \quad \text{para toda } v\in V \text{ y } u\in U
\end{equation}
donde $\mathcal{A}:U\times V \to \mathbb{R}$ es una forma bilineal, y

\begin{equation}
\mathcal{F}(u,v)= \int_a^b f(x)v(x)dx \quad \text{para toda } v\in V
\end{equation}
donde $\mathcal{F}: V \to \mathbb{R}$ es un funcional lineal. Entonces decimos que la formulación débil de \eqref{FF} en (a,b)=(0,1) es:

\begin{equation}
(\text{Formulación Débil})
\begin{aligned}[b]
&\textit{Encontrar } u \in H^1(0,1) \textit{ tal que:} \\
&\left\{
\begin{aligned}
&\mathcal{A}(u,v) = \mathcal{F}(v),  \\
&u(0)=u(1)=1.
\end{aligned}
\quad
\begin{aligned}
&\textnormal{para toda } v \in V \\
& 
\end{aligned} \right.
\end{aligned}
\label{FD}
\end{equation}

\subsubsection*{Triangulación}

Usaremos una triangulación uniforme con seis vértices,

$$ \mathcal{T}^h= \{x_1=0, x_2=\frac{1}{5}, x_3=\frac{2}{5}, x_4=\frac{3}{5}, x_5=\frac{4}{5}, x_6=1\}.$$
con lo cual obtendremos 5 elementos en la triangulación $\mathcal{T}^h$,

$$ K_1=\left(0,\frac{1}{5}\right), K_2=\left(\frac{1}{5},\frac{2}{5}\right), K_3=\left(\frac{2}{5},\frac{3}{5}\right), K_4=\left(\dfrac{3}{5},\frac{4}{5}\right), K_5=\left(\dfrac{4}{5},1\right).$$

\subsubsection*{Montaje de la matriz}

Usando el lema \ref{MM}, montamos la matriz A definida en \eqref{MA}, por lo que construimos las matrices $A_{K_i}$ de \eqref{AKi} y las matrices de restricción $R_{K_i}$, para $i=1,\ldots,5.$

\noindent Entonces, para el primer elemento $K_1=(x_1,x_2)=(0,\frac{1}{5})$, las funciones base son,

$$ \phi_1(x)=\dfrac{x_2-x}{x_2-x_1}=(1-5x) \quad \text{y} \quad \phi_2(x)=\dfrac{x-x_1}{x_2-x_1}=5x, \quad x \in (x_1,x_2). $$
Con esto calculamos $A_{K_1}$, de la siguiente manera,

\begin{align*}
&\mathcal{A}_{K_1}(\phi_1,\phi_1) = \int_{K_1} \phi_1'\phi_1'= \int_{x_1}^{x_2} (-5)(-5)= \int_{0}^{\frac{1}{5}} 25 = 5 \\
&\mathcal{A}_{K_1}(\phi_1,\phi_2) = \mathcal{A}_{K_1}(\phi_2,\phi_1) = \int_{K_1} \phi_1'\phi_2'= \int_{x_1}^{x_2} (-5)(5)= \int_{0}^{\frac{1}{5}} -25 = -5
\\
&\mathcal{A}_{K_1}(\phi_2,\phi_2) = \int_{K_1} \phi_2'\phi_2'= \int_{x_1}^{x_2} (5)(5)= \int_{0}^{\frac{1}{5}} 25 = 5
\end{align*}

De esto, vemos que $A_{K_1}=5\begin{bmatrix}
1 & -1 \\
-1 & 1
\end{bmatrix}$. Note que $R_1=\begin{bmatrix}
1&0&0&0&0&0 \\
0&1&0&0&0&0 \\
\end{bmatrix}$ y

$$R_1^TA_{K_1}R_1 = 5 \begin{bmatrix}
1&-1&0&0&0&0 \\
-1&1&0&0&0&0 \\
0&0&0&0&0&0 \\
0&0&0&0&0&0 \\
0&0&0&0&0&0 \\
0&0&0&0&0&0 
\end{bmatrix}.$$

Para el elemento $K_2=(x_1,x_2)=(\frac{1}{5},\frac{2}{5})$, las funciones base son,

$$ \phi_1(x)=\dfrac{x_2-x}{x_2-x_1}=(2-5x) \quad \text{y} \quad \phi_2(x)=\dfrac{x-x_1}{x_2-x_1}=5x-1, \quad x \in (x_1,x_2). $$

Con esto calculamos $A_{K_2}$, de la siguiente manera,
\begin{align*}
&\mathcal{A}_{K_2}(\phi_1,\phi_1) = \int_{K_2} \phi_1'\phi_1'= \int_{x_1}^{x_2} (-5)(-5)= \int_{\frac{1}{5}}^{\frac{2}{5}} 25 = 5. \\
&\mathcal{A}_{K_2}(\phi_1,\phi_2) = \mathcal{A}_{K_2}(\phi_2,\phi_1) \int_{K_2} \phi_1'\phi_2'= \int_{x_1}^{x_2} (-5)(5)= \int_{\frac{1}{5}}^{\frac{2}{5}} -25 = -5. \\
&\mathcal{A}_{K_2}(\phi_2,\phi_2) = \int_{K_2} \phi_2'\phi_2'= \int_{x_1}^{x_2} (5)(5)= \int_{\frac{1}{5}}^{\frac{2}{5}} 25 = 5.
\end{align*}

De esto, vemos que $A_{K_2}=5\begin{bmatrix}
1 & -1 \\
-1 & 1
\end{bmatrix}$. Note que $R_2=\begin{bmatrix}
0&1&0&0&0&0 \\
0&0&1&0&0&0 \\
\end{bmatrix}$ y

$$R_1^TA_{K_1}R_1 = 5 \begin{bmatrix}
0&0&0&0&0&0 \\
0&1&-1&0&0&0 \\
0&-1&1&0&0&0 \\
0&0&0&0&0&0 \\
0&0&0&0&0&0 \\
0&0&0&0&0&0 
\end{bmatrix}.$$

Realizando los mismos cálculos para $i=3,\ldots,5$, con las funciones base para 

$$K_3 \rightarrow \phi_1(x)= 3-5x \qquad \phi_2(x)= 5x-2$$
$$K_4 \rightarrow \phi_1(x)= 4-5x \qquad \phi_2(x)= 5x-3$$
$$K_5 \rightarrow \phi_1(x)= 5-5x \qquad \phi_2(x)= 5x-4$$

así vemos que $A_{K_i}$ es

$$A_{K_i}=5\begin{bmatrix}
1 & -1 \\
-1 & 1
\end{bmatrix}, \quad i=1,\ldots,5.$$

Con las matrices locales $A_{K_i}$, armamos la matriz global $A$ como

$$A=\sum_{i=1}^{5}R_i^TA_{K_i}R_i= 5 \begin{bmatrix}
 1 & -1  &  0  &  0  &  0  &  0 \\
-1 &(1+1)& -1  &  0  &  0  &  0 \\
 0 &  -1 &(1+1)& -1  &  0  &  0 \\
 0 &  0  & -1  &(1+1)& -1  &  0 \\
 0 &  0  & 0   & -1  &(1+1)& -1 \\
 0 & 0   & 0   &  0  & -1  &  1 
\end{bmatrix} = $$ $$ = 5 \begin{bmatrix}
 1 & -1  &  0  &  0  &  0  &  0 \\
-1 &  2  & -1  &  0  &  0  &  0 \\
 0 &  -1 &  2  & -1  &  0  &  0 \\
 0 &  0  & -1  &  2  & -1  &  0 \\
 0 &  0  & 0   & -1  &  2  & -1 \\
 0 & 0   & 0   &  0  & -1  &  1 
\end{bmatrix}. $$
que se denomina matriz global de Neumann (dimension $6 \times 6$), ya que usamos todos los vértices.

\subsubsection*{Montaje del vector $\boldsymbol{b}$}

Usando la definición del vector $\boldsymbol{b} $ en \eqref{vb} en el lema \ref{MM}, entonces calculamos $\boldsymbol{b}_{K_i} $ de \ref{bi} para los $K_i$ de nuestra triangulación $\mathcal{T}^h$. Vemos que $f(x)=-1$.

\noindent Para $K_1=(x_1,x_2)=(0,\frac{1}{5})$, calculamos $\boldsymbol{b}_{K_1}$, de la siguiente manera,
\begin{align*}
&\mathcal{F}_{K_1}(\phi_1) = \int_{K_1} -1\phi_1(x) = -\int_{x_1}^{x_2} \dfrac{x_2-x}{x_2-x_1} = -\int_{0}^{\frac{1}{5}} 1-5x = \left.\left(\frac{5x^2}{2}-x\right) \right|_0^\frac{1}{5} = -\frac{1}{10}  \\
&\mathcal{F}_{K_1}(\phi_2) = \int_{K_1} -1\phi_2(x) = -\int_{x_1}^{x_2} \dfrac{x-x_1}{x_2-x_1} = -\int_{0}^{\frac{1}{5}} 5x = \left.\left(-\frac{5x^2}{2}\right) \right|_0^\frac{1}{5} = -\frac{1}{10}  \\
\end{align*}
De esto, vemos que $\boldsymbol{b}_{K_1}=-\begin{bmatrix}
\frac{1}{10} \\
\frac{1}{10}
\end{bmatrix}$. 

\noindent Para $K_2=(x_1,x_2)=(\frac{1}{5},\frac{2}{5})$, calculamos $\boldsymbol{b}_{K_2}$, de la siguiente manera,

\begin{align*}
&\mathcal{F}_{K_2}(\phi_1) = \int_{K_1} -1\phi_1(x) = -\int_{x_1}^{x_2} \dfrac{x_2-x}{x_2-x_1} = -\int_{\frac{1}{5}}^{\frac{2}{5}} 2-5x = \left.\left(\frac{5x^2}{2}-2x\right) \right|_\frac{1}{5}^\frac{2}{5} = -\frac{1}{10}  \\
&\mathcal{F}_{K_1}(\phi_2) = \int_{K_1} -1\phi_2 = -\int_{x_1}^{x_2} \dfrac{x-x_1}{x_2-x_1} = -\int_{\frac{1}{5}}^{\frac{2}{5}}2 - 5x = \left.\left(\frac{5x^2}{2}-2x\right) \right|_\frac{1}{5}^\frac{2}{5} = -\frac{1}{10}  \\
\end{align*}
De esto, vemos que $\boldsymbol{b}_{K_2}=-\begin{bmatrix}
\frac{1}{10} \\
\frac{1}{10}
\end{bmatrix}$. 

\noindent Análogamente, vemos que $\boldsymbol{b}_{K_i}$ es

$$\boldsymbol{b}_{K_i}=-\begin{bmatrix}
\frac{1}{10} \\
\frac{1}{10}
\end{bmatrix}, \quad \text{para } i=1,\ldots,5.$$
Por lo que $\boldsymbol{b}$ sera,

$$\boldsymbol{b}=\sum_{i=1}^{5}R_i^T\boldsymbol{b}_{K_i}=  \begin{bmatrix}
-\frac{1}{10} \\
-\frac{1}{10}-\frac{1}{10}\\
-\frac{1}{10}-\frac{1}{10}\\
-\frac{1}{10}-\frac{1}{10}\\
-\frac{1}{10}-\frac{1}{10}\\
-\frac{1}{10}
\end{bmatrix} = - \begin{bmatrix}
1/10 \\   
1/5    \\ 
1/5     \\
1/5     \\
1/5     \\
1/10  
\end{bmatrix}. $$

\subsubsection*{Solución del sistema lineal}

Nos queda solucionar el sistema $A\boldsymbol{u}_h=\boldsymbol{b}$, que se puede resolver en MatLab. Tenemos,

$$A\boldsymbol{u}_h = 5 \begin{bmatrix}
 1 & -1  &  0  &  0  &  0  &  0 \\
-1 &  2  & -1  &  0  &  0  &  0 \\
 0 &  -1 &  2  & -1  &  0  &  0 \\
 0 &  0  & -1  &  2  & -1  &  0 \\
 0 &  0  & 0   & -1  &  2  & -1 \\
 0 & 0   & 0   &  0  & -1  &  1 
\end{bmatrix} \begin{bmatrix}
u_h(x_1) \\
u_h(x_2) \\
u_h(x_3) \\
u_h(x_4) \\
u_h(x_5) \\
u_h(x_6) \\
\end{bmatrix} =  - \begin{bmatrix}
1/10 \\   
1/5    \\ 
1/5     \\
1/5     \\
1/5     \\
1/10  
\end{bmatrix} $$
Pero, por condiciones de frontera sabemos que la aproximación $u_h$ debe cumplir $u_h(x_0)=1$ y $u_h(x_6)=1$ entonces,

$$\begin{bmatrix}
u_h(x_1) \\
u_h(x_2) \\
u_h(x_3) \\
u_h(x_4) \\
u_h(x_5) \\
u_h(x_6) \\
\end{bmatrix} = \begin{bmatrix}
1 \\
0 \\
0 \\
0 \\
0 \\
1 \\
\end{bmatrix} + \begin{bmatrix}
0 \\
u_h(x_2) \\
u_h(x_3) \\
u_h(x_4) \\
u_h(x_5) \\
0 \\
\end{bmatrix} $$
y el termino conocido lo pasamos al otro lado, así

$$ 5 \begin{bmatrix}
 1 & -1  &  0  &  0  &  0  &  0 \\
-1 &  2  & -1  &  0  &  0  &  0 \\
 0 &  -1 &  2  & -1  &  0  &  0 \\
 0 &  0  & -1  &  2  & -1  &  0 \\
 0 &  0  & 0   & -1  &  2  & -1 \\
 0 & 0   & 0   &  0  & -1  &  1 
\end{bmatrix} \begin{bmatrix}
0 \\
u_h(x_2) \\
u_h(x_3) \\
u_h(x_4) \\
u_h(x_5) \\
0 \\
\end{bmatrix} =  - \begin{bmatrix}
1/10 \\   
1/5    \\ 
1/5     \\
1/5     \\
1/5     \\
1/10  
\end{bmatrix} - \begin{bmatrix}
5 \\
-5 \\
0 \\
0 \\
-5 \\
5 
\end{bmatrix} $$
y de aquí, podemos escoger la matriz de los vértices interiores (Matriz de Dirichlet), que es de dimensión $4 \times 4$, la cual es

$$ 5 \begin{bmatrix}
2  & -1  &  0  &  0 \\
-1 &  2  & -1  &  0 \\
0  & -1  &  2  & -1 \\
0  &  0  & -1  &  2 \\
\end{bmatrix} \begin{bmatrix}
u_h(x_2) \\
u_h(x_3) \\
u_h(x_4) \\
u_h(x_5) \\
\end{bmatrix} =  \begin{bmatrix}
24/5    \\ 
-1/5     \\
-1/5     \\
24/5     \\
\end{bmatrix} $$
que tiene como solución $\begin{bmatrix}
u_h(x_2) \\
u_h(x_3) \\
u_h(x_4) \\
u_h(x_5) \\
\end{bmatrix} =  \begin{bmatrix}
23/25    \\ 
22/25     \\
22/25     \\
23/25     \\
\end{bmatrix}. $ Por lo que la aproximación de elementos finitos es,

$$\boldsymbol{u}_h = \begin{bmatrix}
u_h(x_1) \\
u_h(x_2) \\
u_h(x_3) \\
u_h(x_4) \\
u_h(x_5) \\
u_h(x_6) \\
\end{bmatrix} = \begin{bmatrix}
1 \\   
23/25    \\ 
22/25     \\
22/25     \\
23/25     \\
1  
\end{bmatrix}, $$
y $u_h$ es la combinación lineal de las funciones base con el vector $\boldsymbol{u}_h$,

$$ u_h = \phi_1 + \frac{23}{25}\phi_2 + \frac{22}{25}\phi_3 + \frac{22}{25}\phi_4 + \frac{23}{25}\phi_5 + \phi_6.$$
La solución exacta se obtiene integrando dos veces y aplicando las condiciones de contorno y es,

$$ u(x)=\dfrac{t^2}{2} - \frac{t}{2} + 1.$$

\begin{figure}[hrt]
  \centering
    \includegraphics[width=0.93\textwidth]{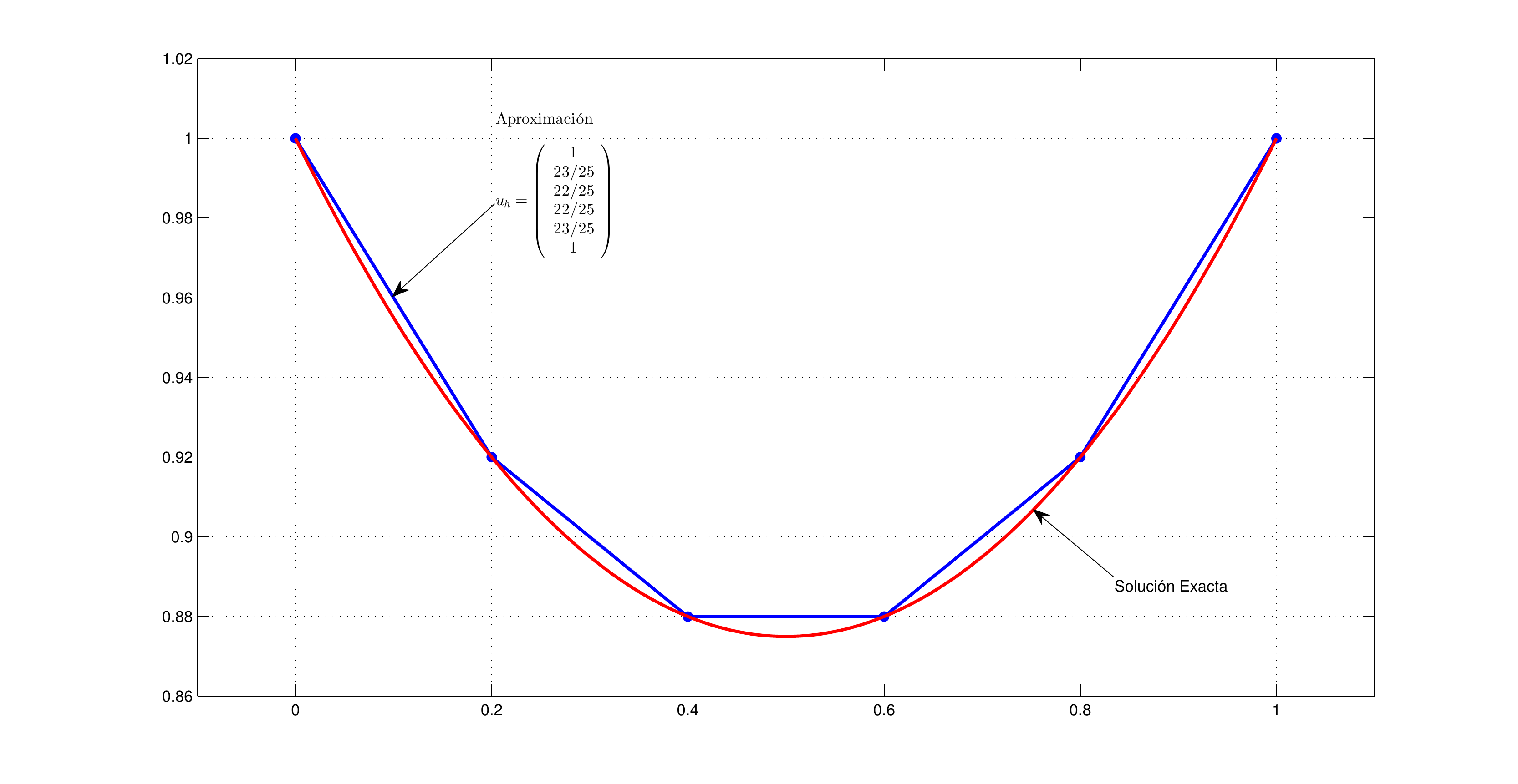}
  \caption{Aproximación $u_h$ vs. Solución Exacta.}
  \label{fL}
\end{figure}

Con esto podemos comparar, que tanto se aproxima el método a su solución, por ejemplo en la Figura \ref{fL}.

Se nota una cierta aproximación con $h=1/6$, sabemos que cuando $h \to 0$ la aproximación va a ser mas exacta, pero el sistema lineal sera de dimension mayor, y por ende, mas difícil su solución.

\chapter{Método de los elementos finitos en dos dimensiones}

En este capítulo veremos el MEF aplicado a una ecuación diferencial en dos dimensión, en un dominio $\Omega$ encontraremos los espacios adecuados para resolver el problema, sus formulaciones, y trabajaremos con el espacio de funciones lineales por partes en dos dimensiones (para facilitar la introducción), pero este puede cambiar dependiendo las necesidades del problema.

El dominio físico $\Omega \subset \mathbb{R}^2$ donde la EDP es formulada es un subconjunto abierto, conexo y \textit{Lipschitz.} Para facilitar la escritura denotaremos $x=(x_1,x_2)$ y $\partial \Omega$ a la frontera del dominio $\Omega$. Los espacios de funciones que trabajaremos serán los definidos en la sección \ref{EP}.

\section{Formulación Débil}

Para deducir la formulación débil de una ecuación diferencial en el dominio $\Omega \subset \mathbb{R}^2$ tenemos que,

\begin{enumerate}
\item Suponer que existe $u$ solución de nuestra ecuación diferencial, multiplicar los dos lados de la ecuación por una \textit{función de prueba} $v \in C_0^\infty(\Omega)$ e integrar a ambos lados de la igualdad.

\item Usar la fórmula de integración por partes en $\mathbb{R}^2$ (primera identidad de Green) y las condiciones iniciales (de contorno) para obtener expresiones que necesiten solamente de derivadas del menor orden posible.

\item Cambiar $v \in C_0^\infty(\Omega)$ por $v \in V = H_0^1(\Omega).$
\end{enumerate}

Al obtener dicha \textit{formulación débil}, podemos garantizar la existencia de la solución de esta formulación, que se conoce como \textit{solución débil} de la ecuación original.
La solución del problema original se denomina \textit{solución fuerte}, y en este caso se cumplen las igualdades en el problema para todo $x \in \Omega$, y por lo tanto se pueden calcular sus derivadas.

Para probar la existencia de \textit{soluciones débiles} se usan resultados como el Teorema de Lax-Milgram, o el teorema de Babuska enunciado en la Sección \ref{LM} y demostrado en \cite{Evans}. 

\begin{ejem}
Vamos a construir la formulación débil, para la EDP definida en \eqref{EE2D}, para esto primero multiplicamos los dos lados de la igualdad en \eqref{EE2D} por la función de prueba $v \in C_0^\infty(\Omega)$ fija pero arbitraria, luego integramos a ambos lados, aplicamos la formula de integración por partes y obtenemos

\begin{equation}
(\text{Formulación Fuerte})
\begin{aligned}[b]
&\textit{Encontrar } u : \Omega \to \mathbb{R} \textit{ tal que:} \\
&\left\{
\begin{aligned}
-\int_\Omega\text{div }(\kappa(x)\nabla u(x))v(x)dx &= \int_\Omega f(x)v(x)dx,  \\
u(x) &= g(x),
\end{aligned}
\quad
\begin{aligned}
&\textnormal{para } x\in \Omega, \\
&\textnormal{para } x\in \partial\Omega,
\end{aligned} \right.
\label{EE2D1}
\end{aligned}
\end{equation}
Luego, al usar la fórmula de integración por partes, y como $v \in C_0^\infty(\Omega)$, entonces $v(x)=0$ para todo $x \in \partial\Omega$, por lo que
\begin{align*}
-\int_\Omega\text{div }(\kappa(x)\nabla u(x))v(x) dx &= \int_\Omega \nabla u(x) \kappa \nabla v(x) dx - \int_{\partial \Omega} v(x) (\nabla u(x) \kappa(x) \cdot \eta(x)) dS  \\
&= \int_\Omega \nabla u(x) \kappa(x) \nabla v(x) dx - 0\\
&= \int_\Omega \nabla u(x) \kappa(x) \cdot \nabla v(x) dx 
\end{align*}
Entonces, \eqref{EE2D1} se escribe como

\begin{equation}
\begin{aligned}[b]
&\textit{Encontrar } u : \Omega \to \mathbb{R} \textit{ tal que:} \\
&\left\{
\begin{aligned}
\int_\Omega \nabla u(x) \kappa(x) \cdot \nabla v(x) dx &= \int_\Omega f(x)v(x) dx,  \\
u(x) &= g(x),
\end{aligned}
\quad
\begin{aligned}
&\textnormal{para toda } v \in C_0^\infty(\Omega) \\
&\textnormal{para } x \in \partial \Omega.
\end{aligned} \right.
\end{aligned}
\label{EE2D2}
\end{equation}
Para que la ecuación \eqref{EE2D2}, sea la formulación débil de \eqref{EE2D1}, falta escoger los espacios de funciones $U=H^1(\Omega)$ y $V=H_0^1(\Omega)$ adecuados para la solución, donde $u \in U$, $v \in V$ y $C_0^\infty(\Omega) \subset V$. Introducimos la notación,

\begin{equation}
\mathcal{A}(u,v)= \int_\Omega \nabla u(x) \kappa(x) \cdot \nabla v(x) dx = \int_\Omega f(x)v(x) dx \quad \text{para toda } v\in V \text{ y } u\in U,
\label{fb1}
\end{equation}
donde $\mathcal{A}:U\times V \to \mathbb{R}$ es una forma bilineal, y

\begin{equation}
\mathcal{F}(v)= \int_\Omega f(x)v(x)dx \quad \text{para toda } v\in V,
\end{equation}
donde $\mathcal{F}: V \to \mathbb{R}$ es un funcional lineal. Entonces decimos que la formulación débil de \eqref{EEFF} es la ecuación \eqref{EEFD1}, pero puesta en los espacios U y V adecuados. La formulación débil puede ser escrita como

\begin{equation}
(\text{Formulación Débil})
\begin{aligned}[b]
&\textit{Encontrar } u \in H_0^1(\Omega) \textit{ tal que:} \\
&\left\{
\begin{aligned}
\mathcal{A}(u,v) &= \mathcal{F}(v),  \\
u &= g,
\end{aligned}
\quad
\begin{aligned}
&\textnormal{para toda } v \in H_0^1(\Omega) \\
&\textnormal{para } x \in \partial \Omega
\end{aligned} \right.
\end{aligned}
\label{EE2DD}
\end{equation}
\end{ejem}

\section{Formulación de Galerkin}

Asumiendo que las formulaciones débiles que trabajaremos son de la forma de la ecuación \eqref{EE2DD}, donde $H_0^1(\Omega)$ es un espacio de dimensión infinita, entonces debemos cambiar el espacio de dimension infinita de \eqref{EE2DD} por espacios de dimensión finita (espacios de elementos finitos).

Considere $V^h$, $h>0$, subespacios de dimension finita tal que $V^h \subset V$. El parámetro $h$ indica el tamaño del espacio, a menor valor $h$, mayor es la dimensión de $V^h$. Lo que se necesita es que $V^h \xrightarrow [h \to 0]\, V$ en algún sentido. Con este $V_h \subset V$ de dimensión finita, la formulación de Galerkin es

\begin{equation}
\begin{aligned}[b]
&\textit{Encontrar } u_h \in V^h \textit{ tal que:} \\
&\left\{
\begin{aligned}
\mathcal{A}(u_h,v_h) &= \mathcal{F}(v_h),  \\
u &= g,
\end{aligned}
\quad
\begin{aligned}
&\textnormal{para toda } v_h \in V^h \\
&\textnormal{para } x \in \partial \Omega
\end{aligned} \right.
\end{aligned}
\label{EE2DG}
\end{equation}
Escogemos
$$\{\phi_1,\phi_2, \ldots , \phi_{N_h} \} \text{ donde } N_h \text{ es la dimensión de } V_h,$$
como una base de $V_h$, y como en el capítulo anterior obtenemos el sistema lineal $N_h \times N_h$
\begin{equation}
\textit{Encontrar } \boldsymbol{u_h} \in \mathbb{R}^{N_h} \textit{tal que: } A\boldsymbol{u_h} = b
\label{FM2D}
\end{equation}
tal que la forma bilineal $\mathcal{A}$ tiene una representación matricial, la matriz $A=[a_{ij}]$ de dimensión $N_h \times N_h$ es de la forma
\begin{equation}
a_{ij}=\mathcal{A}(\phi_i,\phi_j), \quad i,j=1,\ldots,N_h.
\label{A1}
\end{equation}  
y
\begin{equation}
b=\begin{bmatrix}
b_1 \\
b_2 \\
\vdots \\
b_{N_h}
\end{bmatrix} \in \mathbb{R}^{N_h}, 
\quad b_j=\mathcal{F}(\phi_j), \quad j=1,\ldots,n,
\label{b1}
\end{equation} 
y el vector 
$$\boldsymbol{u_h} = \begin{bmatrix}
x_1 \\
x_2 \\
\vdots \\
x_{N_h}
\end{bmatrix} \in \mathbb{R}^{N_h}. $$
$$u_h^*=x_1\phi_1 + x_2\phi_2 + \cdots + x_{N_h}=\displaystyle \sum_{i=1}^{N_h}x_i\phi_i.$$

Ahora necesitamos un espacio $V_h$ adecuado, su elección y la de su base es importante para el método, ya que de esto depende que la matriz resultante sea cómoda para trabajar, así que las funciones base deben tener soporte pequeño. Por lo que para introducir el método usaremos el espacio de funciones lineales por partes basada en una triangulación del dominio de la EDP. 

\subsection{El espacio de elementos finitos de funciones lineales por partes}

Sea $ \Omega \subset \mathbb{R}^2$ un dominio poligonal en donde esta definida una EDP. Llamaremos $\mathcal{T}^h$ a la triangulación del dominio $\Omega$, la cual es formada por subconjuntos disyuntos, llamados \textit{elementos}, tal que 
$$\bigcup_{i=1}^{N_h}\bar{K}_i=\bar{\Omega}, \quad K_i\cap K_j = \emptyset \text{ para } i \neq j \text{ y } h=m\max_{1\leq i \leq N_h} \text{diámetro}(K_i)$$
donde $\mathcal{T}^h $ triangulación geométricamente conforme, es decir que los elementos $K_i$ cumplen,
\begin{equation*}
\begin{aligned}
K_i \cap K_j =& \emptyset \quad \text{Si } i\neq j ,\\
\bar{K_i} \cap \bar{K}_j =& \begin{cases}
\text{Un triángulo completo} \\
\text{Un lado completo} \\
\text{Un vértice} \\
\emptyset, 
\end{cases} \\
\bigcup_{i=0,\ldots,N} \bar{K_i} =& \Omega.
\end{aligned}
\end{equation*}

El espacio de funciones lineales por partes asociado a la triangulación $\mathcal{T}^h$, notado por  $\mathbb{P}^1(\mathcal{T}^h)$ es el definido en \eqref{P1} como,
\begin{equation*}
\mathbb{P}^1(\mathcal{T}^h)=\left\{
\begin{aligned}
 v \in C(\Omega) \quad : & \ v|_K \text{ es un polinomio en dos variables de grado total 1  } \\
&\text{ para todo elemento de la triangulación } \mathcal{T}^h. 
\end{aligned}
\right\}
\end{equation*}
y también definimos
$$\mathbb{P}^1_0(\mathcal{T}^h)=\left\{ v \in \mathbb{P}^1(\mathcal{T}^h)\quad : \quad v(x)=0 \text{ para todo } x \in \partial \Omega \right\}.$$

Por lo que usaremos $V = \mathbb{P}^1(\mathcal{T}^h)$ en \eqref{EE2DG} y esta solución sera la aproximación de elementos finitos en una triangulación $\mathcal{T}^h$ de $\Omega$ y se denota por $u^h$. Ahora, si $v^h \in \mathbb{P}^1(\mathcal{T}^h)$, entonces

$$ v^h(x) = \displaystyle \sum_{i=1}^{N_h} v^h(x_i) \phi_i (x) $$
donde $N_h$ es el número de vértices de la triangulación y $\phi_i, \, i=1,\ldots, N_h$, son las funciones base definidas por,

$$\phi_i(x)= \left\{ \begin{array}{llc}

1, & \text{Si } x=x_i, \\
0, & \text{Si } x=x_j, j \neq i \\
\text{extension lineal}, & \text{Si } x \text{ no es un vértice.}

\end{array} \right. $$

\begin{obs}
El espacio $\mathbb{P}^1(\mathcal{T}^h) $ es generado por las funciones $\{\phi_i\}_{i=1}^{N_h}$ y $\mathbb{P}^1_0(\mathcal{T}^h)$ por $\{\phi_i \ : \ x_i \text{ es un vértice interior}\}$.
\end{obs}

\begin{lema}
Sea $\hat{K}$  el triangulo de referencia y $\varphi_1,\varphi_2,\varphi_3$ sus funciones base definidas como
$$\varphi_1(x)=1-x_1-x_2, \quad \varphi_2=x_1, \quad \varphi_3=x_2$$
para todo $x\in \hat{K}.$
\end{lema}
De aquí podemos notar que el método de los elementos finitos lineales por partes en dos dimensiones tiene tres grados de libertad por cada elemento.

\noindent De lo anterior podemos ver que el vector $\boldsymbol{u}^h \in \mathbb{R}^M$ es dado por los valores de la aproximación $u^h$ en los vértices de la triangulación $\mathcal{T}^h$, es decir,

\begin{equation}
\boldsymbol{u}^h = [u^h(x_1),u^h(x_M), \ldots, u^h(x_M)]^T \in \mathbb{R}^M.
\label{u1}
\end{equation}
Y obtenemos un sistema lineal con funciones del espacio de elementos finitos de funciones lineales por partes, dado por $\boldsymbol{u}^h$ como en \eqref{u1}, $\boldsymbol{b}$ de \eqref{b1}, y $A$ la matriz definida en \eqref{A1}.

\begin{equation}
A\boldsymbol{u}^h=\boldsymbol{b}
\end{equation}

\begin{obs}
\label{MD}
Dependiendo de las condiciones de contorno se elige la matriz a trabajar, puede ser, escogiendo todos los vértices, se obtiene una matriz de $M \times M$ (Matriz de Neumann), o si excluimos los vértices de frontera, obtenemos una matriz de dimensión menor (Matriz de Dirichlet), y de aquí facilitar la solución de dicho sistema.

\end{obs}

Ahora, enunciaremos un lema, que es muy útil para formar la matriz de elementos finitos, pero para eso usaremos la forma bilineal de \eqref{fb1}, y tomamos funciones de elementos finitos por partes.

\begin{lema}
\label{MM1}

Sea $\mathcal{A}$ la forma bilineal definida en \eqref{fb1} y $\mathcal{T}^h$ una triangulación en $\Omega$. Sean $K_i, \, i=1,\ldots,N_h$ los elementos de la triangulación. Defina $R_i$ como la matriz $3\times N_h$ de restricción al elemento $K_i$, es decir, la matriz $R_i$ tiene todas sus entradas nulas, excepto las posiciones $(1,i_1)$,$(2,i_2)$ y $(3,i_3)$ donde su valor es 1. Sea $A_{K_i}$ la matriz local definida por 

\begin{equation}
A_{K_i}=\begin{bmatrix}
\mathcal{A}_{K_i}(\phi_{i_1},\phi_{i_1}) & \mathcal{A}_{K_i}(\phi_{i_1},\phi_{i_2}) & \mathcal{A}_{K_i}(\phi_{i_1},\phi_{i_3})\\
\mathcal{A}_{K_i}(\phi_{i_2},\phi_{i_1}) & \mathcal{A}_{K_i}(\phi_{i_2},\phi_{i_2}) & \mathcal{A}_{K_i}(\phi_{i_2},\phi_{i_3})\\
\mathcal{A}_{K_i}(\phi_{i_3},\phi_{i_1}) & \mathcal{A}_{K_i}(\phi_{i_3},\phi_{i_2}) & \mathcal{A}_{K_i}(\phi_{i_3},\phi_{i_3})\\
\end{bmatrix}
\label{AKi1}
\end{equation}
donde $\mathcal{A}_{K_i}$, la restricción de la forma bilineal $\mathcal{A}$ al elemento $K_i$, es dada por,

\begin{equation}
{A}_{K_i}(v,w)=int_{K_i}\kappa(x)\nabla v(x)\cdot \nabla w(x)dx.
\end{equation}
Y sea $\boldsymbol{b}_{K_i}$, el lado directo local, tal que

\begin{equation}
\boldsymbol{b}_{K_i}=\begin{bmatrix}
\mathcal{F}_{K_i}(\phi_{i_1}) \\
\mathcal{F}_{K_i}(\phi_{i_2}) \\
\mathcal{F}_{K_i}(\phi_{i_3})
\label{bi1}
\end{bmatrix}
\end{equation}
donde $\mathcal{F}_{K_i}$, es la restricción de $\mathcal{F}$ al elemento $K_i$, es dada por,

$$\mathcal{F}_{K_i}(v)=\int_{K_i} f(x)v(x)dx.$$
Tenemos entonces que
\begin{equation}
A=\sum_{i=1}^{N_h}R_i^TA_{K_i}R_i
\label{MA1}
\end{equation}
y
\begin{equation}
\boldsymbol{b}=\sum_{i=1}^{N_h}R_i^T\boldsymbol{b}_{K_i}.
\label{vb1}
\end{equation}

\end{lema}

Este lema permite montar la matriz $A$ de \eqref{MA1} y el vector \eqref{vb1}, usando los aportes de cada elemento, y acomodándolos con las matrices $R_i$, para montar el sistema lineal y así obtener la aproximación por elementos finitos de la EDP en estudio.

\part{Aplicación del Método de Elementos Finitos (MEF) para las Ecuaciones de Stokes y 
de Advección-Difusión}

\chapter{Ecuación de Stokes}

El objetivo de este capítulo es mostrar el desarrollo del MEF para la ecuación de Stokes en dos dimensiones, la cual quiere mostrar o modelar, el movimiento de un fluido en algún dominio especifico, ya sea por una fuerza externa que actúa sobre el, o por la presión que ejerce este sobre el dominio. Para mas detalles sobre esta ecuación, y su aproximación por el método de elementos finitos se puede consultar \cite{Stokes,Rojo,Cla,FreeFem}.

\section{Introducción}

La ecuación de Stokes (también llamada Stokes estacionaria) es un sistema de ecuaciones diferenciales parciales lineales derivada de las ecuaciones de Navier-Stokes (véase \cite{Stokes}), las cuales pretender determinar el flujo de un fluido viscoso (denominado fluido newtoniano) cuando el número de Reynolds \footnote{El número de Reynolds (Re) es un número adimensional utilizado en mecánica de fluidos, diseño de reactores y fenómenos de transporte para caracterizar el movimiento de un fluido, este  relaciona la densidad, viscosidad, velocidad y dimensión típica de un flujo en una expresión adimensional, que interviene en numerosos problemas de dinámica de fluidos.} es muy bajo, aplicando los principios de la mecánica y la termodinámica a un volumen fluido. Esta ecuación se presenta en la atmósfera, corrientes de agua, etc; su relación con las ecuaciones de Navier-Stokes se basa en que la ecuación de Stokes es una linearización estacionaria de estas.
Para esta ecuación no se dispone de una solución analítica general, por lo que recurrimos al método de elementos finitos (MEF) para encontrar una aproximación de la solución.

Para el desarrollo de MEF, presentaremos la ecuación de Stokes en su \textit{formulación fuerte} dada por,
\begin{equation}
\begin{aligned}[b]
&\textit{Encontrar } \vec{u}=(u_1,u_2) \textit{ y p tal que:} \\
&\left\{
\begin{aligned}
-\mu \Delta \vec{u} + \nabla p  = \vec{f}  \quad
 &\text{para } x \in \Omega,& \\
\text{div }\vec{u}=0 \quad
 &\text{para } x \in \Omega,& \\
 \int_\Omega p = 0 \quad
 & \text{Condición de solubilidad} \\
 \vec{u}(x)=\vec{h}_D(x) \quad
 &\text{para }  x  \in \Gamma_D \subset \partial\Omega,	 \\
\nabla \vec{u}(x)\cdot\vec{\eta} = \vec{h}_N(x) \quad
 &\text{para } x \in \Gamma_N \subset \partial \Omega,
\end{aligned}\right.
\label{ESF}
\end{aligned}
\end{equation}
donde $x=(x_1,x_2)$, $\Omega \subset \mathbb{R}^2$, $\vec{u}=(u_1,u_2)$ es la velocidad, $p$ es la presión, $\mu $ es el coeficiente de viscosidad del fluido, y $f$ es una fuerza externa que afecta al fluido. La condición de frontera la separamos en dos partes, la condición de frontera de Dirichlet denotada $\Gamma_D$, y la condición de frontera de Neumann denotada $\Gamma_N$, donde $\partial \Omega = \Gamma_D \cup \Gamma_N.$

Ahora con la formulación \eqref{ESF}, el siguiente paso para el desarrollo del MEF es construir su formulación débil, que sera lo que tratara la siguiente sección.

\section{Formulación débil de la ecuación de Stokes}

En esta sección construiremos la \textit{formulación débil} para la ecuación de Stokes, por lo tanto suponemos que $\vec{u}$, es una solución de la ecuación \eqref{ESF} y escogemos $H^1_0(\Omega)$ como el espacio de prueba. Entonces, multiplicamos la primera ecuación de \eqref{ESF} por $\vec{v}\in [H^0_1(\Omega]^2,$ función de prueba. Obtenemos
$$-\mu \Delta \vec{u} \vec{v} + \nabla p \vec{v}  = \vec{f}\vec{v}.$$
Luego, integramos sobre $\Omega$ vemos que,
$$-\mu \int_\Omega \Delta \vec{u} \vec{v} + \int_\Omega \nabla p \vec{v}  = \int_\Omega \vec{f}\vec{v}.$$
Por la fórmula de Green tenemos que,
$$\mu \int_\Omega \nabla \vec{u} \nabla\vec{v} - \mu \cancelto{0}{\int_{\partial\Omega} (\nabla \vec{u}\cdot \eta)\cdot \vec{v}} - \int_\Omega \text{div } \vec{v} p + \cancelto{0}{\int_{\partial\Omega} \vec{v}\cdot\eta\cdot p}  = \int_\Omega \vec{f}\vec{v},$$
en donde hemos usado el hecho que $v=0$ en $\partial \Omega.$
Después de eliminar los términos nulos, queda la expresión
$$\mu \int_\Omega \nabla \vec{u} \nabla\vec{v} - \int_\Omega \text{div } \vec{v} p 	 = \int_\Omega \vec{f}\vec{v}.$$
El segundo termino de la ecuación \eqref{ESF} es,
$\text{div }\vec{u}=0 $
multiplicando por $q \in L^2(\Omega)$ queda como,
$$\int_\Omega \text{(div } \vec{u}) q  = 0.$$
Con lo anterior, escribimos la formulación débil de la ecuación de Stokes como,

\begin{equation}
\begin{aligned}[b]
&\textit{Encontrar } \vec{u} \in [H^1_0(\Omega)]^2  \textit{ y } p \in L^2(\Omega) \text{ tal que:} \\
&\left\{
\begin{aligned}
\mu \int_\Omega \nabla \vec{u} \nabla\vec{v} - \int_\Omega \text{(div } \vec{v}) p 	 = \int_\Omega \vec{f}\vec{v}  \quad
 &\text{para toda } \vec{v} \in [H^1_0(\Omega)]^2,& \\
-\int_\Omega \text{(div } \vec{u}) q  = 0 \quad
 &\text{para toda } q \in L^2(\Omega),& \\
\int_\Omega p = 0 \quad
 & \text{Condición de solubilidad} \\
 \vec{u}(x)=\vec{h}_D(x) \quad
 &\text{para }  x  \in \Gamma_D \subset \partial\Omega,	 \\
\nabla \vec{u}(x)\cdot\vec{\eta} = \vec{h}_N(x) \quad
 &\text{para } x \in \Gamma_N \subset \partial \Omega.
\end{aligned}\right.
\label{ESD'}
\end{aligned}
\end{equation}
\label{Fbil}
Ahora definimos, las siguientes formas bilineales y funcionales lineales para la \textit{formulación débil} de la ecuación de Stokes como,
\begin{equation}
\mathcal{A}:[H_0^1(\Omega)]^2 \times [H_0^1(\Omega)]^2 \to \mathbb{R} \quad \text{como }\mathcal{A}(u,v)=\mu \int_\Omega \nabla \vec{u} \nabla\vec{v},
\label{FBA}
\end{equation}
\begin{equation}
\mathcal{B}:[H_0^1(\Omega)]^2 \times L_2(\Omega) \to \mathbb{R} \quad \text{como }\mathcal{B}(v,p)= - \int_\Omega \text{(div } \vec{v}) p, 
\label{FBB}
\end{equation} y
\begin{equation}
\mathcal{F}:[H_0^1(\Omega)]^2  \to \mathbb{R} \quad \text{como }\mathcal{F}(v)= \int_\Omega \vec{f}\vec{v}.
\label{FLF}
\end{equation}

Entonces la \textit{formulación débil} de la ecuación de Stokes \eqref{ESF} se escribirá como,
\begin{equation}
\begin{aligned}[b]
&\textit{Encontrar } \vec{u} \in [H^1_0(\Omega)]^2  \textit{ y } p \in L^2(\Omega) \text{ tal que:} \\
&\left\{
\begin{aligned}
\mathcal{A}(u,v) + \mathcal{B}(v,p) = \mathcal{F}(v)  \quad
 &\text{para toda } \vec{v} \in [H^1_0(\Omega)]^2,& \\
\mathcal{B}(u,q)  = 0 \quad
 &\text{para toda } q \in L^2(\Omega),& \\
\int_\Omega p = 0 \quad
 & \text{Condición de solubilidad} \\
 \vec{u}(x)=\vec{h}_D(x) \quad
 &\text{para }  x  \in \Gamma_D \subset \partial\Omega,	 \\
\nabla \vec{u}(x)\cdot\vec{\eta} = \vec{h}_N(x) \quad
 &\text{para } x \in \Gamma_N \subset \partial \Omega.
\end{aligned}\right.
\label{ESD}
\end{aligned}
\end{equation}

Entonces la ecuación \eqref{ESD} es la \textit{formulación débil} de la ecuación de Stokes. A este tipo de formulaciones, también se les llama \textit{formulaciones punto de silla}, y para garantizar la existencia de soluciones débiles para estas formulaciones, usamos la condición de Babuska-Brezzi, mostrada en la Sección \ref{Brezzi}.

Ahora, con la ecuación \eqref{ESD}, que es la \textit{formulación débil} de la ecuación de Stokes, debemos escoger un espacio de dimensión finita para obtener \textit{la formulación de Galerkin,} la cual veremos en la siguiente sección.

\section{Formulación de Galerkin y aproximación de elementos finitos}
\label{1.3}

Vamos a construir la \textit{formulación de Galerkin} de la Ecuación de Stokes, pero para esto, siguiendo el MEF, debemos construir una triangulación $\mathcal{T}^h$ que sea de aspecto regular y cuasi-uniforme en el dominio en el que estamos trabajando. En esta triangulación obtendremos los elementos (triángulos) $K_i$ con los vértices $x_i$, para crear un espacio de dimension finita.

Aquí es donde cada ecuación cambia el espacio en el que se va a resolver, tal que $u_1,u_2$ se buscaran en el espacio de funciones $\mathbb{P}^2(\mathcal{T}^h)$, donde usaremos funciones de prueba $v \in \mathbb{P}^2_0(\mathcal{T}^h)$, y buscaremos $p \in \mathbb{P}^1_0(\mathcal{T}^h)$, donde usaremos funciones de prueba $q \in \mathbb{P}^1_0(\mathcal{T}^h)$. Estos son espacios de dimension finita definidos como,
\begin{align*}
\mathbb{P}^1(\mathcal{T}^h)
&=\left\{ v \in C(\Omega) \quad : \quad v|_K \text{ es un polinomio de grado total 1 para todo } K \in \mathcal{T}^h \right\} \\
&=\left\{ v \in C(\Omega) \quad : \quad v|_K = a+bx+cy \right\}, \\
\mathbb{P}^1_0(\mathcal{T}^h)
&=\left\{ v \in \mathbb{P}^1(\mathcal{T}^h) \quad : \quad v(x) = 0 \quad \text{para todo } x \in  \partial \Omega \right\}, \\
\mathbb{P}^2(\mathcal{T}^h)
&=\left\{ v \in C(\Omega) \quad : \quad v|_K \text{ es un polinomio de grado total 2 para todo } K \in \mathcal{T}^h \right\}  \\
&=\left\{ v \in C(\Omega) \quad : \quad v|_K = a+bx+cy+dxy+ex^2+fy^2 \right\}, \text{ y } \\
\mathbb{P}^2_0(\mathcal{T}^h)
&=\left\{ v \in \mathbb{P}^2(\mathcal{T}^h) \quad : \quad v(x) = 0 \quad \text{para todo } x \in \partial \Omega \right\}. 
\end{align*}
Por lo que la \textit{formulación de Galerkin} se escribe como,
\begin{equation}
\begin{aligned}[b]
&\textit{Encontrar } \vec{u} \in [\mathbb{P}^2(\mathcal{T}^h)]^2  \textit{ y }p \in \mathbb{P}^1(\mathcal{T}^h) \text{ tal que:} \\
&\left\{
\begin{aligned}
\mathcal{A}(u,v) + \mathcal{B}(v,p) = \mathcal{F}(v)  \quad
 &\text{para toda } \vec{v} \in \mathbb{P}^2_0(\mathcal{T}^h),& \\
\mathcal{B}(u,q)  = 0 \quad
 &\text{para toda } q \in \mathbb{P}^1_0(\mathcal{T}^h),& \\
\int_\Omega p = 0 \quad
 & \text{Condición de solubilidad} \\
 \vec{u}(x)=\vec{h}_D(x) \quad
 &\text{para }  x  \in \Gamma_D \subset \partial\Omega,	 \\
\nabla \vec{u}(x)\cdot\vec{\eta} = \vec{h}_N(x) \quad
 &\text{para } x \in \Gamma_N \subset \partial \Omega.
\end{aligned}\right.
\label{ESG}
\end{aligned}
\end{equation} 

\begin{figure}[hrt]
  \centering
    \includegraphics[width=0.6\textwidth]{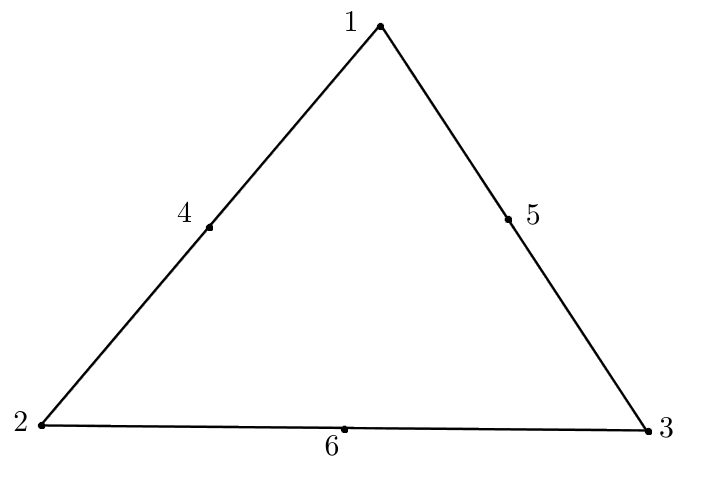}
  \caption{Nodos del triángulo $K \in \mathcal{T}^h$, que se usan para definir las funciones base de $\mathbb{P}^2(\mathcal{T}^h)$ }
  \label{TP2}
\end{figure}

Para el espacio de funciones $\mathbb{P}^2(\mathcal{T}^h)$, tenemos seis grados de libertad, por lo que debemos tomar seis nodos encada triángulo $K \in \mathbb{P}^2(\mathcal{T}^h)$, los cuales son los tres vértices del triángulo, y los tres puntos medios de cada lado del triángulo, como se ve en la Figura \ref{TP2}.

Entonces como función de prueba tomamos las funciones base de  que denotaremos $\varphi_i$ las cuales son de la forma,
$$\varphi_i(x)= \begin{cases}
1, & \text{Si } x=x_i \\
0, & \text{Si } x \neq x_i \\ 
\text{función cuadrática,} & \text{Caso contrario} 
\end{cases}$$
para los nodos $i=1,\ldots,6$. Para mas detalles de las funciones base de $\mathbb{P}^2(\mathcal{T}^h)$ véase la Sección \ref{EP2}, y véase \cite{Cla}.

Por lo que $\vec{u}=(u_1,u_2)$ se puede escribir como combinación lineal de estas funciones base, donde habrán seis para $u_1$, y seis para $u_2$, de modo que 
$$\vec{u}=\sum_{i=1}^{N_h^v}\boldsymbol{\vec{u}}(x_i)\vec{\varphi}_i$$
donde $N_h^v$ es el numero de vértices de la triangulación, y $\boldsymbol{\vec{u}}$ es el vector que representa las coordenadas de la función de elementos finitos, donde $\vec{\varphi_i}$ es de la forma
$$\vec{\varphi_i}=(\varphi_i,0) \quad \text{para } u_1,$$
$$\vec{\varphi}_{i+1}=(0,\varphi_i) \quad \text{para } u_2.$$
donde vemos que cada elemento (triángulo) $K \in \mathcal{T}^h$ posee doce funciones base en la \textit{formulación de Galerkin} de la ecuación de Stokes.

\noindent Por lo que escribimos la formulación matricial como,

\begin{equation}
\begin{aligned}[b]
&\textit{Encontrar } \boldsymbol{\vec{u}} \in [\mathbb{P}^2(\mathcal{T}^h)]^2  \textit{ y } \boldsymbol{p} \in \mathbb{P}^1(\mathcal{T}^h) \textit{ tal que:} \\
&\left\{
\begin{aligned}
A\boldsymbol{\vec{u}} + B\boldsymbol{p}= F \quad
 &\text{para toda } \vec{v} \in \mathbb{P}^2_0(\mathcal{T}^h),& \\
B\boldsymbol{\vec{u}}  = 0 \quad
 &\text{para toda } q \in \mathbb{P}^1_0(\mathcal{T}^h),& \\
\int_\Omega p = 0 \quad
 & \text{Condición de solubilidad} \\
 \vec{u}(x)=\vec{h}_D(x) \quad
 &\text{para }  x  \in \Gamma_D \subset \partial\Omega,	 \\
\nabla \vec{u}(x)\cdot\vec{\eta} = \vec{h}_N(x) \quad
 &\text{para } x \in \Gamma_N \subset \partial \Omega.
\end{aligned}\right.
\label{ESM}
\end{aligned}
\end{equation}
Donde,
$$A=[a_{ij}], \  B=[b_{ij}], \ F=[f_{i}], \ $$ y
$$a_{ij}=\mathcal{A}(\varphi_i,\varphi_j),$$
$$b_{ij}=\mathcal{B}(\varphi_i,\varphi_j),$$
$$ \ f_i=\mathcal{F}(\varphi_i)+ k_i.$$

Donde estos términos son las formas bilineales y los funcionales lineales, definidos en las ecuaciones \eqref{FBA},\eqref{FBB},\eqref{FLF}, en la página \pageref{Fbil} y $k_i$ es un término que surge de la condición de frontera.

Podemos escribir el sistema lineal de la ecuación \eqref{ESM}, de forma matricial como
\begin{equation}
\begin{bmatrix}
A & B^T \\
B & 0
\end{bmatrix}
\begin{bmatrix}
U \\ P
\end{bmatrix} =
\begin{bmatrix}
F \\ 0
\end{bmatrix}
\end{equation}

Con esto se conforma un sistema de ecuaciones que al resolverlo obtendremos la aproximación de la solución de la ecuación de Stokes, pero en una triangulación muy grande este sistema solo se podrá resolver por medio de maquinas con alta capacidad computacional.

\section{Ejemplos usando FreeFem++}
\label{FFSE}
Para resolver el la ecuación de Stokes, en un dominio $\Omega$, la abordaremos con un software especializado en el MEF, denominado FreeFem++. Para esto, escribiremos un código que proporcionara una aproximación de la solución de la ecuación. El código consta de varias partes, las cuales tienen su respectiva función en el desarrollo de la aproximación de la solución, que explicaremos a continuación.

\begin{figure}
  \centering
    \includegraphics[width=0.9\textwidth]{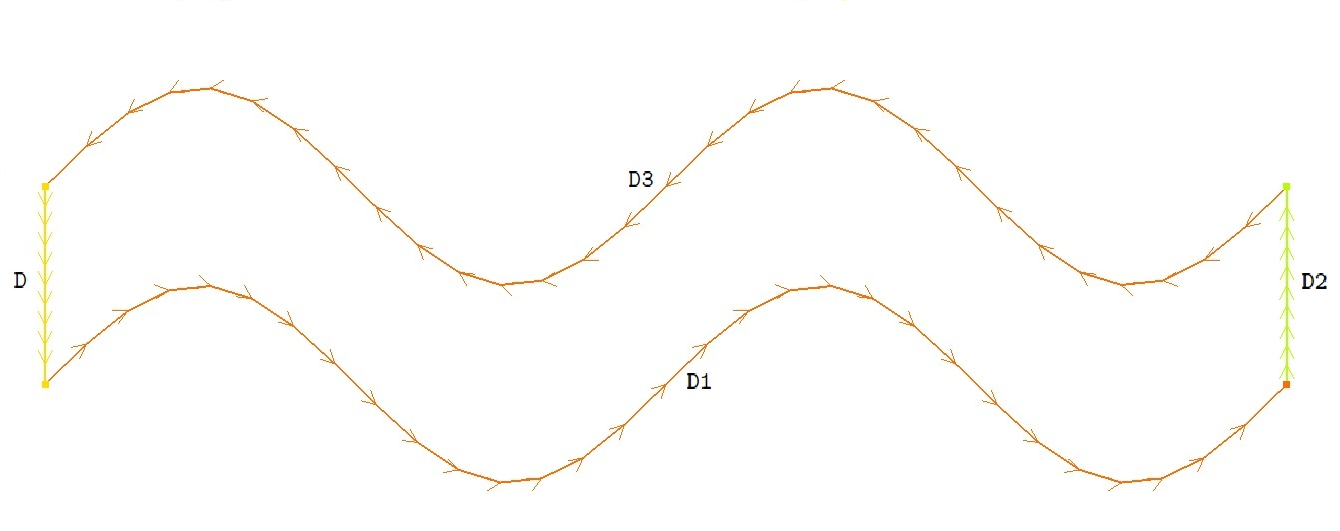}
  \caption{Dominio $\Omega$, obtenido parametrizando las curvas de la frontera $\partial \Omega$ (D,D1,D2,D3) en FreeFem++. Este dominio sera usado para resolver la ecuación de Stokes \eqref{ESF}, usando elementos finitos $\mathbb{P}^2$ para $\vec{u}=(u_1,u_2)$ y $\mathbb{P}^1$ para la presión $p$ }
  \label{Dom}
\end{figure}

Lo primero que haremos es crear un dominio $\Omega$, que sera el dominio en el que se buscara la aproximación del MEF para la ecuación de Stokes, que para este ejemplo buscaba que se asemejara a un dique. En FreeFem++ se define la frontera por pedazos, parametrizando cada pedazo de acuerdo a la forma que se quiera, lo cual se escribe como,

\begin{verbatim}
// Parametrización del Dominio

border D(t=1,-1){x=-2*pi; y=t; label=2;};
border D1(t=-2*pi,2*pi){x=t; y=sin(t)-1; label=1;};
border D2(t=-1,1){x=2*pi; y=t; label=3;};
border D3(t=2*pi,-2*pi){x=t; y=sin(t)+1; label=1;};
\end{verbatim}
donde \verb+D, D1, D2, D3+ serian los lados de la frontera. Podemos visualizar el dominio con el comando,

\begin{verbatim}
//Plotar el Dominio

plot(D(10)+D1(30)+D2(10)+D3(30));
\end{verbatim}
El programa genera una imagen que representa el dominio, en este ejemplo es la Figura \ref{Dom}.

\begin{figure}[hrt]
  \centering
    \includegraphics[width=0.9\textwidth]{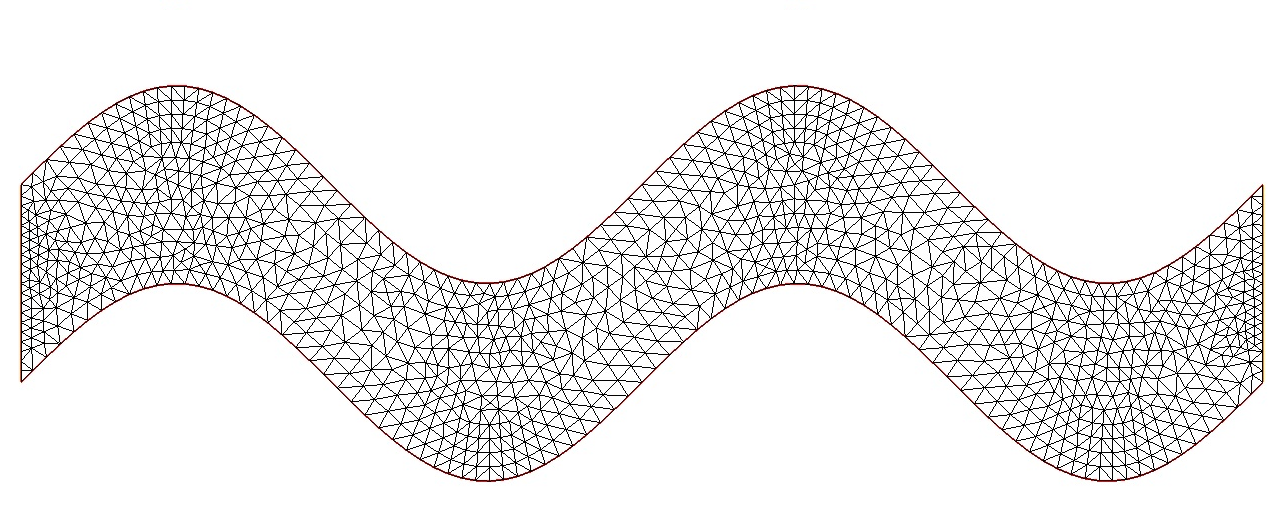}
  \caption{Triangulación $\mathcal{T}^h$ generada por FreeFem++, la cual contiene 2176 triángulos y 1199 vértices.}
  \label{Tri}
\end{figure}

Una vez definido el dominio, el paso a seguir con FreeFem++ es crear la triangulación, en la cual se va a trabajar, y así poder definir el espacio de elementos finitos de acuerdo a la dimension de la triangulación (según su número de elementos). La triangulación se construye con el comando

\begin{verbatim}
//Construcción de la triangulación

mesh Th= buildmesh (D(20)+D1(90)+D2(20)+D3(90));
\end{verbatim}
donde los números dentro de los paréntesis determinan el numero de vértices sobre la frontera, en este ejemplo en el borde \verb+D+ y \verb+D2+ habrán 20 vértices, mientras que en el borde \verb+D1+ y \verb+D3+ habrán 90 vértices, y de acá FreeFem++ creara la triangulación geométricamente conforme \footnote[1]{Una triangulación es llamada \textit{geométricamente conforme } si la intersección de dos elementos diferentes es un lado completo o un vértice común a los dos elementos.}, la cual según su análisis contiene 2176 triángulos y 1199 vértices.
\begin{verbatim}
--  mesh:  Nb of Triangles =   2176, Nb of Vertices 1199
\end{verbatim}
Esta triangulación se ilustra en la Figura \ref{Tri}. Para imprimirla usamos el comando,
\begin{verbatim}
//Plotar la triangulación

plot(Th,wait=1,ps="Triangulacion.eps");
\end{verbatim}

Una vez obtenida la triangulación, ahora debemos definir el espacio de elementos finitos  en el que vamos a trabajar, para el caso de la ecuación de Stokes, el espacio  $\mathbb{P}^2(\mathcal{T}^h)$ para el campo de velocidades, y usamos $\mathbb{P}^1(\mathcal{T}^h)$ para definir la presión $p$.
Con Freefem++, basta definir los espacios en la triangulación de la siguiente manera,

\begin{verbatim}
//Definición del Espacio de Elementos Finitos

fespace Vh(Th,P2);
fespace Mh(Th,P1);
\end{verbatim}
Como la ecuación de Stokes maneja 3 incógnitas, dos coordenadas para la velocidad y la presión, definimos $u_1, u_2 \in\mathbb{P}^2(\mathcal{T}^h) $ y $p\in\mathbb{P}^1(\mathcal{T}^h)$ que son $V_h,M_h$ respectivamente. En el código que estamos trabajando, luego de esto, definimos todos los términos que vamos a utilizar en cada espacio de elementos finitos de la siguiente manera,

\begin{verbatim}
//Declaración de funciones de elementos finitos de prueba y de forma

Vh u1,v1,u2,v2;
Mh p,q;
\end{verbatim}
es decir, en $V_h=\mathbb{P}^2(\mathcal{T}^h)$ definiremos $\vec{u}=(u_1,u_2)$, $\vec{v}=(v_1,v_2)$ que son velocidades solución y de prueba, y en $M_h=\mathbb{P}^1(\mathcal{T}^h)$ definimos $p,q$ que son presión solución y de prueba.

Luego escribimos la \textit{formulación de Galerkin} de la Ecuación de Stokes, en donde asumiremos $\mu$ constante $\mu = 0.1$, en FreeFem++ se escribe,

\begin{verbatim}
func mus=0.1; //Coeficiente de Viscosidad

//Formulación del problema de Stokes

solve Stokes (u1,u2,p,v1,v2,q,solver=Crout) =  
      int2d(Th)( mus*(dx(u1)*dx(v1)+dy(u1)*dy(v1)+dx(u2)*dx(v2)+dy(u2)*dy(v2) ) 
	              - p*q*(0.00000001) - p*dx(v1) - p*dy(v2) - dx(u1)*q -dy(u2)*q )
				  +on(1,u1=0,u2=0) +on(2,u1=-1.5*(y-1)*(y+1),u2=1);
\end{verbatim}
donde \verb+int2d(Th)+ es la integral sobre el dominio omega de las siguientes expresiones,
$$\text{mus*(dx(u1)*dx(v1)+dy(u1)*dy(v1)+dx(u2)*dx(v2)+dy(u2)*dy(v2) )  }= \mu \nabla \vec{u} \nabla \vec{v},$$
$$ \text{- p*dx(v1) - p*dy(v2)} = -\text{(div }\vec{v})p, $$
$$ \text{- dx(u1)*q - dy(u2)*q} = -\text{(div }\vec{u})q$$ 
y \verb+p*q*(0.00000001)+ se refiere a un término de estabilización, para fijar la parte constante de la presión (condición de solubilidad). Para mas detalles, véase \cite{FreeFem}.

La ultima linea define la condición de frontera, que en nuestro ejemplo las condiciones de frontera en \verb+D1+ y \verb+D3+ son Dirichlet cero, es decir no entra ni sale nada del fluido; en \verb+D+ hay una función de entrada de fluido, la cual hace que se perciba una especie de recorrido del fluido en el dominio; y en \verb+D2+ hay condición de frontera Neumann cero, la cual muestra que el fluido no tiene ningún tipo de obstrucción, que de ahí en adelante sigue su recorrido regular.

\begin{figure}[hrt]
  \centering
    \includegraphics[width=0.9\textwidth]{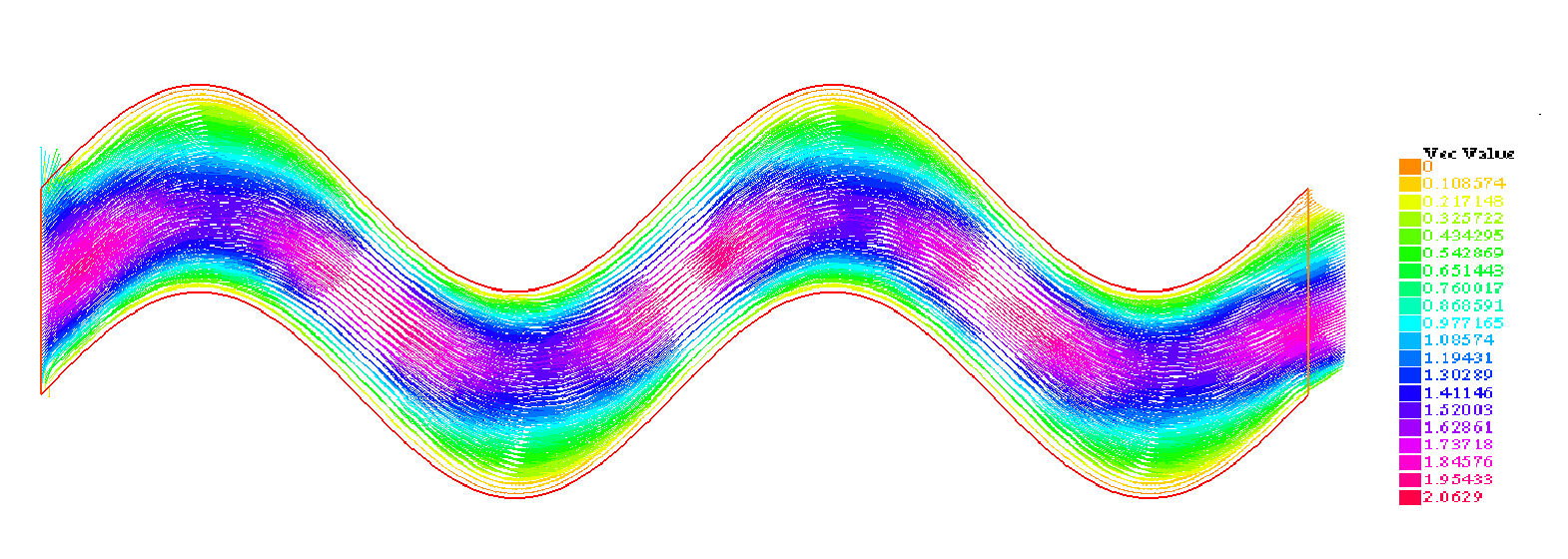}
  \caption{Aproximación por MEF del campo de velocidades $\vec{u}=(u_1,u_2)$ de la ecuación de Stokes \eqref{ESF}, con coeficiente de viscosidad $\mu=0.1.$}
  \label{Stokes}
\end{figure}

Con esto basta para encontrar la aproximación por MEF de la solución de la ecuación de Stokes, para mostrarla, imprimimos la solución $\vec{u}=(u_1,u_2)$ con 

\begin{verbatim}
// Plotar la velocidad u=(u_1,u_2)

plot([u1,u2],wait=1,dim=2,fill=1,value=1);
\end{verbatim}
que se muestra en la Figura \ref{Stokes}, la cual muestra el campo de velocidades de la ecuación de Stokes, que aproxima como se mueve un fluido con un coeficiente de viscosidad $\mu$ que en este ejemplo es $\mu=0.1$.

\begin{figure}[H]
  \centering
    \includegraphics[width=0.9\textwidth]{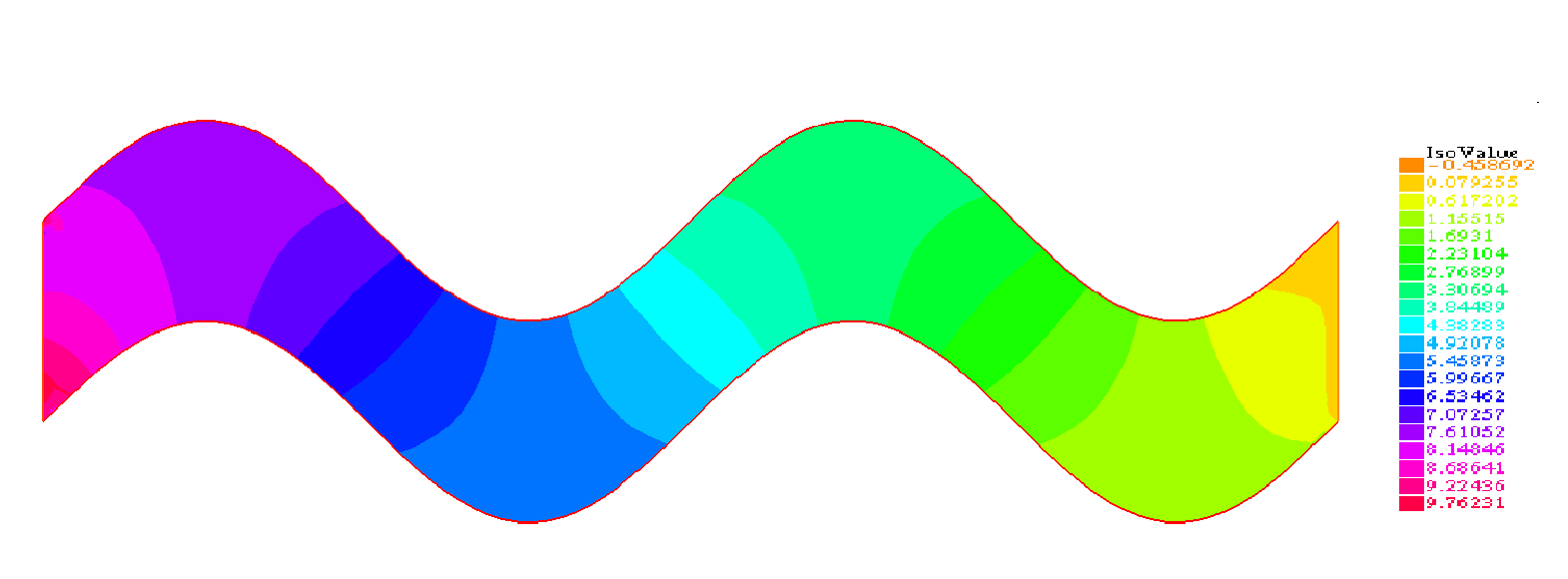}
  \caption{Aproximación por MEF de la presión $p$ de la ecuación de Stokes}
  \label{Pres}
\end{figure}

\noindent Por último, usamos el comando
\begin{verbatim}
//Plotar la presión p

plot(p,wait=1,dim=2,fill=1,value=1);
\end{verbatim}
para imprimir la presión $p$ en el dominio $\Omega$, que se mostrara en la Figura \ref{Pres}.

Esta solución mostraría como se mueve un fluido en un dique con la forma del dominio $\Omega$. Con este ejemplo trabajaremos mas adelante para mostrar como se usa esta ecuación en diversos problemas de la física y de la vida cotidiana.

\chapter{Ecuación de advección-difusión}

En este capitulo mostrare el desarrollo del método de elementos finitos  para la ecuación de advección-difusión, la cual describe el fenómeno por el cual una sustancia se transforma, o en otras palabras, se mueve debido a los procesos de advección y de difusión. Para mas detalles sobre esta ecuación, y su aproximación por el método de elementos finitos se puede consultar \cite{Cla,FreeFem}.
\section{Introducción}

La ecuación de advección-difusión, es una EDP de tipo parabólico, y describe el proceso de la advección y la difusión de una sustancia en un fluido, donde la advección es una especie de transporte de la sustancia, es decir, la sustancia es transportada por el efecto de un campo vectorial, generado por una energía, por ejemplo, el transporte de una sustancia en un río, y la difusión es que tanto se propaga la sustancia en el medio que la rodea.

Esta ecuación aparece en el estudio de la propagación de sustancias contaminantes en diferentes medios, su solución representa la concentración de la sustancia que se dispersa y es una información muy importante en el estudio de diversos problemas de diferente naturaleza (física, química, económica, ambiental, etc.). Debido a la complejidad de estos problemas no es posible en general obtener una solución analítica para la ecuación de advección-difusión, por lo que usamos el MEF para hallar una solución aproximada.

La ecuación de advección-difusión esta dada por, 
\begin{equation}
\begin{aligned}[b]
&\textit{Encontrar } u : \Omega \subset \mathbb{R}^2 \to \mathbb{R} \textit{ tal que:} \\
&\left\{
\begin{aligned}
\dfrac{\partial u}{\partial t}-\mu \Delta u + \beta_1 \dfrac{\partial u}{\partial x_1} + \beta_2 \dfrac{\partial u}{\partial x_2}  = f \quad
&\text{para } x \in \Omega,& \\
u(x,t)=g_D(x) \quad
 &\text{para }  x  \in \Gamma_D \subset \partial\Omega, \ t>0 \\
\nabla u(x,t)\eta(x) = g_N(x) \quad
 &\text{para } x \in \Gamma_N \subset \partial \Omega, \ t>0 \\
u(0,x)=u_0(x) \quad
& \text{para } x \in \Omega,
\end{aligned}\right.
\label{EADF}
\end{aligned}
\end{equation}
donde $x=(x_1,x_2)$, $\frac{\partial u}{\partial t}$ es el diferencial a través del tiempo, $\Delta u$ es el termino correspondiente a la difusión, con coeficiente de difusión $\mu>0 $ y las primeras derivadas parciales corresponden a la advección en la dirección $\beta=(\beta_1,\beta_2)$.

Para la aplicación del MEF a la ecuación de advección-difusión, denominaremos \eqref{EADF} como la \textit{formulación fuerte} de  esta, hallaremos su respectiva \textit{formulación débil}, y así seguiremos con el método para llegar a la solución de la ecuación. Podemos ver que esta ecuación depende de el tiempo y del espacio.

\section{Formulación débil de la ecuación de advección-difusión}
\label{FDAD}

Para hallar la \textit{formulación débil} de la ecuación de advección-difusión, comenzaremos suponiendo que $u$ es una solución a la ecuación \eqref{EADF}, tomamos $v \in H^1_0(\Omega)$ como función de prueba y la multiplicamos por el primer término de \eqref{EADF}, pero antes sabemos que $\beta=(\beta_1,\beta_2)$ y $\nabla u = (\frac{\partial u}{\partial x_1},\frac{\partial u}{\partial x_2})$ entonces,
$$\beta_1 \dfrac{\partial u}{\partial x_1} + \beta_2 \dfrac{\partial u}{\partial x_2}= \beta\nabla u.$$
Por lo tanto la ecuación de advección-difusión se escribe como,
$$\dfrac{\partial u}{\partial t}-\mu \Delta u + \beta\nabla u = f$$
al multiplicarla por $v$ obtenemos,
$$\dfrac{\partial u}{\partial t}v-\mu \Delta u v + \beta\nabla u v = fv.$$
Integramos sobre el dominio $\Omega \subset \mathbb{R}^2$, 
$$\int_\Omega\dfrac{\partial u}{\partial t}v-\mu \int_\Omega \Delta u v + \int_\Omega \beta \nabla u v = \int_\Omega fv.$$
Aplicamos la fórmula de green y obtenemos,
$$\int_\Omega\dfrac{\partial u}{\partial t}v+\mu \int_\Omega \nabla u \nabla v - \cancelto{0}{\int_{\partial\Omega}v\nabla u \cdot \eta} + \int_\Omega \beta\nabla u v = \int_\Omega fv$$
de aquí,
$$\int_\Omega\dfrac{\partial u}{\partial t}v+\mu \int_\Omega \nabla u \nabla v + \int_\Omega \beta\nabla u v  = \int_\Omega fv.$$
Note que usamos el hecho que $v=0$ en $\partial \Omega.$

Con lo anterior, escribimos la formulación débil de la ecuación de advección-difusión como,

\begin{equation}
\begin{aligned}[b]
&\textit{Encontrar } u \in H_0^1(\Omega) \textit{ tal que:} \\
&\left\{
\begin{aligned}
\int_\Omega\dfrac{\partial u}{\partial t}v+\mu \int_\Omega \nabla u \nabla v + \int_\Omega \beta\nabla u v  = \int_\Omega fv \quad
&\text{para toda } v \in H_0^1(\Omega),& \\
u(x,t)=g_D(x) \quad
 &\text{para }  x  \in \Gamma_D \subset \partial\Omega, \ t>0 \\
\nabla u(x,t)\eta(x) = g_N(x) \quad
 &\text{para } x \in \Gamma_N \subset \partial \Omega, \ t>0 \\
u(0,x)=u_0(x) \quad
& \text{para } x \in \Omega
\end{aligned}\right.
\label{EADD'}
\end{aligned}
\end{equation}

%%%Ya con esto, para facilitar la escritura de \ref{EADD'} definimos las siguientes formas bilineales y funcionales lineales como
%%%$$\mathcal{T}(u,v)=\int_\Omega\dfrac{\partial u}{\partial t}v$$
%%%$$\mathcal{A}(u,v)=\mu \int_\Omega \nabla u \nabla v$$
%%%$$\mathcal{M}(u,v)= \int_\Omega\beta \nabla u v $$
%%%$$\mathcal{L}(u,v)= \int_\Omega uv $$
%%%$$\mathcal{F}(v)= \int_\Omega fv$$ 
%%%entonces la \textit{formulación débil} se escribirá como,
%%%
%%%\begin{equation}
%%%\begin{aligned}[b]
%%%&\textit{Encontrar } u \in H^0_1(\Omega)  \textit{ tal que:} \\
%%%&\left\{
%%%\begin{aligned}
%%%\mathcal{T}(u,v)+\mathcal{A}(u,v) + \mathcal{M}(u,v) + \mathcal{L}(u,v) = \mathcal{F}(v)  \quad
%%% &\text{para toda } v \in H_0^1(\Omega),& \\
%%%\text{Condicion de frontera} \quad
%%%&\text{para } x=(x_1,x_2) \in \Gamma \subset \partial\Omega,&	 \\
%%%\text{Condición inicial} \quad
%%%&u=u_0 \text{ en la frontera }
%%%\end{aligned}\right.
%%%\label{EADD}
%%%\end{aligned}
%%%\end{equation}

Entonces la ecuación \eqref{EADD'} es la \textit{formulación débil} de la ecuación de advección-difusión. Para garantizar la existencia de una solución débil para esta formulación, usamos el teorema de Lax-Milgram, enunciado en la Sección \ref{LM}. 

Ahora, con la ecuación \eqref{EADD'}, que es la \textit{formulación débil} de la ecuación de advección-difusión, debemos elegir un espacio de dimensión finita adecuado, para construir la \textit{formulación de Galerkin}, la cual es el paso a seguir en el MEF.

\section{Formulación de Galerkin y aproximación de elementos finitos}

La \textit{formulación de Galerkin} es un método para discretizar el espacio $\Omega$ en el que estamos trabajando, pero el tiempo sigue siendo continuo. Luego de realizar esta discretización del espacio, veremos un método para discretizar el tiempo.

Sabemos que para esta formulación, lo primero que necesitamos es construir una triangulación $\mathcal{T}^h$, con aspecto regular y cuasi-uniforme en el dominio $\Omega$, con esta triangulación obtenemos los elementos $K_i$ y los vértices $x_i$, que crearan un espacio de dimension finita. Para la ecuación de advección difusión usamos el espacio de funciones $\mathbb{P}^2_0(\mathcal{T}^h)$ (definido en la Sección \ref{1.3}),
\begin{align*}
\mathbb{P}^2(\mathcal{T}^h)
&=\left\{ v \in C(a,b) \quad : \quad v|_K \text{ es un polinomio de grado total 2 para todo } K \in \mathcal{T}^h \right\}  \\
&=\left\{ v \in C(a,b) \quad : \quad v|_K = a+bx+cy+dxy+ex^2+fy^2 \right\}, \text{ y } \\
\mathbb{P}^2_0(\mathcal{T}^h)
&=\left\{ v \in \mathbb{P}^2(\mathcal{T}^h) \quad : \quad v(x) = 0 \quad \text{para todo } x \in \partial \Omega \right\}. 
\end{align*}
es decir este sera el espacio de funciones de forma y funciones de prueba.

Por lo que la \textit{formulación de Galerkin} se escribe como,
\begin{equation}
\begin{aligned}[b]
&\textit{Encontrar } u \in \mathbb{P}^2_0(\mathcal{T}^h)  \textit{ tal que:} \\
&\left\{
\begin{aligned}
\int_\Omega\dfrac{\partial u}{\partial t}v+\mu \int_\Omega \nabla u \nabla v + \int_\Omega \beta\nabla u v  = \int_\Omega fv  \quad
&\text{para toda } v \in \mathbb{P}^2_0(\mathcal{T}^h),& \\
u(x,t)=g_D(x) \quad
 &\text{para }  x  \in \Gamma_D \subset \partial\Omega, \ t>0 \\
\nabla u(x,t)\eta(x) = g_N(x) \quad
 &\text{para } x \in \Gamma_N \subset \partial \Omega, \ t>0 \\
u(0,x)=u_0(x) \quad
& \text{para } x \in \Omega.
\end{aligned}\right.
\label{EADG}
\end{aligned}
\end{equation}
Con esto tenemos la \textit{formulación de Galerkin} semidiscreta, ya que el tiempo aun esta continuo, para obtener la formulación discreta usaremos un método para discretizar el tiempo.

Ahora, escribiremos la formulación matricial semidiscreta (tiempo continuo), por lo que usaremos las funciones base $\varphi_i$ de $\mathbb{P}^2(\mathcal{T}^h)$, tal que $u$ es combinación lineal de estas $\varphi_i$.

$$u=\sum_{i=1}^{N_h^v}\boldsymbol{\alpha}(x_i)\varphi_i$$
donde $N_h^v$ es el numero de vértices de la triangulación, y $\boldsymbol{u}$ es el vector que representa las coordenadas de la función de elementos finitos. 

Por lo que escribimos la \textit{formulación matricial } semi-discreta (tiempo continuo) como,

\begin{equation}
\begin{aligned}[b]
&\textit{Encontrar } \alpha \in \mathbb{P}^2_0(\mathcal{T}^h)  \textit{ tal que:} \\
&\left\{
\begin{aligned}
(\boldsymbol{T} + \boldsymbol{A} + \boldsymbol{V})\alpha=\boldsymbol{b} \quad
&\text{para toda } v \in \mathbb{P}^2_0(\mathcal{T}^h),& \\
u(x,t)=g_D(x) \quad
 &\text{para }  x  \in \Gamma_D \subset \partial\Omega, \ t>0 \\
\nabla u(x,t)\eta(x) = g_N(x) \quad
 &\text{para } x \in \Gamma_N \subset \partial \Omega, \ t>0 \\
u(0,x)=u_0(x) \quad
& \text{para } x \in \Omega.
\end{aligned}\right.
\label{EADM'}
\end{aligned}
\end{equation}
donde,
$$ \boldsymbol{T}=[t_{ij}], \ \boldsymbol{A}=[a_{ij}], \ \boldsymbol{V}=[v_{ij}], \ \boldsymbol{b}=[b_{i}], $$ y
$$t_{ij}=\int_\Omega\dfrac{\partial \varphi_i}{\partial t}\varphi_j $$
$$a_{ij}=\mu \int_\Omega \nabla \varphi_i \nabla \varphi_j$$
$$v_{ij}=\int_\Omega \beta\nabla \varphi_i \varphi_j $$
$$ \ b_i=\int_\Omega f\varphi_i + l_i $$
son las matrices formadas por las funciones base de $\mathbb{P}^2_0(\mathcal{T}^h)$ y $l_i$ es un término que surge de acuerdo a la condición de frontera.

Ahora, esta formulación \eqref{EADM'} no es posible solucionarla en un software para el desarrollo del MEF, ya que el tiempo esta continuo, por lo que estudiaremos un método para discretizar el tiempo y así poder introducirla en un software como FreeFem++, para encontrar una aproximación de la solución.

\subsection{Discretización del tiempo}

Para obtener la discretización del tiempo, escribimos la serie de Taylor de $u$ como,
$$u(t+\Delta t,x)=u(t,x)+\Delta t \dfrac{\partial u}{\partial t}(t,x)+\Delta t^2r(x,t),$$
despejando $\dfrac{\partial u}{\partial t}$ y asumiendo que $r(x,t)$ es acotado, decimos que
$$\dfrac{\partial u}{\partial t}(t,x)\approx\dfrac{u(t+\Delta t,x)-u(t,x)}{\Delta t}$$
con un error de aproximación de orden $\Delta t$.
Para simplificar la escritura, denotamos esta expresión como,
\begin{equation}
\dfrac{\partial u}{\partial t}(t,x)\approx\dfrac{u(t+\Delta t,x)-u(t,x)}{\Delta t}=\dfrac{u_m - u_{m-1}}{\Delta t}
\label{MEG}
\end{equation}
donde $t_m=t+\Delta t$, $t_{m-1}=t$, $u_m=u(t+\Delta t,x)$ y $u_{m-1}=u(t,x).$

Con esto, aplicamos el método de Euler-Galerkin, el cual consiste en reemplazar \eqref{MEG} en \eqref{EADG} y  así obtenemos
\begin{equation}
\int_\Omega\dfrac{u_m - u_{m-1}}{\Delta t}v+\mu \int_\Omega \nabla u \nabla v + \int_\Omega \beta\nabla u v = \int_\Omega fv. 
\label{EADG'}
\end{equation}
Ahora, como $u_{m-1}=u(t,x)$ es decir en el tiempo $t$, despejamos la ecuación \eqref{EADG'} de la siguiente manera,
$$\dfrac{1}{\Delta t}\int_\Omega u_mv + \mu \int_\Omega \nabla u_{m} \nabla v + \int_\Omega\beta \nabla u_{m} v   =  \int_\Omega fv + \dfrac{1}{\Delta t}\int_\Omega u_{m-1}v $$
tenemos que,
\begin{equation}
\begin{aligned}[b]
&\textit{Dado } u_0 \in \mathbb{P}^2_0(\mathcal{T}^h) \textit{ encontrar } u_1,u_2,\ldots,u_m,\ldots  \textit{ tal que:} \\
&\left\{
\begin{aligned}
\frac{1}{\Delta t}{ \int_\Omega} u_mv + \mu \int_\Omega \nabla u_{m} \nabla v+ \int_\Omega \beta\nabla u_{m} v =    \quad
&\text{para toda } v \in \mathbb{P}^2_0(\mathcal{T}^h),& \\ 
  \int_\Omega fv + \frac{1}{\Delta t}\int_\Omega u_{m-1}v \quad & \\
u(x,t)=g_D(x) \quad
 &\text{para }  x  \in \Gamma_D \subset \partial\Omega, \ t>0 \\
\nabla u(x,t)\eta(x) = g_N(x) \quad
 &\text{para } x \in \Gamma_N \subset \partial \Omega, \ t>0 \\
u(0,x)=u_0(x) \quad
& \text{para } x \in \Omega
\end{aligned}\right.
\label{EADGD}
\end{aligned}
\end{equation} 
es la \textit{formulación de Galerkin} discreta y implícita, así ya podemos hacer la aproximación por elementos finitos.

Ahora tomamos las funciones base $\varphi_i$ de $\mathbb{P}^2(\mathcal{T}^h)$, tal que $u$ es combinación lineal de estas $\varphi_i$, es decir,

$$u=\sum_{i=1}^{N_h^v}\boldsymbol{u}(x_i)\varphi_i$$
donde $N_h^v$ es el numero de vértices de la triangulación, y $\boldsymbol{u}$ es el vector que representa las coordenadas de la función de elementos finitos. Ahora, definimos las formas bilineales y funcional lineal como, 

\begin{equation}
\mathcal{A}:H_0^1(\Omega) \times H_0^1(\Omega) \to \mathbb{R} \quad \text{como }\mathcal{A}(u,v)=\mu \int_\Omega \nabla u \nabla v,
\label{FBA1}
\end{equation}
\begin{equation}
\mathcal{V}:H_0^1(\Omega) \times H_0^1(\Omega) \to \mathbb{R} \quad \text{como }\mathcal{V}(u,v)= \int_\Omega\beta \nabla u v,
\label{FBV}
\end{equation}
\begin{equation}
\mathcal{M}:H_0^1(\Omega) \times H_0^1(\Omega) \to \mathbb{R} \quad \text{como }\mathcal{M}(u,v)= \int_\Omega uv,
\label{FBM}
\end{equation}
y
\begin{equation}
\mathcal{F}:H_0^1(\Omega) \to \mathbb{R} \quad \text{como }\mathcal{F}(v)= \int_\Omega fv.
\label{FLF1}
\end{equation}
Esto se hace para que sea mas cómoda la escritura. Entonces podemos escribir la \textit{formulación de Galerkin} discreta como sigue, 
\begin{equation}
\begin{aligned}[b]
&\textit{Dado } u_0 \in \mathbb{P}^2_0(\mathcal{T}^h) \textit{ encontrar } u_1,u_2,\ldots,u_m,\ldots  \textit{ tal que:} \\
&\left\{
\begin{aligned}
\textstyle{\frac{1}{\Delta t}\mathcal{M}(u_m,v) + \mathcal{A}(u_{m},v)- \mathcal{V}(u_{m},v) =  \mathcal{F}(v) + \frac{1}{\Delta t}\mathcal{M}(u_{m-1},v)} \ 
&\text{para toda } v \in \mathbb{P}^2_0(\mathcal{T}^h),& \\ 
u(x,t)=g_D(x) \
 &\text{para }  x  \in \Gamma_D \subset \partial\Omega,  t>0 \\
\nabla u(x,t)\eta(x) = g_N(x)  \
 &\text{para } x \in \Gamma_N \subset \partial \Omega,  t>0 \\
u(0,x)=u_0(x) \
& \text{para } x \in \Omega.
\end{aligned}\right.
\end{aligned}
\end{equation} 
Para escribir la formulación matricial implícita como
\begin{equation}
\begin{aligned}[b]
&\textit{Dado } \boldsymbol{\vec{u}_0} \in \mathbb{R}^2 \textit{ encontrar } \boldsymbol{\vec{u}_1},\boldsymbol{\vec{u}_2},\ldots,\boldsymbol{\vec{u}_m},\ldots \textit{ tal que:} \\
&\left\{
\begin{aligned}
\frac{1}{\Delta t}M\boldsymbol{\vec{u}_m} + A\boldsymbol{\vec{u}_{m}} + V\boldsymbol{\vec{u}_{m}} =  b  + \frac{1}{\Delta t}M\boldsymbol{\vec{u}_{m-1}}  \quad
 &\text{para toda } \vec{v} \in \mathbb{P}^2_0(\mathcal{T}^h),& \\
u(x,t)=g_D(x) \quad
 &\text{para }  x  \in \Gamma_D \subset \partial\Omega, \ t>0 \\
\nabla u(x,t)\eta(x) = g_N(x) \quad
 &\text{para } x \in \Gamma_N \subset \partial \Omega, \ t>0 \\
\vec{u}(0,x)=\vec{u}_0(x) \quad
& \text{para } x \in \Omega.
\end{aligned}\right.
\label{EADM}
\end{aligned}
\end{equation}
donde,
$$M=[m_{ij}], \ A=[a_{ij}], \ V=[v_{ij}], \ b=[b_{i}], $$ y tenemos que
$$m_{ij}=\mathcal{M}(\varphi_i,\varphi_j),$$
$$a_{ij}=\mathcal{A}(\varphi_i,\varphi_j),$$
$$v_{ij}=\mathcal{V}(\varphi_i,\varphi_j),$$ y
$$ \ b_i=\mathcal{F}(\varphi_i)+ k_i.$$
Donde estos términos son las formas bilineales y los funcionales lineales definidos anteriormente en las ecuaciones \eqref{FBM},\eqref{FBA1},\eqref{FBV},\eqref{FLF1} y $k_i$ es un término que surge de acuerdo a la condición de frontera.

La ecuación \eqref{EADM} es la \textit{formulación matricial} de la ecuación de advección-difusión, para hallar una aproximación por elementos finitos, basta con resolver este sistema, el cual en una triangulación muy fina, tendrá un alto costo computacional.

Para esto diseñamos algunos códigos en distintos programas para resolver esta ecuación que veremos en lo que resta de este capítulo.

\section{Ejemplos en una dimensión usando MatLab}

Escribimos un código en Matlab para resolver la ecuación de advección-difusión en una dimensión, en donde usamos la formulación matricial de esta, armando las matrices localmente, y luego ensamblando las matrices para obtener la aproximación. En Matlab, creamos estas matrices con variables tipo \textit{estructura}, las cuales asignan a cada elemento de la partición $K_i$ toda la información necesaria para construir la matriz.

Debemos saber que para este código trabajamos, en el intervalo $[0,1]$, con el espacio de funciones lineales por partes $\mathbb{P}^1_0(\mathcal{T}^h)$, en vez de $\mathbb{P}^2_0(\mathcal{T}^h)$,y que vamos a usar la condición inicial $u_0=0$. Así que lo primero que hacemos es definir una serie de elementos que vamos a usar a lo largo del código, entonces creamos las siguientes variables
\begin{verbatim}
%Parámetros

N=200;     %Numero de elementos
T=5;       %Tiempo final
dt=0.1;    %Incremento de tiempo 
u_0=zeros(N+1,1); %Vector de ceros, representando la condición inicial

%Formula de Cuadratura

omega_aux(1)=5/9; %Pesos de cuadratura
omega_aux(2)=8/9;
omega_aux(3)=5/9;
p_aux(1)=-sqrt(15)/5; %Puntos de cuadratura en [-1,1]
p_aux(2)=0;
p_aux(3)=+sqrt(15)/5;
\end{verbatim}
donde \verb+N+ es el numero de elementos de la partición, \verb+dt+ es el intervalo de tiempo, $u_0$ es la condición inicial, y los términos terminados en \verb+aux+, son usados para la integración numérica, por medio de la cuadratura de Gauss, la cual consiste en dar una aproximación de una integral definida, de una función dada. Una cuadratura n de Gauss, selecciona unos puntos de evaluación $x_i$ llamados puntos de cuadratura y unos coeficientes o pesos $\omega_i$ para $i=1,...,n$. El dominio de tal cuadratura es $[-1, 1]$ y esta dada por,

$$\displaystyle\int_{-1}^{1} f(x) \, dx\approx \sum_{i=1}^n \omega_i f(x_i),$$
pero como nuestras integrales deben ser calculadas en un dominio $[a,b]=[0,1]$ diferente a $[-1,1]$ entonces debemos hacer el cambio,

$$\displaystyle\int_{a}^{b} f(x) \, dx=\frac{b-a}{2}\displaystyle \int_{-1}^{1} f\left( \frac{b-a}{2}x+\frac{a+b}{2}\right) \, dx,$$
y al aplicar cuadratura obtenemos,

$$\displaystyle\int_{a}^{b} f(x) \, dx\approx \frac{b-a}{2} \sum_{i=1}^n \omega_i f\left( \frac{b-a}{2}x_i+\frac{a+b}{2}\right),$$
y esta es la aproximación que usamos en el código para calcular las integrales. Los detalles de esta cuadratura se pueden revisar en \cite{Rojo,Verde,Gauss}.

A continuación, se calcula la información correspondiente a los elementos finitos y las funciones base de elementos finitos, así como los datos necesarios para calcular las integrales locales que aparecen en la \textit{formulación de Galerkin} de la ecuación \eqref{EADGD}. Tenemos el código,	
A continuación, creamos las matrices locales, que en Matlab se le asigna elemento por elemento la información necesaria para calcular las integrales de  A estas variables se les llama estructuras, las cuales asignan a cada elemento $K_i$ una parte de la información de la matriz. Con el siguiente código 
\begin{verbatim}
%Elementos

for i=1:N
    xini=(i-1)/N;
    xfin=i/N;
    E(i).xini=xini;
    E(i).xfin=xfin;
    omega=0.5*(xfin-xini)*omega_aux;
    x_q=xini+ 0.5*(1+p_aux)*(xfin-xini);
    
    E(i).omega=omega;
    E(i).x_q=x_q;
    E(i).b1= (xfin-x_q)./(xfin-xini);
    E(i).b2= (x_q-xini)./(xfin-xini);
    E(i).b1d= [-1,-1,-1]./(xfin-xini);
    E(i).b2d=[1,1,1]./(xfin-xini);
    E(i).dof=[i,i+1];
end
\end{verbatim}
Ahora, creamos las matrices locales de M,A,V y el vector b, que en Matlab se le asigna elemento por elemento la información necesaria. Para esto usamos el siguiente código,
\begin{verbatim}
%%% Calculo de matrices locales

for i=1:N;
    b1=E(i).b1;
    b2=E(i).b2;
    b1d=E(i).b1d;
    b2d=E(i).b2d;
    x_q=E(i).x_q;
    w=E(i).omega;
    coef=0.05; %% mu(x)
    coefv=0.5; %% beta(x)
      
    %Calculo matriz local A   
    a11= sum(coef.*w.*b1d.*b1d);
    a12= sum(coef.*w.*b1d.*b2d);
    a21= sum(coef.*w.*b2d.*b1d);
    a22= sum(coef.*w.*b2d.*b2d);
    E(i).Alocal=[a11,a12;a21,a22];

    %Calculo matriz local M    
    m11= sum(w.*b1.*b1);
    m12= sum(w.*b1.*b2);
    m21= sum(w.*b2.*b1);
    m22= sum(w.*b2.*b2);
    E(i).Mlocal= [m11,m12;m21,m22];
    
    %Calculo matriz local V
    v11= sum(coefv.*w.*b1.*b1d);
    v12= sum(coefv.*w.*b1.*b2d);
    v21= sum(coefv.*w.*b2.*b1d);
    v22= sum(coefv.*w.*b2.*b2d);
    E(i).Vlocal=[v11,v12;v21,v22];
    
    %Calculo vector b
    f=0;
    f1=sum(w.*f.*b1);
    f2=sum(w.*f.*b2);
    E(i).flocal=[f1,f2];
end
\end{verbatim}
donde escribimos la integral para calcularla numéricamente localmente por el método de la cuadratura de Gauss. Para este ejemplo tomamos los valores de parámetros $\mu=0.05,\beta=0.5$ y $f=0$. Ahora para crear las matrices globales, creamos matrices dispersas nulas, y por cada elemento vamos añadiendo a cada matriz su contribución, que se ve en las siguientes lineas,
\begin{verbatim}
%%% Matriz Global
A=sparse(N+1,N+1);
M=sparse(N+1,N+1);
V=sparse(N+1,N+1);
b=sparse(N+1,1);

for k=1:N;
    Alocal=E(k).Alocal;
    flocal=E(k).flocal;
    Mlocal=E(k).Mlocal;
    Vlocal=E(k).Vlocal;
    dof=E(k).dof;
    i1=dof(1);
    i2=dof(2);

    %Contribución del elemento A(i)
    A(i1,i1)=A(i1,i1)+Alocal(1,1);
    A(i1,i2)=A(i1,i2)+Alocal(1,2); 
    A(i2,i1)=A(i2,i1)+Alocal(2,1);
    A(i2,i2)=A(i2,i2)+Alocal(2,2);

    %Contribución del elemento M(i)
    M(i1,i1)=M(i1,i1)+Mlocal(1,1);
    M(i1,i2)=M(i1,i2)+Mlocal(1,2); 
    M(i2,i1)=M(i2,i1)+Mlocal(2,1);
    M(i2,i2)=M(i2,i2)+Mlocal(2,2);
    
    %Contribución del elemento V(i)
    V(i1,i1)=V(i1,i1)+Vlocal(1,1);
    V(i1,i2)=V(i1,i2)+Vlocal(1,2); 
    V(i2,i1)=V(i2,i1)+Vlocal(2,1);
    V(i2,i2)=V(i2,i2)+Vlocal(2,2);
    
    %Contribución vector F(i)
    b(i1,1)= b(i1,1)+flocal(1);
    b(i2,1)= b(i2,1)+flocal(2);

end
\end{verbatim}
con esto obtenemos las matrices globales que definen el sistema global en la ecuación \eqref{EADM} y cuya solución aproxima la solución de la ecuación de advección-difusión. Aun falta discretizar el tiempo, observe que la derivada temporal aparece en el argumento de la matriz M, de la \textit{formulación matricial implícita} \eqref{EADM}, por lo que creamos una matriz nueva con el código
\begin{verbatim}
%Matriz discreta
B = (1/dt * M) + V + A ;
\end{verbatim}
y como tenemos condición de Dirichlet $u(0)=1$ y $u(1)=0$ la escribimos como
\begin{verbatim}
%Condición de Dirichlet
x_d=zeros(N+1,1);
x_d(1)=1;
x_d(N+1)=0  ;

b_d=b-B*x_d;
\end{verbatim}
donde la ultima linea impone la condición de Dirichlet. Ahora, como tenemos condición de Dirichlet solo necesitamos resolver para los grados de libertad interiores (véase la Observación \eqref{MD}), y por tanto solo necesitamos de la submatriz de grados de libertad interiores. Esto se extrae con 
\begin{verbatim}
%Matriz de Dirichlet
int=2:N;
B_int=B(int,int);
M_int=M(int,int);
u_0int=u_0(int);
b_int=b_d(int);
\end{verbatim}
con esto, solo creamos una iteración, que reemplaza $u_0$ por el ultimo resultado para lograr notar su desarrollo a través del tiempo, imprimiendo su resultado en una figura con las siguientes lineas
\begin{verbatim}
%Iteracion en el tiempo
for j=1:dt:T
x_sol=x_d;
b_ant=b_int + (1/dt * M_int * u_0int);
x_aux=B_int\b_ant;
u_0int=x_aux ;
x_sol(int)=x_aux;

%Plotar la solución
for k=1:N
    x=[E(k).xini, E(k).xfin];
   y=x_sol(E(k).dof);
  plot(x,y)
 hold on
end
pause(0.4)
end
\end{verbatim}

\begin{figure}[H]
  \centering
    \includegraphics[width=0.61\textwidth]{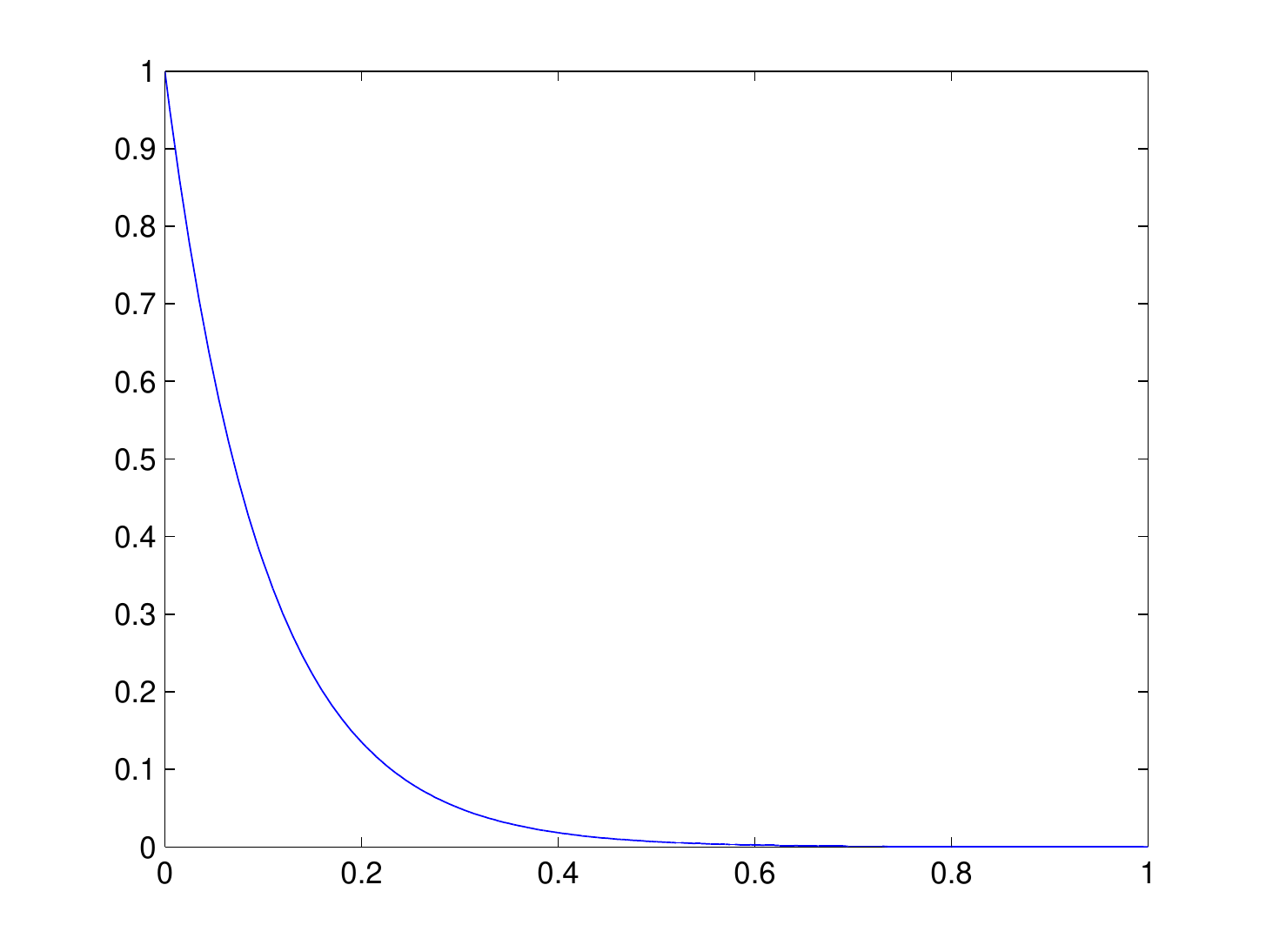}
  \caption{Aproximación de la ecuación de advección-difusión \eqref{EADF} en el instante  $t=\Delta t$, con incremento de tiempo $\Delta t=0.1$. Para esto, usamos la \textit{formulación de Galerkin} \eqref{EADGD} cuando $\mu=0.05,\beta=0.5$ y $f=0$.}
  \label{1D1}
 \end{figure}

\begin{figure}[H]
  \centering
    \includegraphics[width=0.61\textwidth]{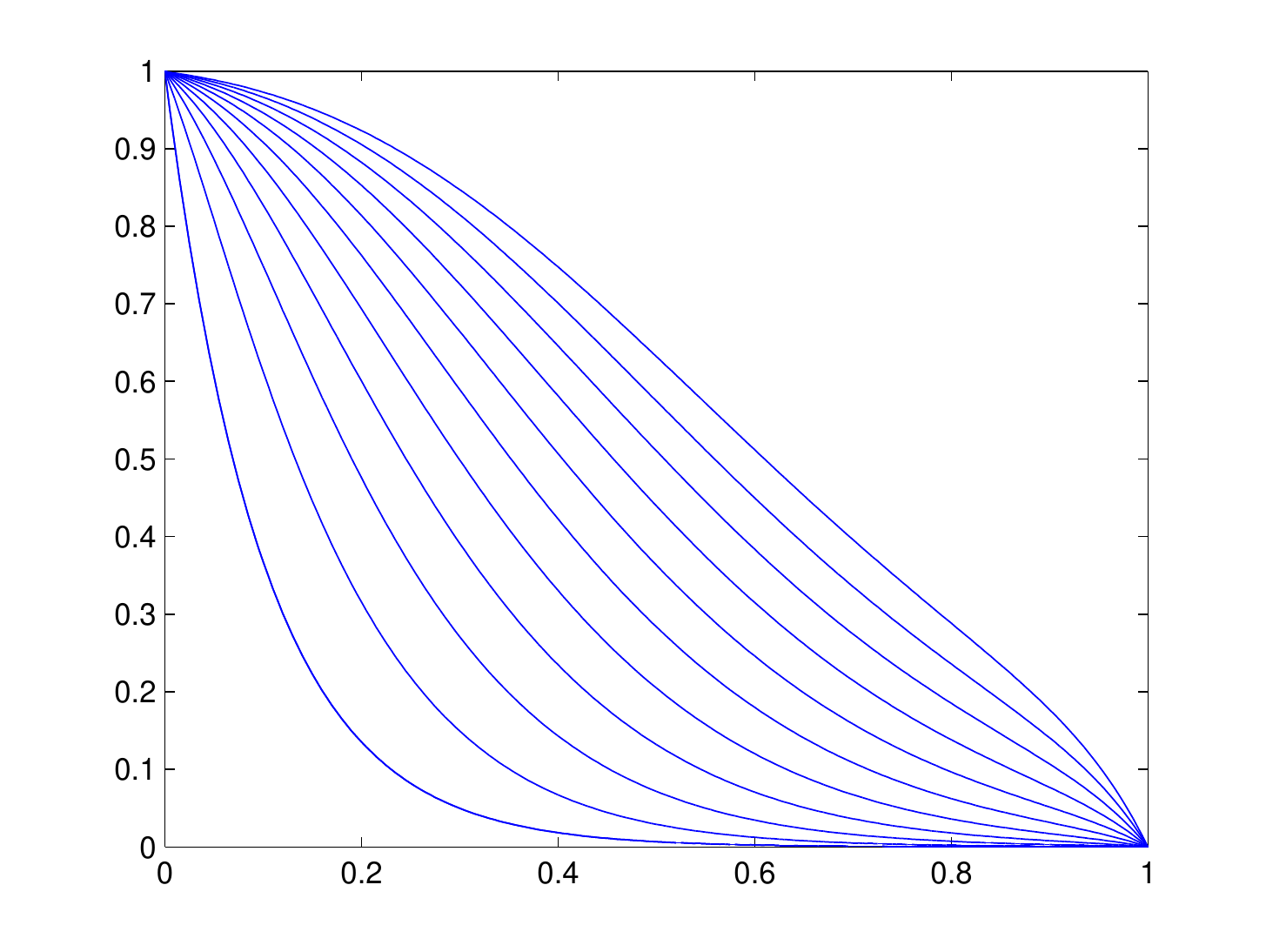}
  \caption{Aproximación de la ecuación de advección-difusión \eqref{EADF} en el instante $t=1.5$, con incremento de tiempo $\Delta t=0.1$. Para esto, usamos la \textit{formulación de Galerkin} \eqref{EADGD} cuando $\mu=0.05,\beta=0.5$ y $f=0$.}
  \label{1D2}
\end{figure}

\begin{figure}[H]
  \centering
    \includegraphics[width=0.61\textwidth]{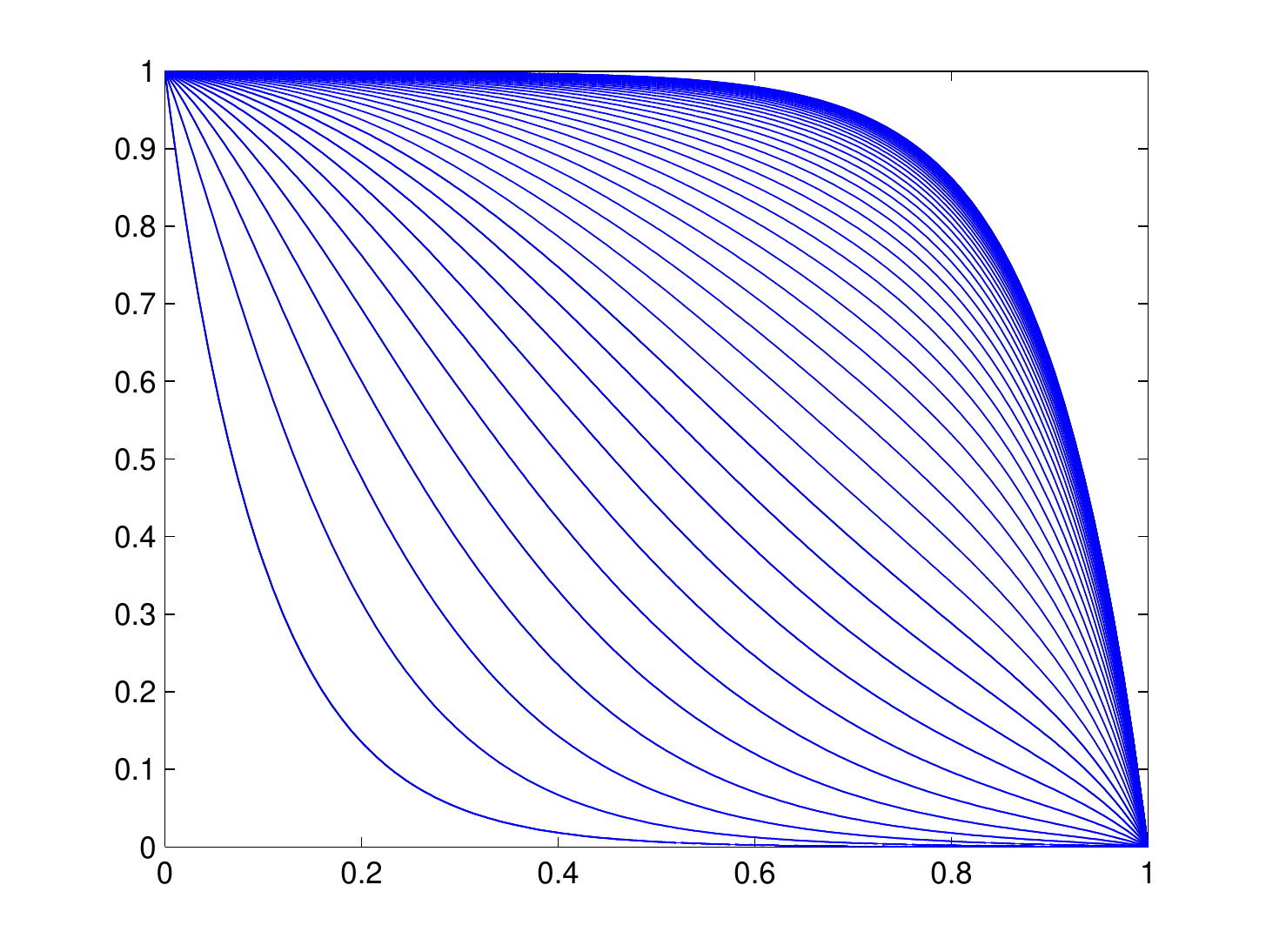}
  \caption{Aproximación de la ecuación de advección-difusión \eqref{EADF} en el tiempo final $t=5$, con incremento de tiempo $\Delta t=0.1$. Para esto, usamos la \textit{formulación de Galerkin} \eqref{EADGD} cuando $\mu=0.05,\beta=0.5$ y $f=0$.}
  \label{1D3}
\end{figure}

Con esto logramos la aproximación de la ecuación de advección-difusión por elementos finitos en una dimension, las figuras \ref{1D1}, \ref{1D2}, \ref{1D3} nos muestra la evolución de la aproximación de la solución en los tiempos $t=0.1$, $t=1.5$ y $t=5$ (tiempo final) respectivamente.

Con esto, tenemos una aproximación de la ecuación de advección-difusión, con coeficiente de difusión $\mu=0.05$ y velocidad de transporte $\beta=0.5$, con intervalos de tiempo $\Delta t=0.1$.
 
\section{Ejemplos en dos dimensiones usando FreeFem++}
\label{FFADE}

Con FreeFem++ también escribimos un código que nos permite visualizar la aproximación de elementos finitos de la ecuación de advección-difusión. 

\begin{figure}[hrt]
  \centering
    \includegraphics[width=0.7\textwidth]{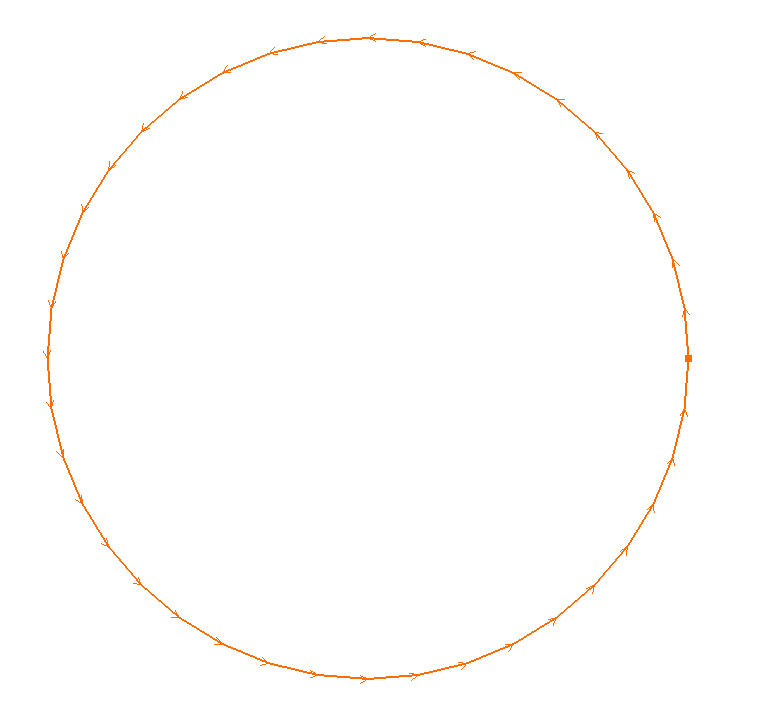}
  \caption{Dominio $\Omega$ definido en FreeFem++. Este dominio sera usado para resolver la ecuación de advección-difusión \eqref{EADF} usando elementos finitos $\mathbb{P}^2$ para aproximar la solución $u$.}
  \label{circ}
\end{figure}

Como sabemos, lo primero que tenemos que hacer en FreeFem++, es definir un dominio, que en este ejemplo sera un circulo, así que escribimos,
\begin{verbatim}

//Parametrización del dominio

border D(t=0,2*pi){x=cos(t);y=sin(t);label=1;};

//Plotar el dominio

plot(D(40),wait=1);
\end{verbatim}
El resultado de estas lineas se muestra en la figura \ref{circ}.

\begin{figure}[H]
  \centering
    \includegraphics[width=0.7\textwidth]{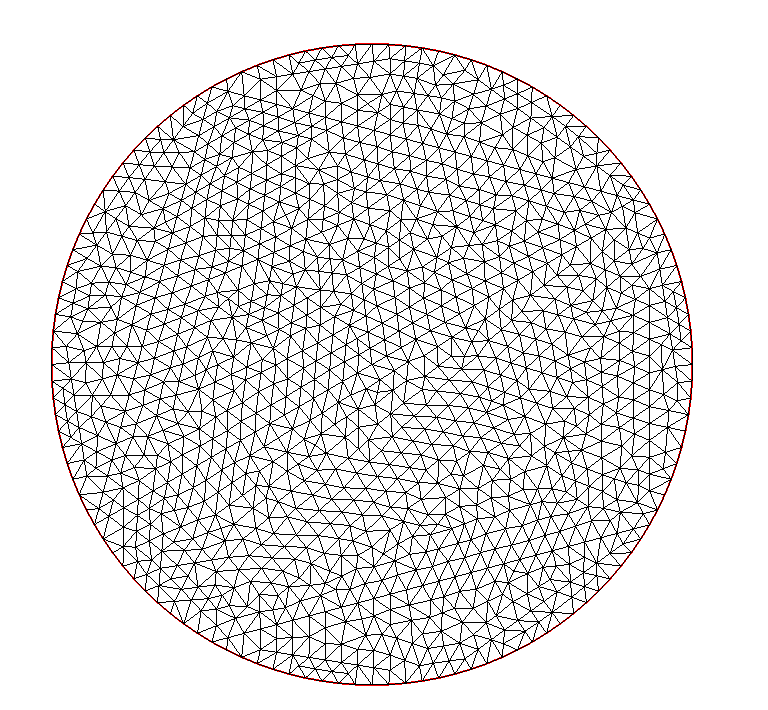}
  \caption{Triangulación $\mathcal{T}^h$ generada por FreeFem++, la cual contiene 2562 triángulos y 1342 vértices.}
  \label{Tri2}
\end{figure}

A continuación tomamos una triangulación en la que vamos a trabajar, para así crear nuestro espacio de elementos finitos y escribir el problema de la ecuación de advección-difusión, la triangulación se obtiene en FreeFem++ usando el siguiente comando,
\begin{verbatim}
// Crear la triangulación

mesh Th= buildmesh (D(120));
\end{verbatim}
y se puede visualizar usando el comando,
\begin{verbatim}
//Plotar la triangulación

plot(Th,wait=1);
\end{verbatim}
cuyo resultado se muestra en la Figura \ref{Tri2}. La triangulación generada tiene 2562 elementos(triángulos), y 1342 vértices.

\begin{figure}[hrt]
  \centering
    \includegraphics[width=0.7\textwidth]{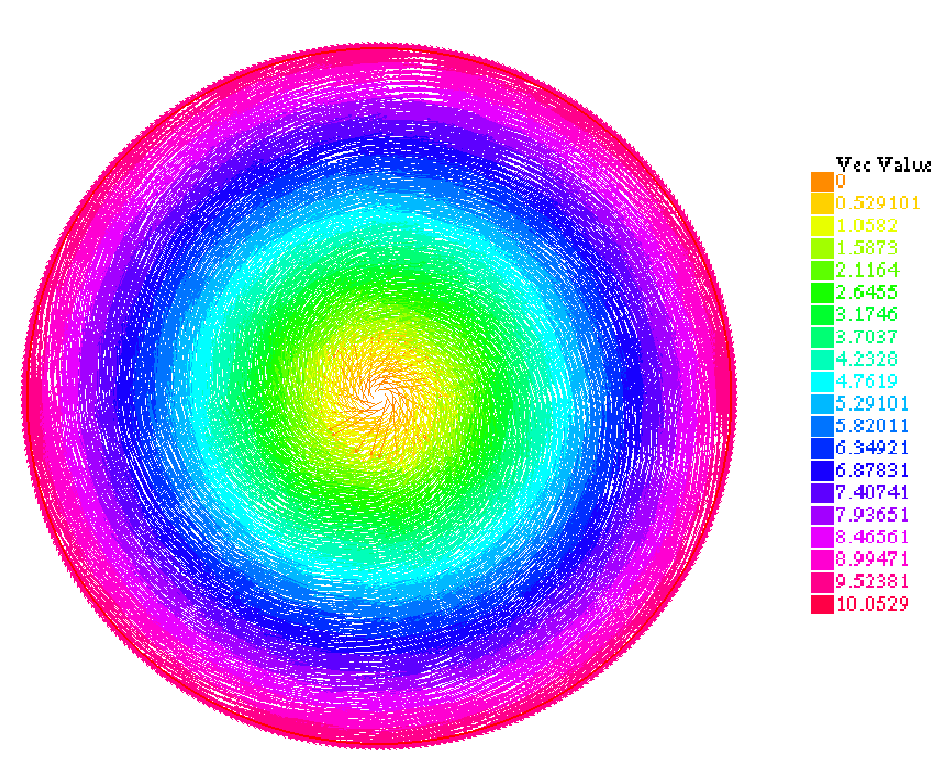}
  \caption{Campo de velocidades $\beta=(\beta_1,\beta_2)$ generado por FreeFem++, para aproximar el movimiento del fluido, es decir describe el proceso de advección para el fluido en el dominio $\Omega$.}
  \label{vel}
\end{figure}

Ahora, debemos definir el espacio en el que vamos a trabajar, que sera el espacio de funciones cuadráticas por partes $\mathbb{P}^2_0(\mathcal{T}^h)$. Luego, debemos definir todas las variables necesarias para formular el problema en este espacio, para esto escribimos las siguientes lineas,
\begin{verbatim}
fespace Vh(Th,P2); //Espacio de funciones P2

Vh u, v,uant,ud=0,betaxh,betayh, fh; //Declaración funciones de elementos finitos

func u0 =0; //Condición inicial
func mu=0.05; //Coeficiente de Viscosidad
func betax=10*y;  //Velocidad en x
func betay=-10*x; //Velocidad en y
func f=  (x>0.4)*(x<0.5)*(y>-0.1)*(y<0.1); //función de entrada de fluido

fh= f;	
betaxh=betax;
betayh=betay;

plot([betaxh,betayh], wait=1); //Plotar el campo de velocidades

real N=1000; //Numero de partes del tiempo
real  T=10; //Tiempo Final
dt=T/N; //Incremento de tiempo
uant=u0;
\end{verbatim}
donde \verb+u,v+ son las funciones de forma y de prueba respectivamente, uant se usa para iterar en el tiempo donde \verb+u_0=0, ud+ se usa para la condición de frontera, \verb+betaxh, betayh+ es el campo de velocidades $\beta=(\beta_1,\beta_2)$ para el coeficiente de transporte, \verb+fh+ es una función de entrada, \verb+mu+ ($\mu$) es el coeficiente de difusión, y \verb+dt+ es el intervalo de tiempo $\Delta t$. En la Figura \ref{vel} se muestra el campo de velocidades generado por las funciones \verb+betax, betay+.

\begin{figure}[H]
  \centering
  \subfloat[$t=0.1.$]{\label{AD1}
  \includegraphics[width=0.48\textwidth]{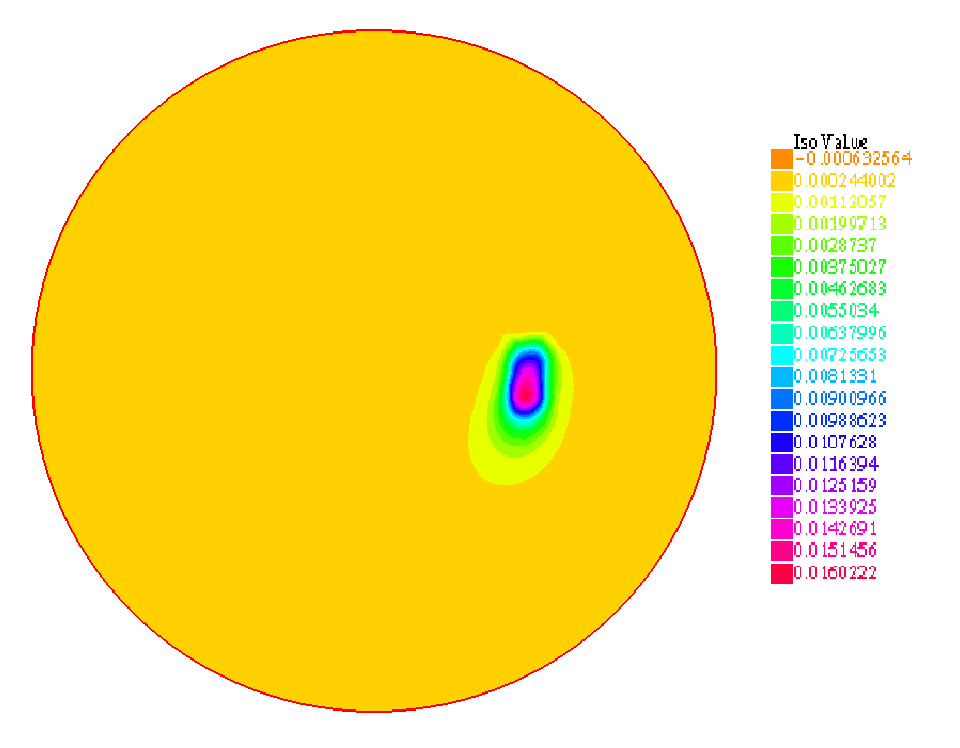}} 
  \subfloat[$t=1.5.$]{\label{AD2}
  \includegraphics[width=0.48\textwidth]{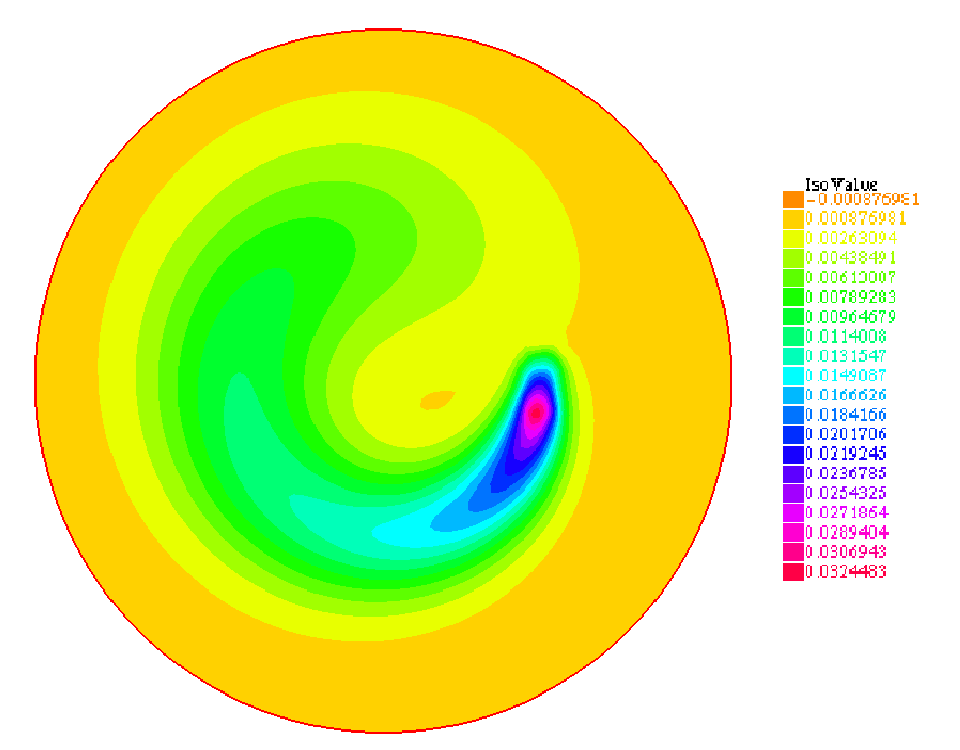}} \\
  \subfloat[$t=5.$]{\label{AD3}
  \includegraphics[width=0.48\textwidth]{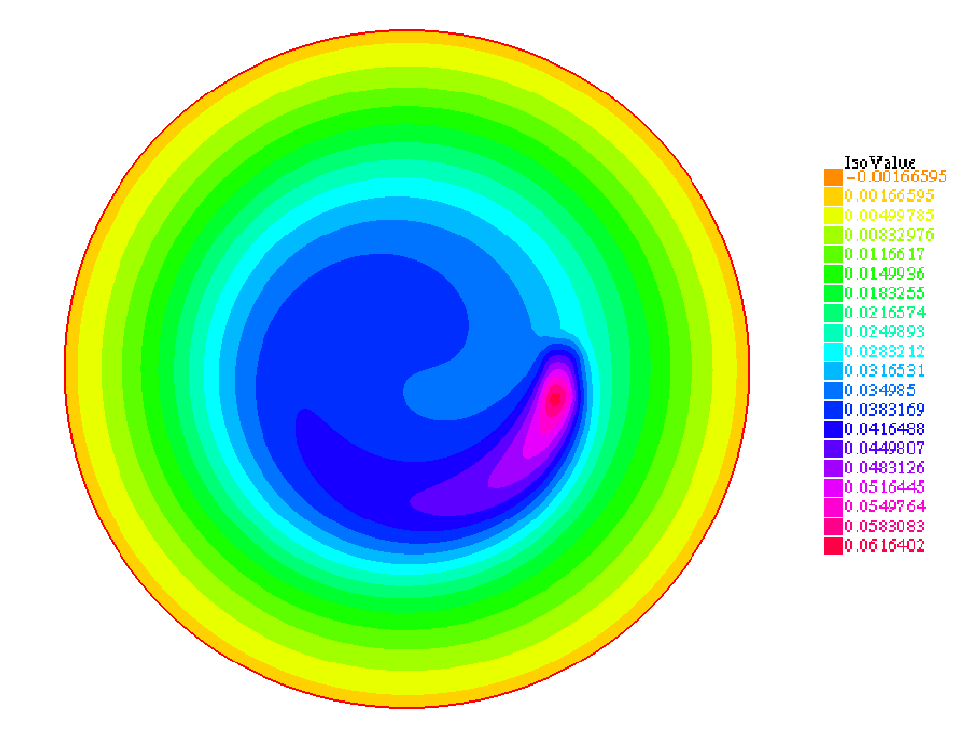}} 
  \subfloat[$t=10.$]{\label{AD4}
  \includegraphics[width=0.48\textwidth]{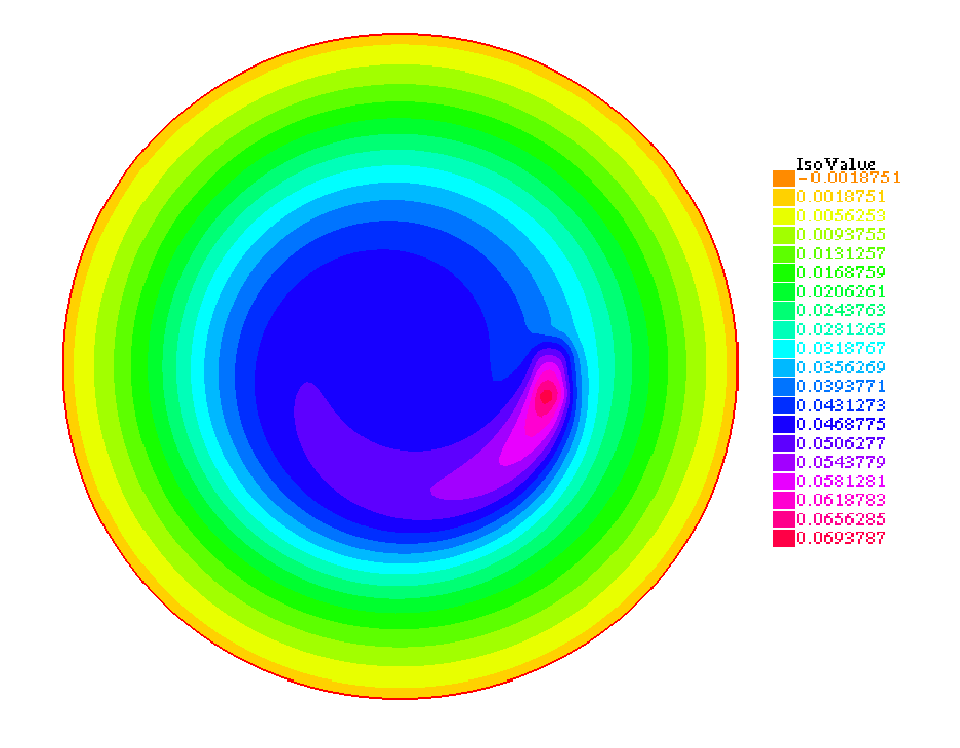}} \\
  \caption{Aproximación por MEF de la ecuación de advección-difusión \eqref{EADF}, en la triangulación $\mathcal{T}^h$ (véase Figura \ref{Tri2}), que contiene 2562 triángulos y 1342 vértices. Usando la \textit{formulación de Galerkin} \eqref{EADGD}, con el campo de velocidades $\beta=(\beta_1,\beta_2)$ de la Figura \ref{vel}, coeficiente de difusión $\mu=0.05$, para cada paso de tiempo $\Delta t=0.01$, en el instante $t$.}
\end{figure}

Ya con esto, escribimos el problema de la ecuación de advección-difusión en su formulación de Galerkin implícita como, 

\begin{verbatim}
//Formulación del problema de Advección-difusión

problem ecudifu(u,v)= int2d(Th)(mu*(dx(u)*dx(v) + dy(u)*dy(v)) + 
                               (betaxh*dx(u)+betayh*dy(u))*v + u*v/dt ) 
                      - int2d(Th)(uant*v/dt+fh*v)
                      + on(1,u=ud) ;
\end{verbatim}
donde \verb+problem ecudifu(u,v)+ crea un problema llamado \verb+ecudifu+ que contiene la  información que se dedujo en la Sección \ref{FDAD} y detallada en la ecuación \eqref{EADD'}, con condición de frontera Dirichlet cero, es decir que en la frontera no hay paso del fluido, la cual sera resuelta mas adelante. Esta es una de las ventajas de FreeFem++ que simula formulaciones de Galerkin que son idénticas a las formulaciones débiles. Ahora, falta resolver la ecuación de advección -difusión para cada paso de tiempo y visualizar la solución. Creamos una iteración de la siguiente manera,
\begin{verbatim}
//Iteracion en el tiempo dt

for(real tt=0;tt<T;tt+=dt){
    
    ecudifu; //Solución del problema ecudifu
    uant=u;

    plot(cmm="tiempo = "+tt,u,wait=0,dim=2,fill=1,value=0);
};
\end{verbatim}
esto va mostrando la aproximación de la ecuación de advección-difusión por elementos finitos a lo largo del tiempo. En FreeFem++, basta con escribir el nombre del problema \verb+ecudifu+ formulado anteriormente, para que resuelva el sistema lineal de la ecuación \eqref{EADM}, y el \verb+for+ hace que lo resuelva para cada paso de tiempo \verb+dt+, por medio de una factorización LU. En las figuras \ref{AD1}, \ref{AD2}, \ref{AD3}, \ref{AD4} se mostrara la aproximación por MEF de la ecuación de advección difusión, en los instantes $t=0.01$, $t=0.5$, $t=5$, y $t=10$ (tiempo final) respectivamente. 

Lo que nos muestra como se difundiría aproximadamente alguna sustancia, en un fluido que se este agitando circularmente, asumiendo que el coeficiente de difusión es $\mu=0.05$, el campo de velocidades es $\beta=(\beta_1,\beta_2)=$\verb+(betax,betay)+ de la Figura \ref{vel}, la sustancia es derramada en el punto $(x,y)=(0.3,0)$ desde el tiempo inicial $t=0$ hasta el tiempo final $t=10$.

\chapter{Ecuación de advección con velocidades de Stokes: ejemplos usando FreeFem++}

Para este capítulo usaremos la ecuación de Stokes \eqref{ESF}, la cual es
\begin{equation*}
\begin{aligned}[b]
&\textit{Encontrar } \vec{u}=(u_1,u_2) \textit{ y p tal que:} \\
&\left\{
\begin{aligned}
-\mu \Delta \vec{u} + \nabla p  = \vec{f}  \quad
 &\text{para } x \in \Omega,& \\
\text{(div }\vec{u})=0 \quad
 &\text{para } x \in \Omega,& \\
 \int_\Omega p = 0 \quad
 & \text{Condición de solubilidad} \\
 \vec{u}(x)=\vec{h}_D(x) \quad
 &\text{para }  x  \in \Gamma_D \subset \partial\Omega,	 \\
\nabla \vec{u}(x)\cdot\vec{\eta} = \vec{h}_N(x) \quad
 &\text{para } x \in \Gamma_N \subset \partial \Omega,
\end{aligned}\right.
\end{aligned}
\end{equation*}
que al aplicar el MEF, se obtiene la \textit{formulación de Galerkin} \eqref{ESG},
\begin{equation*}
\begin{aligned}[b]
&\textit{Encontrar } \vec{u} \in [\mathbb{P}^2(\mathcal{T}^h)]^2  \textit{ y }p \in \mathbb{P}^1(\mathcal{T}^h) \text{ tal que:} \\
&\left\{
\begin{aligned}
\mathcal{A}(u,v) + \mathcal{B}(v,p) = \mathcal{F}(v)  \quad
 &\text{para toda } \vec{v} \in \mathbb{P}^2_0(\mathcal{T}^h),& \\
\mathcal{B}(u,q)  = 0 \quad
 &\text{para toda } q \in \mathbb{P}^1_0(\mathcal{T}^h),& \\
\int_\Omega p = 0 \quad
 & \text{Condición de solubilidad} \\
 \vec{u}(x)=\vec{h}_D(x) \quad
 &\text{para }  x  \in \Gamma_D \subset \partial\Omega,	 \\
\nabla \vec{u}(x)\cdot\vec{\eta} = \vec{h}_N(x) \quad
 &\text{para } x \in \Gamma_N \subset \partial \Omega.
\end{aligned}\right.
\end{aligned}
\end{equation*} 
para obtener el campo de velocidades $\vec{u}=(u_1,u_2)=(\beta_1,\beta_2)=\vec{\beta}$ de la ecuación de advección-difusión \eqref{EADF},
\begin{equation*}
\begin{aligned}[b]
&\textit{Encontrar } u : \Omega \subset \mathbb{R}^2 \to \mathbb{R} \textit{ tal que:} \\
&\left\{
\begin{aligned}
\dfrac{\partial u}{\partial t}-\mu \Delta u + \beta_1 \dfrac{\partial u}{\partial x_1} + \beta_2 \dfrac{\partial u}{\partial x_2}  = f \quad
&\text{para } x \in \Omega,& \\
u(x,t)=g_D(x) \quad
 &\text{para }  x  \in \Gamma_D \subset \partial\Omega, \ t>0 \\
\nabla u(x,t)\eta(x) = g_N(x) \quad
 &\text{para } x \in \Gamma_N \subset \partial \Omega, \ t>0 \\
u(0,x)=u_0(x) \quad
& \text{para } x \in \Omega
\end{aligned}\right.
\end{aligned}
\end{equation*}
que tiene \textit{formulación de Galerkin} \eqref{EADGD} como sigue,
\begin{equation*}
\begin{aligned}[b]
&\textit{Dado } u_0 \in \mathbb{P}^2_0(\mathcal{T}^h) \textit{ encontrar } u_1,u_2,\ldots,u_m,\ldots  \textit{ tal que:} \\
&\left\{
\begin{aligned}
\frac{1}{\Delta t}{ \int_\Omega} u_mv =  \int_\Omega fv - \mu \int_\Omega \nabla u_{m-1} \nabla v- \int_\Omega\beta \nabla u_{m-1} v \quad
&\text{para toda } v \in \mathbb{P}^2_0(\mathcal{T}^h),& \\ 
  + \frac{1}{\Delta t}\int_\Omega u_{m-1}v \quad\quad \quad \quad \quad\quad\quad\quad\quad\quad\quad \\
u(x,t)=g_D(x) \quad
 &\text{para }  x  \in \Gamma_D \subset \partial\Omega, \ t>0 \\
\nabla u(x,t)\eta(x) = g_N(x) \quad
 &\text{para } x \in \Gamma_N \subset \partial \Omega, \ t>0 \\
u(0,x)=u_0(x) \quad
& \text{para } x \in \Omega
\end{aligned}\right.
\end{aligned}
\end{equation*}

Este capitulo esta dedicado a acoplar estas dos ecuaciones vistas en los anteriores capítulos. Específicamente, usamos el resultado de la aproximación de velocidades de la ecuación de Stokes \eqref{ESF}, como el campo de velocidades $\beta=(\beta_1,\beta_2)$ que genera la advección y se usa para calcular la aproximación de la ecuación de advección-difusión \eqref{EADF}.
 
Lo que haremos es mostrar un código en FreeFem++, que implementa este procedimiento. Volveremos al ejemplo del dique en la Sección \ref{FFSE}, para mostrar como la aproximación de elementos finitos de las velocidades de la ecuación de Stokes \eqref{ESF} en este dominio, se puede usar para calcular la aproximación de la ecuación de advección-difusión \eqref{EADF}. Primero, definimos de nuevo el mismo dominio, y creamos la misma triangulación, así
\begin{verbatim}
// Parametrización del dominio

border D(t=1,-1){x=-2*pi; y=t; label=2;};
border D1(t=-2*pi,2*pi){x=t; y=sin(t)-1; label=1;};
border D2(t=-1,1){x=2*pi; y=t; label=3;};
border D3(t=2*pi,-2*pi){x=t; y=sin(t)+1; label=1;};

plot(D(10)+D1(30)+D2(10)+D3(30));

//Crear  la triangulación

mesh Th= buildmesh (D(20)+D1(90)+D2(20)+D3(90));

plot(Th,wait=1,ps="Triangulacion.eps");
\end{verbatim}
las cuales generan las Figuras \ref{Dom} y \ref{Tri}. Ahora definiremos los espacios de funciones en los que vamos a trabajar, que en este caso serán los mismos que en la \textit{formulación de Galerkin} de la ecuación de Stokes \eqref{ESG}, $\mathbb{P}^2_0(\mathcal{T}^h)$ para las velocidades, y $\mathbb{P}^1_0(\mathcal{T}^h)$ para la presión $p$, lo que notaremos acá es que el espacio \verb+V_h+$=\mathbb{P}^2(\mathcal{T}^h)$ también  sera el espacio en el que se aplicara el MEF para la ecuación de advección-difusión, por lo que se definirán las variables de las dos ecuaciones en dicho espacio, así que escribiremos
\begin{verbatim}
// Definición de los espacios de elementos finitos

fespace Vh(Th,P2);
fespace Mh(Th,P1);

// Declaración de variables, funciones de forma y de prueba 

Vh u1,v1,u2, v2, betaxh, betayh, u, v,uant,g,ud=0;
Mh p,q;

func u0 =0; //Condición inicial
func mu=0.05; //Coeficiente de difusión
func mus=0.1; //Coeficiente de viscosidad de Stokes
func f=0; //Función externa
\end{verbatim}
donde \verb+u1,v1,u2,v2+ son las funciones de forma y de prueba de la ecuación de Stokes \eqref{ESG}, mientras que \verb+betaxh, betayh+ son los vectores del campo de velocidades $\beta=(\beta_1,\beta_2)$ de la ecuación de advección-difusión \eqref{EADF}, y \verb+f+ sera una fuerza externa; \verb+u,v+ son las funciones de forma y de prueba de la ecuación de advección difusión \eqref{EADGD}, \verb+uant+ se usa para la iteración en el tiempo, \verb+g+ es la función de entrada por la frontera \verb+D+, \verb+ud+ se usa para imponer la condición de frontera en \verb+D1,D3+ , y \verb+u_0=0+ es la condición inicial $u_0=0$; \verb+mu+ es el coeficiente de difusión de la ecuación de advección-difusión \eqref{EADF} y \verb+mus+ es el coeficiente de viscosidad de la ecuación de Stokes \eqref{ESF}. 

\begin{figure}[H]
  \centering
  \subfloat[t=0.1]{
  \includegraphics[width=0.65\textwidth]{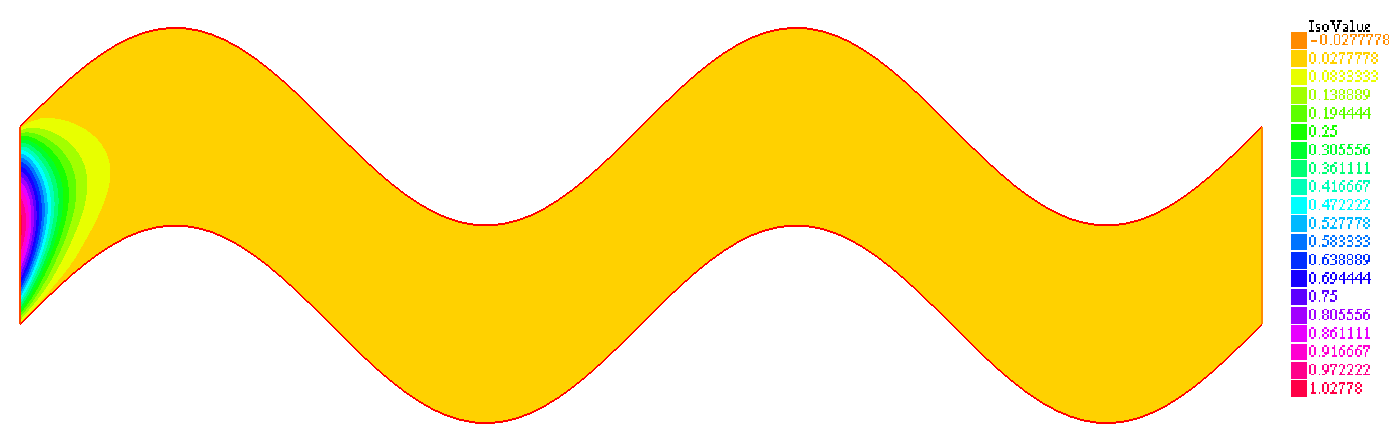}} \\
  \subfloat[t=1.5]{
  \includegraphics[width=0.65\textwidth]{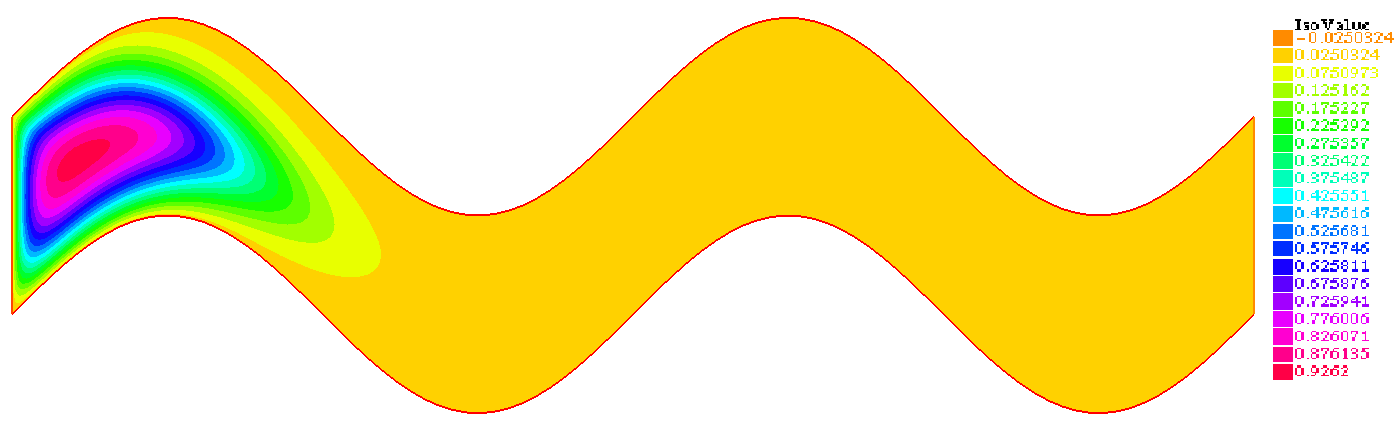}} \\
  \subfloat[t=5]{
  \includegraphics[width=0.65\textwidth]{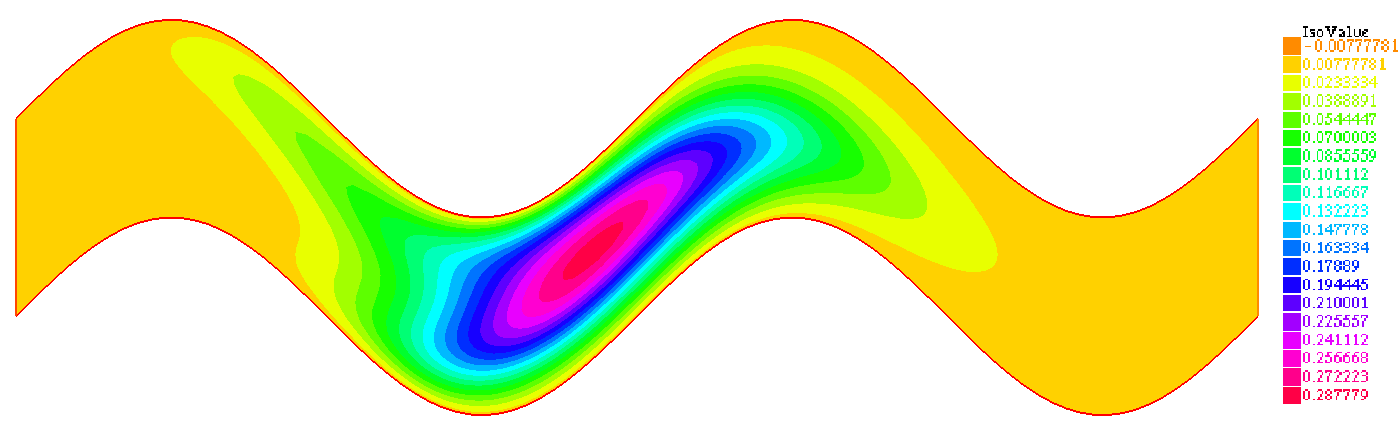}} \\
  \subfloat[t=10]{
  \includegraphics[width=0.68\textwidth]{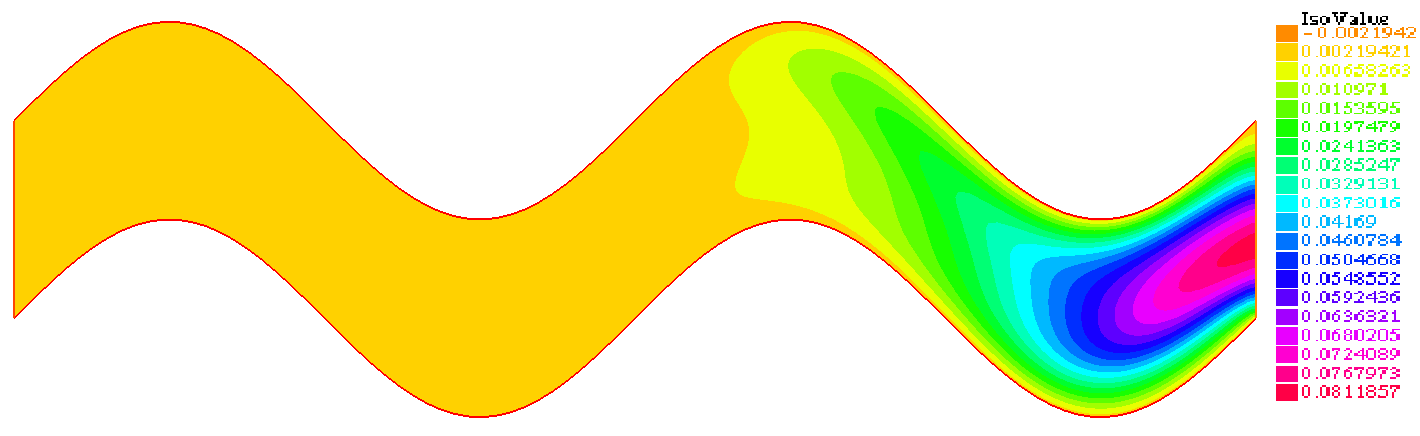}} \\
  \caption{Aproximación por MEF de la ecuación de advección-difusión \eqref{EADF}, con coeficiente de difusión $\mu=0.05$, en la triangulación $\mathcal{T}^h$ (véase Figura \ref{Tri}), que contiene 2176 triángulos y 1199 vértices. Usando la \textit{formulación de Galerkin} \eqref{EADGD}, con el campo de velocidades,  $\beta=(\beta_1,\beta_2)$ de la Figura \ref{Stokes}, obtenido al solucionar la ecuación de Stokes \eqref{ESF}, con coeficiente de viscosidad $\mu_{Stokes}=0.1$, para cada paso de tiempo $\Delta t=0.01$, en el instante $t$.}
  \label{SAD}
\end{figure}

Con esto, lo que sigue es resolver la ecuación de Stokes, con las mismas condiciones de frontera que en la Sección \ref{FFSE}, usando el comando \verb+solve+ que soluciona el problema escrito, que acá es la \textit{formulación débil} de la ecuación de Stokes \eqref{ESD'}, y escogemos cual sera el método para resolver el sistema, en esta caso \verb+solver=Crout+.
\begin{verbatim}
//Resolver el problema de Stokes

solve Stokes (u1,u2,p,v1,v2,q,solver=Crout) =
      int2d(Th) (mus*(dx(u1)*dx(v1)+dy(u1)*dy(v1)+dx(u2)*dx(v2)+dy(u2)*dy(v2)) 
	              - p*q*(0.00000001) - p*dx(v1) - p*dy(v2) - dx(u1)*q -dy(u2)*q )
				  +on(1,u1=0,u2=0)+on(2,u1=-1.5*(y-1)*(y+1),u2=1);
\end{verbatim}
Ahora, asignaremos la solución obtenida \verb+u_1,u_2+ a $\beta=(\beta_1,\beta_2),$ es decir,la solución $(u_1,u_2)$ sera el campo de velocidades de la ecuación de advección-difusión, que en FreeFem++ se escribe,
\begin{verbatim}
betaxh=u1 ; //Asignar u_1 al vector velocidad betaxh
betayh=u2 ; //Asignar u_2 al vector velocidad betayh

//Plotar el campo de velocidades
			
plot([betaxh,betayh], wait=1);
\end{verbatim} 
la cual obviamente, dará el mismo campo de velocidades visto en la Figura \ref{Stokes}, ya con esto basta introducir el código  visto en la Sección \ref{FFADE}, es decir crear $\Delta t=0.01$, y escribir la formulación de Galerkin implícita como el problema de advección difusión \eqref{EADGD}. De aquí,
\begin{verbatim}
real N=1000; //Numero de partes del tiempo
real  T=10;  //Tiempo final
dt=T/N;      //Incremento del tiempo
uant=u0;     

//Problema de advección-difusión

problem ecudifu(u,v)= int2d(Th)(mu*(dx(u)*dx(v) + dy(u)*dy(v)) 
                               + (betaxh*dx(u)+betayh*dy(u))*v+u*v/dt) 
                      - int2d(Th)(uant*v/dt+f*v)
                      + on(1,u=ud) + on(2,u=g);
\end{verbatim}
FreeFem++ creara la aproximación de elementos finitos de la ecuación de advección difusión con condición de frontera Dirichlet cero en \verb+D1+ y \verb+D3+, Neumann cero en \verb+D2+, y en \verb+D+ una condición variable en el tiempo. Por ultimo, falta resolver en cada instante de tiempo y visualizar la solución. En este ejemplo consideraremos que la fuente de la sustancia esta activa desde el tiempo $t=0$ hasta el tiempo $t=1.5$, luego de esto acabara la entrada de sustancia, de modo que
\begin{verbatim}
for(real tt=0;tt<T;tt+=dt){
    g=-(y-1)*(y+1) ;
	  if (tt>1.5)
      g=0;
     ecudifu;
    
    uant=u;
    
    plot(cmm="tiempo = "+tt,u,wait=0,dim=2,fill=1,value=0);
};
\end{verbatim}
Con esto, obtenemos la aproximación de elementos finitos de la ecuación de  advección-difusión, utilizando el campo de velocidades de Stokes, la cual se vera a través del tiempo algo como en la Figura \ref{SAD}.

Con esto, mostramos la aproximación de elementos finitos de la ecuación de advección-difusión, con coeficiente de difusión $\mu = 0.05$ y con el campo de velocidades de la ecuación de Stokes, con coeficiente de viscosidad $\mu_{stokes}=0.1$.

\end {document}